\documentclass[preprint,12pt]{elsarticle}
\usepackage{framed,multirow}
\usepackage{amssymb}
\usepackage{latexsym}
\usepackage{url}
\usepackage{xcolor}
\usepackage{float,epsfig,epstopdf,bm,floatflt,here,relsize}
\usepackage{rotating}
\usepackage{amssymb}
\usepackage{amsmath}
\usepackage{mathtools}
\usepackage{tikz}
\usepackage{tkz-euclide}
\usepackage{multirow}
\usepackage{booktabs}
\usepackage{array}
\usepackage{graphicx}\newcolumntype{H}{>{\setbox0=\hbox\bgroup}c<{\egroup}@{}}
\usepackage{arydshln}
\usepackage{caption}
\usepackage{subcaption}
\usepackage{cleveref}
\usepackage{empheq}
\usepackage[left=2cm,right=2cm]{geometry}
\allowdisplaybreaks

\definecolor{newcolor}{rgb}{.8,.349,.1}
\journal{}

\begin{document}
%\verso{Sen and Sheu}
\begin{frontmatter}



\title{Fourth order compact scheme for the Navier-Stokes equations on time deformable domains}


\author[a1]{Shuvam Sen\corref{cor1}}
\ead{shuvam@tezu.ernet.in}
\author[a2,a3,a4]{Tony W. H. Sheu\corref{cor1}}
\ead{twhsheu@ntu.edu.tw}
\address[a1]{Department of Mathematical Sciences, Tezpur University,\\ Tezpur 784028, India}
\cortext[cor1]{Corresponding author}
\address[a2]{Department of Engineering Science \& Ocean Engineering, National Taiwan University,\\ Taipei 10617, Taiwan}
\address[a3]{Center for Advanced Study in Theoretical Science, National Taiwan University,\\ Taipei 10617, Taiwan}
\address[a4]{Institute of Mathematical Sciences, National Taiwan University,\\ Taipei 10617, Taiwan}

\begin{abstract}
In this work, we report the development of a spatially fourth order temporally second order compact scheme for incompressible Navier-Stokes (N-S) equations in time-varying domain. Sen [J. Comput. Phys. 251 (2013) 251-271] put forward an implicit compact finite difference scheme for the unsteady convection-diffusion equation. It is now further extended to simulate fluid flow problems on deformable surfaces using curvilinear moving grids. The formulation is conceptualized in conjunction with recent advances in numerical grid deformations techniques such as inverse distance weighting (IDW) interpolation and its hybrid implementation. Adequate emphasis is provided to approximate grid metrics up to the desired level of accuracy and freestream preserving property has been numerically examined. As we discretize the non-conservative form of the N-S equation, the importance of accurate satisfaction of geometric conservation law (GCL) is investigated. To the best of our knowledge, this is the first higher order compact method that can directly tackle non-conservative form of N-S equation in single and multi-block time dependent complex regions. Several numerical verification and validation studies are carried out to illustrate the flexibility of the approach to handle high-order approximations on evolving geometries. 

\end{abstract}

\begin{keyword}
%\KWD
Higher order compact\sep Incompressible Navier-Stokes equations\sep Geometric conservation\sep Moving grids \sep Deformable domain



\end{keyword}

\end{frontmatter}


\section{Introduction}\label{}
Many problems in science and engineering involve solving partial differential equations (PDEs) in complex geometries that might even evolve in time. Here, it is necessary to properly account for the time variation of the solution domain and could require a new discretization at each time step. In this context improving the accuracy of a solution with simultaneous reduction of computational costs has been an ever-lasting goal in the computational fluid dynamics (CFD) community. Higher order methods enjoy unique interest largely because of their ability to produce highly accurate solutions often at a comparably lower computational cost. A detailed discussion of the advantages and challenges associated with higher order schemes can be found in \cite{wan_fid_abg_13}. Nevertheless, for efficient capture of small amplitude high wave-number fluctuations numerical schemes with higher truncation order and lower dissipation error are preferred. Fourth order compact schemes are thus inherently advantageous to correctly predict high frequency waves with a fewer number of grid points. These schemes carry lesser dispersion and anisotropy errors. Moreover, in compact schemes, the computational cost of inverting matrices is low as their smaller stencils translate to a banded system with smaller bandwidth. Typically these schemes are used to examine flow fields on simple static domains discretized by Cartesian meshes. Pioneering investigation using high order compact method on a time-varying domain was carried out by Visbal and Gaitonde \cite{vis_gai_02}. Here mapping is constructed between a fixed reference domain and the real time-varying geometry. Subsequently, the original conservation law is first transformed to the reference configuration and then solved using a high order compact difference scheme in time dependent curvilinear coordinates. Apart from exploring the virtues of compact discretization the process has several advantages. Steger \cite{ste_78} whose early work pertains to the second order accurate discretization of conservative N-S equations listed some of those leverages. These include the boundary surfaces in the physical plane being mapped onto rectangular surfaces in the transformed plane, incorporation of unsteady motion into the governing equation, and clustering of grid points in regions with high flow gradients. Visbal and Gaitonde \cite{vis_gai_02} worked with the conservative form of the N-S equations and following them many important works have discretized conservative forms \cite{gha_far_12, lia_miy_zha_14, li_hen_ban_16, cox_lia_ple_16} to simulate flow using moving and deforming grids. In this approach, efficient discretization strategies are designed to approximate first order derivatives. The second order viscous terms are often tackled by the repeated application of the approximation procedure developed for the first order derivatives or at times by central differencing.

Working with a structured grid, finite volume approximation of first order advective terms in conservative formulation for variable meshes generally requires due diligence. The extra care required for efficient implementation and to avoid loss of accuracy is well pointed out by Hirt and Nichols \cite{hir_nic_81} and Li \emph{et al.} \cite{li_hes_zie_05}. For general grid motion, two geometrical relations called GCL ensure that each cell in the flow region is closed and safeguard conservation of volume for time-varying cells. As GCL represents geometric identities, its differential and integral forms are trivial, but unfortunately, in discretized form, automatic satisfaction of the GCL is not guaranteed. A detailed argument on conservation laws and their finite volume and finite difference formulation with emphasis on the specialized procedure adopted for grid metric calculation can be found in the work of Sj\"{o}green \emph{et al.} \cite{sjo_yee_vin_14} and Abe \emph{et al.} \cite{abe_hag_non_16}. Nevertheless, a brief discussion on the development of GCL and its limitation in ensuring only up to second or third order accuracy for space and time using these metric formulations was discussed by Persson \emph{et al.} \cite{per_bon_per_09}. 

Although there is a growing interest in higher order compact finite difference methods for fluid problems another challenge hindering its growth in time deforming domain is the absence of an effective grid deformation algorithm. It is well acknowledged that space-time accurate solutions of the N-S equations require good quality structured grid with smooth grid metrics as the solution domain continues to evolve. Maintaining smoothness and appropriate spacing, even when the initial mesh is comprehensive, are of primary importance in the numerical simulation of evolving physical phenomena. Failure to resolve important length scales of the flow at each time step might contribute to aliasing error. Traditional mesh deformation techniques based on transfinite interpolation (TFI) \cite{wan_prz_94}, although quite fast, are not always sufficient to obtain a reasonable quality of the deformed mesh. The generality of the transformed equations could be amply realized only with a fairly automatic method of generating smoothly varying grids that allow clustering and fit into arbitrary domains. In this context, contemporary grid deformations techniques using radial basis function (RBF) and inverse distance weighting (IDW) could play a pivotal role. Modern mesh deformation techniques based on RBFs \cite{boe_sch_bij_07, ren_all_08} produce high quality mesh but are computationally expensive. Another recently introduced explicit grid movement strategy established by Witteveen and Bijl \cite{wit_bij_09} that utilizes weighted sums of the reciprocal distance of the boundary node displacements to deform interior nodes was found to be quite effective and fast. This method is popularly known as the IDW method. Subsequently, this procedure was significantly improved by Luke \emph{et al.} \cite{luk_col_bla_12}. Towards this end, a comprehensive study involving modern grid deformation techniques was carried out by Sen \emph{et al.} \cite{sen_nay_bre_17}. In that study, the authors were able to bring an improvement over the current state-of-the-art by combining approaches of fundamental mesh deformation methods. It was shown that hybrid methodology can deform a grid with high skew quality metric value in comparatively lesser CPU time. The procedure was indeed found to be effective in simulating fluid--structure interaction (FSI) problems relying on large-eddy simulation (LES) \cite{apo_nay_ble_19, nay_bre_woo_20}. 

The purpose of this work is to present a fourth order compact finite difference method for simulating viscous incompressible fluids in general deformable boundaries. The formulation discretizes the non-conservative form of the incompressible N-S equations and is different from the strategy adopted for the conservative form \cite{vis_gai_02}. The fourth order method advocated here is a pinnacle of a series of independent yet intertwined studies starting with Sen \cite{sen_13} and going onto Sen \cite{sen_16}, Sen and Sheu \cite{sen_she_17} and finally Sen \emph{et al.} \cite{sen_nay_bre_17}. It is a realization of compact discretization strategy for parabolic convection-diffusion equation in time-varying boundaries. The schemes, where transport variable and its first derivatives are considered as the unknowns, coalesce virtues of compact discretization and Pad\'{e} approximation and carry good numerical dispersion and dissipation characteristics. Working with the current methodology it seems that satisfaction of GCL might not be necessary and does not play a significant role. As the grid evolves, the use of modern hybrid grid deformation techniques \cite{luk_col_bla_12, sen_nay_bre_17} in the absence of an analytic process helps maintain excellent mesh quality. Temporally second order accurate, the scheme is probably the first higher order compact discretization of non-conservative N-S equations in deformable geometries as we failed to locate any advancement towards such an approximation in the literature. Without any prejudice, we shall like to mention the work of Chen and Xie \cite{che_xie_16} where a second order accurate discretization based on the vorticity-vector potential formulation was discussed for incompressible flows in time-dependent curvilinear coordinates.

In the next section, we present the underlying mathematical background associated with various stages necessary for the development of the scheme. Subsequently in section 3, numerical verification and validation studies are carried out. Finally in section 4 conclusions are jotted down.

\section{Mathematical formulation}
\subsection{Governing equation}
We consider fluid to occupy a planner domain $\bm x\in\Omega(t)$ at time $t\in(0,T]$. The motion is assumed to be governed by the unsteady incompressible N-S equations. In the Eulerian frame using Cartesian coordinate the nondimensional primitive variable form of this equation together with the equation of continuity is
\begin{subequations}\label{N_S_pri}
\begin{empheq}[left=\empheqlbrace]{align}
&\cfrac{\partial \bm u}{\partial t}+\left( \bm u \cdot \nabla \right) \bm u-\cfrac{1}{Re}\Delta\bm u = -\nabla p,\;\;\;\bm x\in\Omega(t) \label{N_S_pri1}\\
&\nabla \cdot \bm u =0,\;\;\;\;\;\;\;\;\;\;\;\;\;\;\;\;\;\;\;\;\;\;\;\;\;\;\;\;\;\;\;\;\;\;\;\;\;\;\bm x\in\Omega(t) \label{N_S_pri2}
\end{empheq}
\end{subequations}
where $\bm u(\bm x,t)$ denotes the velocity and $p(\bm x,t)$ is the pressure. This flow is characterized by the dimensionless Reynolds number $Re$. Apart from initial conditions, the above equations require appropriate boundary conditions on $\partial \Omega(t)$. The elliptic pressure Poisson equation 
\begin{align}\label{pre_poi}
-\Delta p&=\nabla\cdot((\bm u\cdot\nabla)\bm u)\nonumber\\
		 &=\nabla\bm u:(\nabla\bm u)^T,\;\;\;\bm x\in\Omega(t)
\end{align}
which requires an additional boundary condition 
\begin{equation}\label{pre_poi_bou}
\bm n\cdot\nabla p=\bm n\cdot\left(\frac{1}{Re}\Delta\bm u-(\bm u\cdot\nabla)\bm u\right),\;\;\;\bm x\in\partial\Omega(t)
\end{equation}
is an able substitute of Eq. (\ref{N_S_pri2}). An alternative formulation, economical in two-dimension (2D) includes streamfunction $(\psi)$ and vorticity $(\omega)$ and is often referred to as $\psi-\omega$ formulation. It is 
\begin{subequations}\label{N_S_psi}
\begin{empheq}[left=\empheqlbrace]{align}
&\cfrac{\partial \omega}{\partial t}+(\bm u \cdot \nabla) \omega
                                -\cfrac{1}{Re} \Delta \omega = 0,\;\;\;\bm x\in\Omega(t) \label{N_S_psi1} \\
&-\Delta \psi = \omega,\;\;\;\;\;\;\;\;\;\;\;\;\;\;\;\;\;\;\;\;\;\;\;\;\;\;\;\;\;\;\bm x\in\Omega(t) \label{N_S_psi2}
\end{empheq}
\end{subequations}
where $\omega=(\nabla \times \bm u)\cdot \hat{k}$ is the out-of-plane component of vorticity and $\bm u=\nabla \times \Psi$ with $\Psi=(0, 0, \psi)$. A prototype equation that can be used to describe both primitive and $\psi-\omega$ forms of N-S systems is the convection-diffusion equation
\begin{equation}\label{con_dif}
\frac{\partial \phi}{\partial t}-a(\bm x,t)\Delta \phi+\bm c(\bm x, t)\cdot\nabla\phi=s(\bm x, t),\;\;\;(\bm x, t)\in \Omega(t)\times (0, T]
\end{equation}
with the initial condition
\begin{equation}\label{con_dif_i}
\phi(\bm x, 0)=\phi_0(\bm x),\;\;\;\;\;\bm x\in \Omega(0)
\end{equation}
and boundary condition
\begin{equation}\label{con_dif_b}
b_1(\bm x, t)\phi+b_2(\bm x, t)\frac{\partial \phi}{\partial \bm n}=g(\bm x, t),\;\;\;\;\;\bm x\in\partial \Omega(t),\;\;t\in(0,T].
\end{equation}
Here $\phi(\bm x, t)$ is the transport variable with convection velocity $c(\bm x, t)=(c_1(\bm x, t), c_2(\bm x, t))^T$, and $a(\bm x, t)>0$ having unit boundary normal vector $\bm n$.
\subsection{Time-dependent coordinate transformation}
We introduce an iso-parametric invertible time dependent mapping 
\begin{equation}\label{transformation}
x=x(\xi,\eta,\tau),\;\;\;y=y(\xi,\eta,\tau),\;\;\;t=\tau
\end{equation}
to transform the spatial domain $\Omega(t)$ to a topologically equivalent domain $\mathcal{D}$ thereby resulting in the generation of curvilinear meshes. The Jacobian of the transformation is 
\begin{equation}\label{transformation_jacobian}
\mathcal{J}=\frac{\partial (x,y,t)}{\partial (\xi,\eta,\tau)}=\left[\begin{array}{c c c}
x_{\xi} & x_{\eta}  & x_{\tau}\\
y_{\xi} & y_{\eta}  & y_{\tau}\\
0 		& 0			& 1
\end{array} \right].
\end{equation}
The Jacobian of the inverse transform then corresponds to
\begin{equation}\label{transformation_inv}
\mathcal{J}^{-1}=\frac{\partial (\xi,\eta,\tau)}{\partial (x,y,t)}=\left[\begin{array}{c c c}
\xi_x 	& \xi_y  & \xi_t\\
\eta_x 	& \eta_y & \eta_t\\
0 		& 0			& 1
\end{array} \right]=\frac{1}{J}\left[\begin{array}{c c c}
y_{\eta} & -x_{\eta}& -J_1\\
-y_{\xi} & x_{\xi}  & -J_2\\
0 		 & 0		& J
\end{array} \right]
\end{equation}
where $\displaystyle J=\left|\frac{\partial(x,y)}{\partial(\xi,\eta)}\right|\neq0$ and $\displaystyle J_1=\left|\frac{\partial(x,y)}{\partial(\tau,\eta)}\right|$, $\displaystyle J_2=\left|\frac{\partial(x,y)}{\partial(\xi,\tau)}\right|$. Note that $J_1$ and $J_2$ are influenced by grid velocity in addition to grid metrics. Directly using the chain rule the convection-diffusion equation (\ref{con_dif}) along with auxiliary conditions (\ref{con_dif_i}) and (\ref{con_dif_b}) can be rewritten as
\begin{subequations}\label{con_dif_tra}
\begin{empheq}[left=\empheqlbrace]{align}
&\partial_{\tau}\phi(\xi,\eta,\tau)+A\phi(\xi,\eta,\tau)=s(\xi,\eta,\tau),\;\;\;\;\;\;(\xi,\eta,\tau)\in\mathcal{D}\times(0,T] \label{con_dif_tra1}\\
&\phi(\xi,\eta,0)=\phi_0(\xi,\eta),\;\;\;\;\;\;\;\;\;\;\;\;\;\;\;\;\;\;\;\;\;\;\;\;\;\;\;\;\;\;\;\;(\xi,\eta)\in\mathcal{D}
\label{con_dif_tra2}\\
&\vartheta_1(\xi,\eta,\tau)\phi+\vartheta_2(\xi,\eta,\tau)\partial_{\bm \nu}\phi=g(\xi,\eta,\tau),\;(\xi,\eta)\in\partial\mathcal{D},\;\tau\in(0,T].\label{con_dif_tra3}
\end{empheq}
\end{subequations}
with partial differential operator $[A]$ being
\begin{eqnarray}\label{operator}
[A]\equiv[-\alpha_1\partial_{\xi\xi}-\beta\partial_{\xi\eta}-\alpha_2\partial_{\eta\eta}+\chi_1\partial_{\xi}+\chi_2\partial_{\eta}].
\end{eqnarray}
Eq. (\ref{con_dif_tra1}) possesses positive definiteness of the diffusion matrix 
\textit{i.e.} $\alpha_1>0$, $\alpha_2>0$, $|\beta|^2\le4\alpha_1\alpha_2$ $\forall$ $(\xi,\eta,\tau)\in\mathcal{D}\times(0,T]$ where
\begin{subequations}\label{tran_coeff}
\begin{empheq}[left=\empheqlbrace]{align}
\alpha_1&=\frac{a}{J^2}(x^2_{\eta}+y^2_{\eta}),\\
\alpha_2&=\frac{a}{J^2}(x^2_{\xi}+y^2_{\xi}),\\
\beta&=-\frac{2a}{J^2}(x_{\xi}x_{\eta}+y_{\xi}x_{\eta}),\\
\chi_1&=\frac{1}{J}(-J_1+c_1y_{\eta}-c_2x_{\eta})-\frac{a}{J^3}\bigg(J_{\eta}(x_{\xi}x_{\eta}+y_{\xi}y_{\eta})\nonumber\\
&-J_{\xi}(x^2_{\eta}+y^2_{\eta})+J(x_{\eta}x_{\xi\eta}+y_{\eta}y_{\xi\eta}-x_{\xi}x_{\eta\eta}-y_{\xi}y_{\eta\eta}\bigg),\\
\chi_2&=\frac{1}{J}(-J_2-c_1y_{\xi}+c_2x_{\xi})-\frac{a}{J^3}\bigg(J_{\xi}(x_{\xi}x_{\eta}+y_{\xi}y_{\eta})\nonumber\\
&-J_{\eta}(x^2_{\xi}+y^2_{\xi})+J(x_{\xi}x_{\xi\eta}+y_{\xi}y_{\xi\eta}-x_{\eta}x_{\xi\xi}-y_{\eta}y_{\xi\xi}\bigg).
\end{empheq}
\end{subequations}

\subsection{Geometric conservation law}
As discretization of the governing equation is carried out using a temporally varying grid the mesh geometry changes with time. To better focus the role of GCL, convection and diffusion terms of the 2D convection-diffusion equation with its coefficients kept fixed are expressed as inviscid and viscous fluxes, respectively. The following equation
\begin{eqnarray}\label{conservativeI}
\frac{\partial\phi}{\partial t}-\frac{\partial}{\partial x}\left(a\frac{\partial \phi}{\partial x} \right)-\frac{\partial}{\partial y}\left(a\frac{\partial \phi}{\partial y} \right)+\frac{\partial }{\partial x}\left(c_1 \phi \right)+\frac{\partial }{\partial y}\left(c_2 \phi \right)=s
\end{eqnarray}
can be further written as
\begin{eqnarray}\label{conservative}
\frac{\partial\phi}{\partial t}-\frac{\partial f_v}{\partial x}-\frac{\partial g_v}{\partial y}+\frac{\partial f_i}{\partial x}+\frac{\partial g_i}{\partial y}=s,
\end{eqnarray}
where the fluxes are defined as $f_i=c_1\phi$, $g_i=c_2\phi$, $f_v=a\frac{\partial\phi}{\partial x}$, $g_v=a\frac{\partial\phi}{\partial y}$. The time dependent mapping provided in Eq. (\ref{transformation}) allows one to transform the above equation to yield a strongly conservative formulation using the chain rule. The transformed governing equation in the computational plane is presented using the newly defined solution and fluxes and it takes the following form	
\begin{align}\label{conservative_tra}
\frac{\partial\phi^c}{\partial \tau}-\frac{\partial f^c_v}{\partial \xi}-\frac{\partial g^c_v}{\partial \eta}+\frac{\partial f^c_i}{\partial \xi}+&\frac{\partial g^c_i}{\partial \eta}=s+\phi\left(\frac{\partial J}{\partial \tau}-\frac{\partial J_1}{\partial \xi}-\frac{\partial J_2}{\partial \eta}\right)\nonumber\\
&+\left(c_1\phi-a\frac{\partial \phi}{\partial x}\right)\left(\frac{\partial (J\xi_x)}{\partial \eta}+\frac{\partial (J\eta_x)}{\partial \eta} \right)\nonumber\\
&+\left(c_2\phi-a\frac{\partial \phi}{\partial x}\right)\left(\frac{\partial (J\xi_y)}{\partial \eta}+\frac{\partial (J\eta_y)}{\partial \eta} \right)
\end{align}
where $\;\phi^c=J\phi$, $\;f_i^c=(c_1\xi_x+c_2\xi_y+\xi_t)\phi^c$, $\;f_i^c=(c_1\eta_x+c_2\eta_y+\eta_t)\phi^c$, $f_v^c=J(f_v\xi_x+g_v\xi_y)$, $\;g_v^c=J(f_v\eta_x+g_v\eta_y)$. From Eq. (\ref{conservative_tra}), it is clear that for the system to remain conservative it is imperative that
\begin{subequations}\label{GCL}
\begin{empheq}[left=\empheqlbrace]{align}
&\frac{\partial J}{\partial \tau}-\frac{\partial J_1}{\partial \xi}-\frac{\partial J_2}{\partial \eta}=0, \label{GCL1}\\
&\frac{\partial (J\xi_x)}{\partial \eta}+\frac{\partial (J\eta_x)}{\partial \eta}=0,
\label{GCL2}\\
&\frac{\partial (J\xi_y)}{\partial \eta}+\frac{\partial (J\eta_y)}{\partial \eta}=0.\label{GCL3}
\end{empheq}
\end{subequations}
The above set of equations ensure strong conservation for flow equations and they are referred to as GCL equations. Adherence to these equations rule out any apparent numerical source term in Eq. (\ref{conservative_tra}) and hence are important to accurately represent the transformed governing equation. Thus in any discretization of conservative form, these identities are to be satisfied numerically to avoid numerical instability. This has led to specialized methods for computing grid metrics that duly accommodate GCL relations \cite{sjo_yee_vin_14, abe_hag_non_16}. On the other hand, the transformed equation cast in non-conservative form as reported in Eq. (\ref{con_dif_tra1}) \textit{vis-a-vis} Eq. (\ref{con_dif}) might be immune to such numerical considerations and could turn out to be the one that is much more advantageous. Nevertheless, the computation of temporal and spatial grid metrics up to the respective desired levels of accuracy remain vital.
\subsection{Compact discretization}
The following central difference operators 
$$\delta_{\xi}\phi_{i,j}\equiv\frac{\phi_{i+1,j}-\phi_{i-1,j}}{2h},\;\;\;\delta_{\eta}\phi_{i,j}\equiv\frac{\phi_{i,j+1}-\phi_{i,j-1}}{2k},$$
$$\delta^2_{\xi}\phi_{i,j}\equiv\frac{\phi_{i+1,j}-2\phi_{i,j}+\phi_{i-1,j}}{h^2},\;\;\;\delta^2_{\eta}\phi_{i,j}\equiv\frac{\phi_{i,j+1}-2\phi_{i,j}+\phi_{i,j-1}}{k^2},$$
are used for ease in description. Here $h$ and $k$ are spatial interval in $\xi$ and $\eta$-direction, respectively, in the structured computational mesh and $\phi_{i,j}=\phi(\xi_i,\eta_j)$, with $\xi_i=\xi_{\min}+ih$, $\eta_j=\eta_{\min}+jk$, $0\le i \le N_{\xi}$ and $0\le j\le N_{\eta}$.

The second order unary and mixed partial derivatives appearing in the diffusive terms in the operator $[A]$ of Eq. (\ref{operator}) can be compactly approximated with fourth-order accuracy following their development in the work of Sen \cite{sen_13, sen_16}. These approximations amenable with variable coefficients possess good numerical characteristics and are given by
\begin{subequations}\label{com_aprox_2nd}
\begin{empheq}[left=\empheqlbrace]{align}
&\partial_{\xi\xi}\phi_{i,j}=2\delta^2_{\xi}\phi_{_{i,j}}-\delta_{\xi}\phi_{\xi_{i,j}}+O(h^4),  \\
&\partial_{\eta\eta}\phi_{i,j}=2\delta^2_{\eta}\phi_{_{i,j}}-\delta_{\eta}\phi_{\eta_{i,j}}+O(k^4), \\
&\partial_{\xi\eta}\phi_{i,j}=\delta_{\xi}\phi_{\eta_{i,j}}+\delta_{\eta}\phi_{\xi_{i,j}}-\delta_{\xi}\delta_{\eta}\phi_{i,j}+O(h^2k^2).
\end{empheq}
\end{subequations}
We thus define spatially fourth order accurate discrete operator $A_{h,k}$ as
\begin{eqnarray}\label{operator_dis}
A_{h,k}\phi_{i,j}=(-2\alpha_{1_{i,j}}\delta^2_{\xi}-2\alpha_{2_{i,j}}\delta^2_{\eta}+\beta_{_{i,j}}\delta_{\xi}\delta_{\eta})\phi_{_{i,j}}\nonumber\\
+(\alpha_{1_{i,j}}\delta_{\xi}-\beta_{_{i,j}}\delta_{\eta}+\chi_{1_{i,j}})\phi_{\xi_{i,j}}+(\alpha_{2_{i,j}}\delta_{\eta}-\beta_{_{i,j}}\delta_{\xi}+\chi_{2_{i,j}})\phi_{\eta_{i,j}}.
\end{eqnarray}
The above operators further require compatible fourth order accurate approximations of the gradients of $\phi$ and are accomplished by applying Pad\'{e} approximations 
\begin{subequations}\label{pade}
\begin{empheq}[left=\empheqlbrace]{align}
&(1+\frac{h^2}{6}\delta^2_{\xi})\phi_{\xi_{i,j}}=\delta_{\xi}\phi_{i,j}+O(h^4),  \\
&(1+\frac{k^2}{6}\delta^2_{\eta})\phi_{\eta_{i,j}}=\delta_{\eta}\phi_{i,j}+O(k^4).
\end{empheq}
\end{subequations}
Finally, second order accurate temporal discretization of the Eq. (\ref{con_dif_tra1}) is realized using Crank-Nicolson method. Thus, a fully discrete scheme for grid point $(i, j)$ at time level $(n)$ is
\begin{eqnarray}\label{fully_discrete}
\left(1+\frac{\delta \tau}{2} A^{(n+1)}_{h,k}\right)\phi^{(n+1)}_{i,j}=\left(1-\frac{\delta \tau}{2} A^{(n)}_{h,k}\right)\phi^{(n)}_{i,j}+\frac{\delta \tau}{2} \left(s^{(n+1)}_{i,j}+s^{(n)}_{i,j}\right).
\end{eqnarray}
The above formulation can be effectively implemented to discretize parabolic systems presented via Eq. (\ref{N_S_pri1}) and Eq. (\ref{N_S_psi1}) with convection velocity set to $\bm c=\bm u$. For steady elliptic equation 
\begin{equation}\label{elliptic}
-\Delta\phi=s(x,y)
\end{equation}
representing pressure Poisson and stream function equations (\ref{pre_poi}) and (\ref{N_S_psi2}) respectively at each temporal step the discrete form relates to
\begin{eqnarray}\label{elliptic_discrete}
\tilde{A}_{h,k}\phi_{i,j}=s_{i,j}
\end{eqnarray}
where discrete operator $\tilde{A}_{h,k}$ is
\begin{align}\label{elliptic_discrete}
\tilde{A}_{h,k}\phi_{i,j}=&(-2\tilde{\alpha}_{1_{i,j}}\delta^2_{\xi}-2\tilde{\alpha}_{2_{i,j}}\delta^2_{\eta}+\tilde{\beta}_{_{i,j}}\delta_{\xi}\delta_{\eta})\phi_{_{i,j}}\nonumber\\
&+(\tilde{\alpha}_{1_{i,j}}\delta_{\xi}-\tilde{\beta}_{_{i,j}}\delta_{\eta}+\tilde{\chi}_{1_{i,j}})\phi_{\xi_{i,j}}
+(\tilde{\alpha}_{2_{i,j}}\delta_{\eta}-\tilde{\beta}_{_{i,j}}\delta_{\xi}+\tilde{\chi}_{2_{i,j}})\phi_{\eta_{i,j}}
\end{align}
with
$$\tilde{\alpha}_1=\frac{1}{J^2}(x^2_{\eta}+y^2_{\eta}),\;\;
\tilde{\alpha}_2=\frac{1}{J^2}(x^2_{\xi}+y^2_{\xi}),\;\;
\tilde{\beta}=-\frac{2}{J^2}(x_{\xi}x_{\eta}+y_{\xi}x_{\eta}),$$
$$\tilde{\chi}_1=-\frac{1}{J^3}\bigg(J_{\eta}(x_{\xi}x_{\eta}+y_{\xi}y_{\eta})-J_{\xi}(x^2_{\eta}+y^2_{\eta})+J(x_{\eta}x_{\xi\eta}+y_{\eta}y_{\xi\eta}-x_{\xi}x_{\eta\eta}-y_{\xi}y_{\eta\eta}\bigg),$$
$$\tilde{\chi}_2=-\frac{1}{J^3}\bigg(J_{\xi}(x_{\xi}x_{\eta}+y_{\xi}y_{\eta})-J_{\eta}(x^2_{\xi}+y^2_{\xi})+J(x_{\xi}x_{\xi\eta}+y_{\xi}y_{\xi\eta}-x_{\eta}x_{\xi\xi}-y_{\eta}y_{\xi\xi}\bigg).$$
\subsection{Computation of grid metrics}
For the discretizations presented in Eqs. (\ref{fully_discrete}) and (\ref{elliptic_discrete}) to be able to retain spatial fourth order precision it is important to compute grid metrics to the same order of accuracy. In the absence of analytic expression grids in the physical plane are required to be generated numerically and hence discrete computation of grid metrics becomes imperative. As symmetric-conservative metric evaluation procedure ensuring automatic satisfaction of GCL identities is not necessitated for our discretized non-conservative system, computation of grid metrics is rather straightforward once the deforming grid is generated for a particular time step. To begin with, first order derivatives of the transformation are estimated using Pad\'{e} approximations at the interior points. It is well known in the literature that for a $p$th order interior scheme, boundary accuracy can be $(p-1)$st order without compromising global accuracy. Hence we work with one sided third order wide stencil approximations at boundary points. A similar process is repeated for computing second order derivatives appearing in Eq. (\ref{tran_coeff}). With the grid metrics available at all nodal points, the Jacobian can be explicitly found, following which gradients of the Jacobian can be estimated employing the Pad\'{e} approximation once again. Calculation of components of grid velocity having higher accuracy demands the storage of successive previous grids and is memory intensive. Except for one verification case simulations in this study are carried out by storing only one preceding grid thereby restricting computations to the first order accuracy although the overall discretization is temporally second order accurate.
\subsection{Grid movement}
The numerical grid movement is carried out using TFI and IWD procedure in the absence of analytic expression. For a few cases, the hybrid procedure was found to be best suited. The main steps associated with TFI and IDW interpolation are delineated below, details of which could be found in \cite{wan_prz_94, luk_col_bla_12, sen_nay_bre_17}.

Consider that there are a $Nb$ and $Ni$ number of boundary and interior nodes respectively in the physical domain $\partial\Omega\cup\Omega$. Given the displacements of the boundary points $\bm d(\bm x_b)$, TFI transmutation $\bm d(\bm x)$ at an interior point $\bm x$ is given by
\begin{equation}\label{TFI1}
  \bm d(\bm x)=\bm {F_1\oplus F_2}=\bm {F_1+F_2-F_1F_2}
\end{equation}
where $\bm F_i$'s are univariate projectors and $\bm F_i\bm F_j$ refer to their tensor product. Univariate projectors are obtained from $\bm d(\bm x_b)$ using transfinite formula based on arc-length to maintain the grid distribution laws of the initial grid details of which can be found in \cite{wan_prz_94, sen_nay_bre_17}. 

In contrast, IDW essentially interpolates translational and rotational displacements of the boundary nodes to the interior nodes. If $\bm M_b$ is the resultant rotation matrix and $\bm t_b$ is the associated translation vector for a boundary point having coordinate $\bm x_b$ at any temporal step of the evolving grid, then the displacement field $\bm d(\bm x_b)$ for that time step at that point could be written as
\begin{equation}\label{IDW1}
\bm x_b+\bm d(\bm x_b)=\bm M_b\bm x_b+\bm t_b \; .
\end{equation}
The previous position vector $\bm x_b$, subsequent position vector $\bm x_b+\bm d(\bm x_b)$ and the translation vector $\bm t_b$ can be used to define quaternions $R_{p_b}=[0, \bm x_b]$, $R_{s_b}=[0, \bm x_b+\bm d(\bm x_b)]$ and $T_b=[0, \bm t_b]$. These three quaternions are consequently used to define rotation quaternions $Q_b$ at each boundary point and are given by
\begin{equation}\label{IDW2}
R_{s_b} = Q_b \, R_{p_b} \, Q_b^* + T_b \;.
\end{equation}   
As rotation quaternions and the translation quaternions are estimated for all boundary grid points they are transferred to the interior grids by using the IDW interpolation methodology. This procedure interpolates a function at a given point $\bm
x$ by taking a weighted average of known functional values. For example, the displacements at the interior nodes can be calculated as
\begin{equation}\label{IDW3}
  \bm d(\bm x)=\frac{\sum\limits_{n=1}^{{Nb}} w_n(\bm x) \, \bm
   d(\bm x^n_b)}{\sum\limits_{n=1}^{{Nb}}w_n(\bm x)} \; ,
\end{equation}
where 
\begin{equation}\label{IDW4}
  w_n(\bm x)=A_n\bigg[\bigg(\frac{L_\text{def}}{\parallel\bm x-\bm
    x^n_b\parallel}\bigg)^a+\bigg(\frac{\alpha L_\text{def}}{\parallel\bm
    x-\bm x^n_b\parallel}\bigg)^b\bigg] \, ,
\end{equation}
is the weighting function. Here $A_n$ is the area weight of the $n$th boundary node. The best values for parameters $a$ and $b$ as 3 and 5, respectively. Luke \emph{et al.} \cite{luk_col_bla_12} defined $L_\text{def}$ and $\alpha$ as
\begin{equation}\label{IDW5}
  L_\text{def}=\max_{n=1}^{Nb}\parallel \bm x^n_b-\bm x_\text{centroid}\parallel \qquad \text{where} \qquad \bm x_\text{centroid}=\frac{1}{Nb+Ni}\sum_{n=1}^{Nb+Ni}\bm x^{n} \; ,
\end{equation}
and
\begin{equation}\label{IDW6}
\alpha=\frac{\eta}{L_\text{def}}\max_{n=1}^{Nb}\parallel \bm d(\bm x^n_b)-\bm d_\text{mean}\parallel \quad \text{with}
   \quad \bm d_\text{mean}=\sum_{n=1}^{Nb}a_n  \, \bm d(\bm x^n_b) 
\end{equation}
respectively. Further suggested values are $a_n=A_n\bigg/\sum\limits_{n=1}^{Nb}A_n$ and $\eta=5$. In this work $A_n$ is kept uniformly at one.  It is to be noted that Eq.~(\ref{IDW6}) is an option and could be altered by a fixed value of $\alpha$ appropriate to the deformation sought.

The above procedure applied to interpolate boundary rotation and translation quaternions is quite effective to deform mesh in situations involving a combination of translation and rotation and results in the generation of the smoothly varying grid. Here adequate care needs to be taken to normalize components of the rotation quaternions at the interior mesh points before finally estimating the displacement vectors at those grids.

\section{Numerical examples}
\subsection{Convection-diffusion of Gaussian pulse}
Verification of the newly developed scheme is carried out by simulating convection-diffusion of an initial Gaussian pulse 
\begin{equation}\label{GP1}
\phi_0(x,y)=e^{-a[(x-x_c)^2+(y-y_c)^2]}
\end{equation}
of unit height centred at $(x_c, y_c)$.
\subsubsection{Moving and deforming domain}
To an initial $[0,2]\times[0,2]$ domain discretized using Cartesian grid time-varying deformations are prescribed on all sides. On the left and right boundaries, identical displacements using $A(t)\sin(\pi\eta)$ are enforced. The top boundary is misshapen according to $2A(t)\sin(0.5\pi\xi)$ whereas the bottom boundary is deformed using bi-quadratic polynomial $2A(t)\xi^2(2-\xi)^2$. Time varying magnitude of deformation is set at $A(t)=0.2\sin(0.5\pi\tau/T)$. Moreover, all boundary points uniformly follow the pulse centred at $(0.5,0.5)$ with translational velocity $(0.7, 0.7)\tau/T$ resulting in a continuous displacement of the domain of interest. In the absence of any analytic approach, numerical grid deformation is envisaged in this test case. TFI is preferred as the problem involves the convex domain. 
\begin{figure}[!h]
\begin{minipage}[b]{.55\linewidth}\hspace{-1cm}
\centering\psfig{file=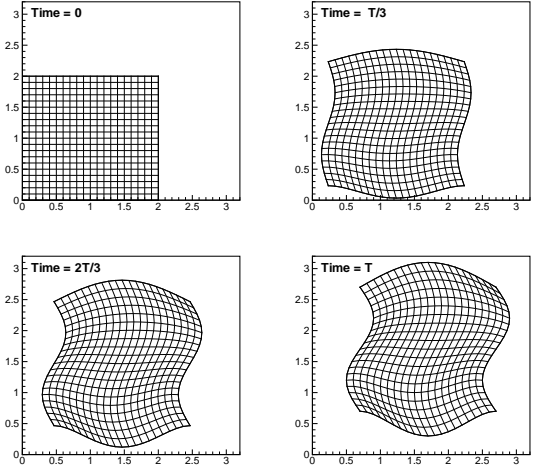,width=0.8\linewidth}
 \\(a)
\end{minipage}
\begin{minipage}[b]{.55\linewidth}\hspace{-1cm}
\centering\psfig{file=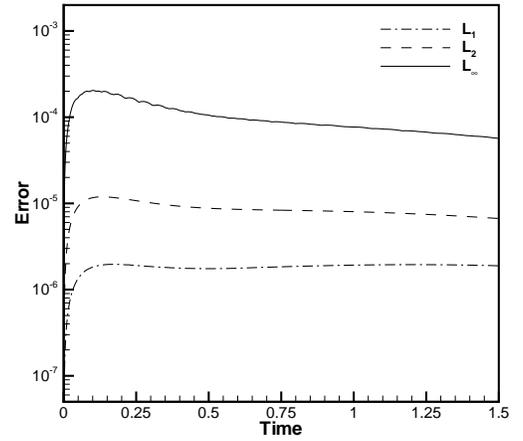,width=0.8\linewidth}
 \\(b)
\end{minipage}
\begin{center}
\caption{{\sl Problem 1: (a) Deformed and translating gird at four different time. (b) Time evolution of errors corresponding to three different norms.} }
    \label{fig:GP_grid_conve}
\end{center}
\end{figure}

Four stages of the resultant time-varying grid are presented in Fig. \ref{fig:GP_grid_conve}(a). In this figure, a significant deterioration of the orthogonality and translation of the grid can be noticed. To check convergence of the scheme we compute upto $T=1.5$ with $h^2=k^2=(1/40)^2=\delta\tau$ taking diffusion coefficient and convection coefficients fixed at $a=100$ and $c_1=80=c_2$ respectively. This choice corresponds to cell Peclet number $Pe=2$ in the computational plane. Isomorphic temporally decaying nature of error computed using three different norms could be found in Fig. \ref{fig:GP_grid_conve}(b). For this problem grid metrics are calculated using a fourth order central scheme at interior points and one sided wide stencil schemes at boundary points. As this problem is designed for verification we store three consecutive grids and estimate grid velocities with second order accuracy.
\begin{table}[!h]
\begin{center}
\caption{\sl {Problem 1: $L_1$-, $L_2$-, $L_{\infty}$- norm error and spatial order of convergence with $\delta \tau=h^2=k^2$.}}
%\footnotesize
{\begin{tabular}{ccccccc} \hline \hline
 Time  		&      	&$21\times21$ 		 &Order &$41\times41$  		 &Order &$81\times81$\\
\hline
$\tau=0.5$  &$L_1$  &6.612$\times10^{-4}$&4.85  &2.283$\times10^{-5}$&3.71  &1.750$\times10^{-6}$\\
         	&$L_2$  &2.275$\times10^{-3}$&4.32  &1.137$\times10^{-4}$&3.69  &8.786$\times10^{-6}$\\
     &$L_{\infty}$	&1.891$\times10^{-2}$&3.72  &1.435$\times10^{-3}$&3.76  &1.056$\times10^{-4}$\\
       		&       &                    &      &                    &      &                    \\
$\tau=1.0$  &$L_1$  &3.279$\times10^{-4}$&4.46  &1.490$\times10^{-5}$&2.95  &1.924$\times10^{-6}$\\
         	&$L_2$  &1.123$\times10^{-3}$&4.07  &6.706$\times10^{-5}$&3.06  &8.018$\times10^{-6}$\\
     &$L_{\infty}$  &1.007$\times10^{-2}$&3.81  &7.180$\times10^{-4}$&3.22  &7.689$\times10^{-5}$\\
\hline \hline
\end{tabular}}
\label{table:GP_order_spatial}
\end{center}
\end{table}
\begin{figure}[!h]
\begin{minipage}[b]{.55\linewidth}\hspace{-1cm}
\centering\psfig{file=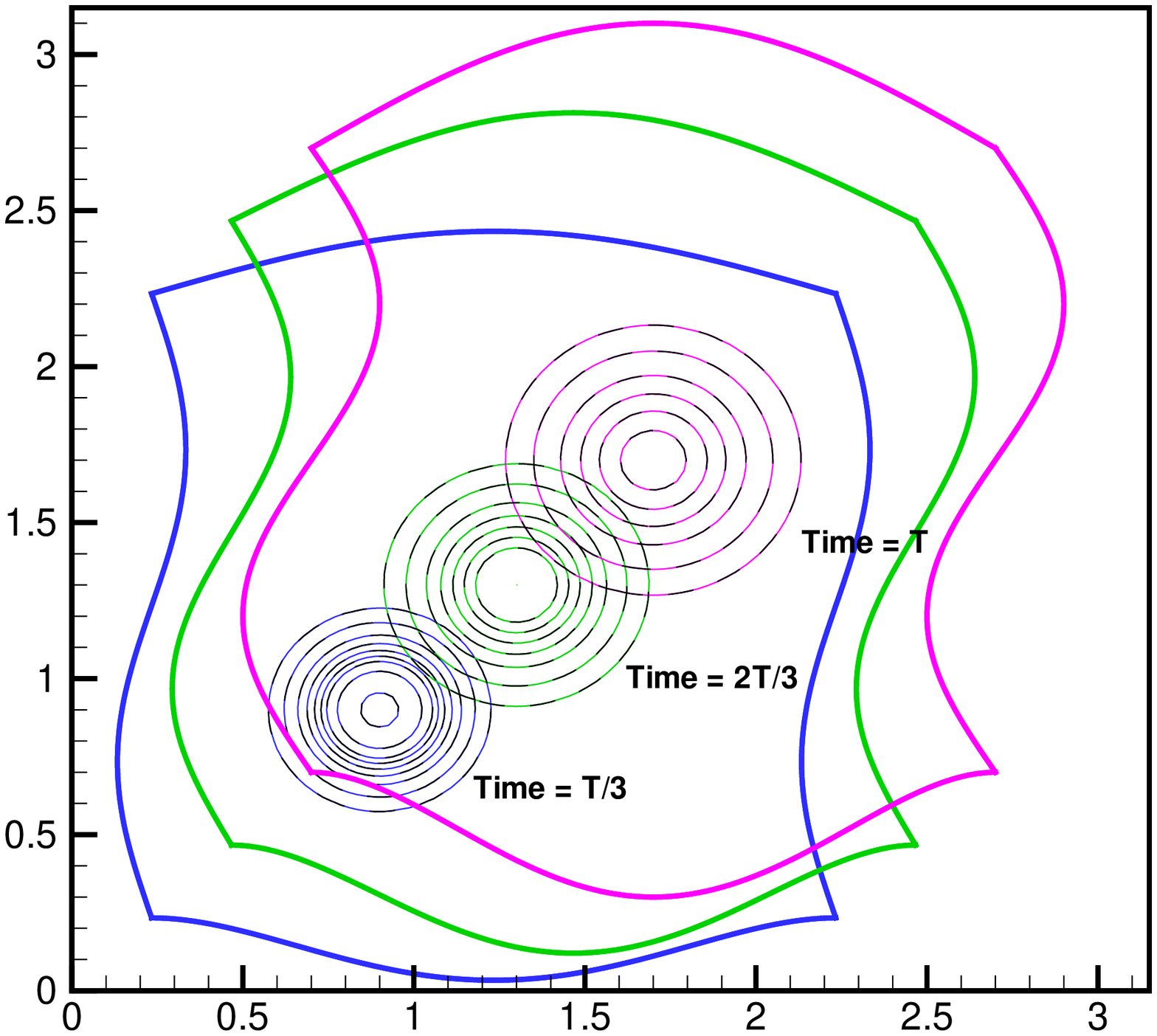,width=0.95\linewidth}
 \\(a)
\end{minipage}
\begin{minipage}[b]{.55\linewidth}\hspace{-1cm}
\centering\psfig{file=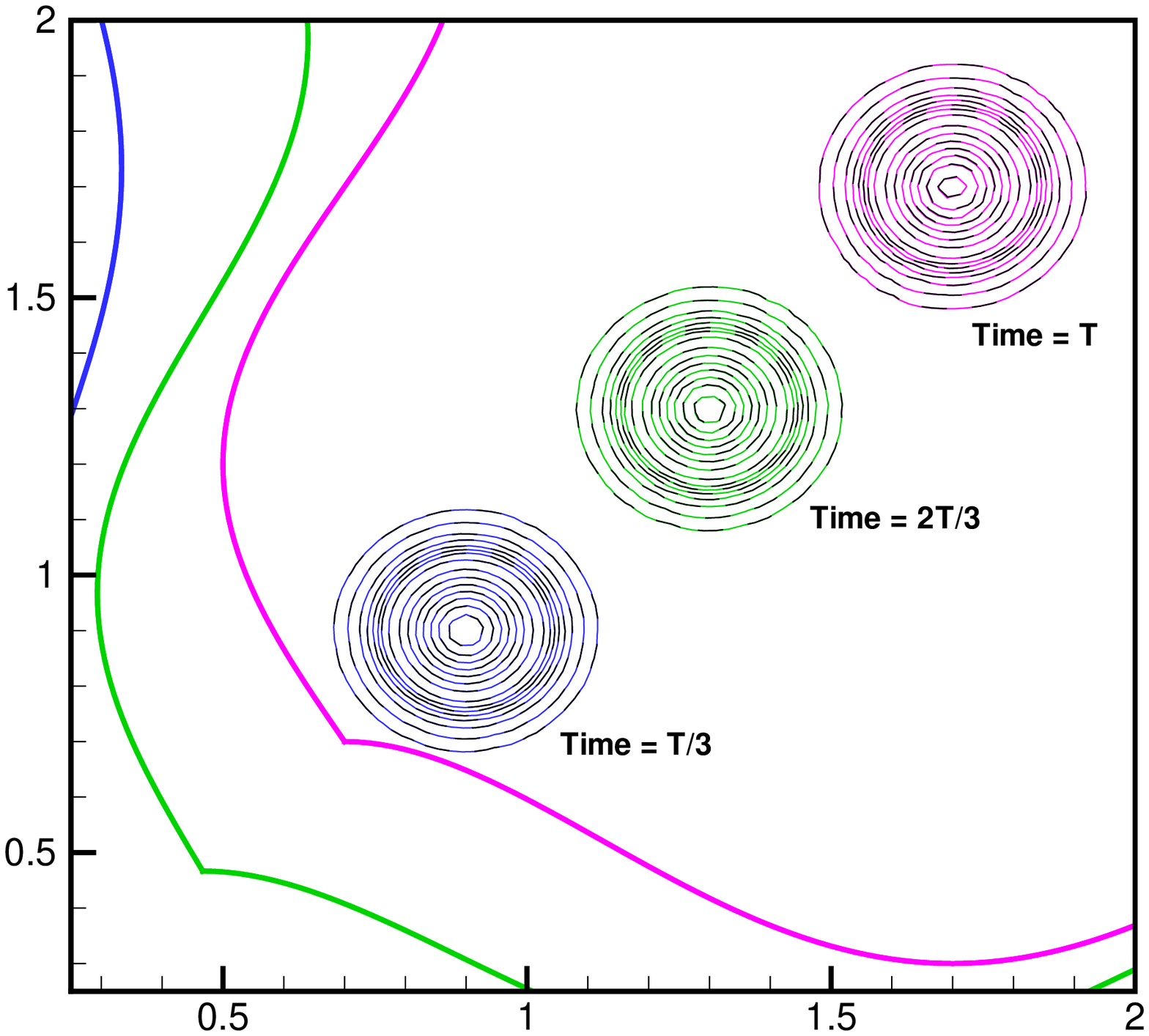,width=0.95\linewidth}
 \\(b)
\end{minipage}
\begin{center}
\caption{{\sl Problem 1: Comparison of  numerical solution (solid lines) and exact solution (black dashed lines) at three different times for (a) $Pe=2$ and (b) $Pe=200$.} }
    \label{fig:GP_num_com}
\end{center}
\end{figure}
\begin{table}[!h]
\begin{center}
\caption{\sl {Problem 1: $L_1$-, $L_2$-, $L_{\infty}$- norm error and temporal order of convergence.}}
%\footnotesize
{\begin{tabular}{ccccccc} \hline \hline
 Time  		&      	&$\delta\tau=0.05$ 	 &Order &$\delta\tau=0.025$  &Order &$\delta\tau=0.0125$ \\
\hline
$\tau=0.5$  &$L_1$  &1.366$\times10^{-3}$&1.85  &3.778$\times10^{-4}$&1.64  &1.211$\times10^{-4}$\\
          	&$L_2$  &6.557$\times10^{-3}$&1.79  &1.890$\times10^{-3}$&1.50  &6.672$\times10^{-4}$\\
     &$L_{\infty}$ 	&7.609$\times10^{-2}$&1.82  &2.150$\times10^{-2}$&1.48  &7.683$\times10^{-3}$\\
       		&       &                    &      &                    &      &                    \\
$\tau=1.0$  &$L_1$  &1.042$\times10^{-3}$&1.80  &2.992$\times10^{-4}$&1.54  &1.031$\times10^{-4}$\\
         	&$L_2$  &4.050$\times10^{-3}$&1.73  &1.219$\times10^{-3}$&1.43  &4.513$\times10^{-4}$\\
     &$L_{\infty}$  &3.799$\times10^{-2}$&1.77  &1.117$\times10^{-2}$&1.40  &4.240$\times10^{-3}$\\
\hline \hline
\end{tabular}}
\label{table:GP_order_temporal}
\end{center}
\end{table}
A numerical verification of the spatial order of accuracy in evolving domain is carried out in Table \ref{table:GP_order_spatial} using three different grids. Computations indeed support the theoretical order of convergence especially at the early time but are seen to report reduced order with time advancing. We further compare the numerical and analytic solutions in Fig. \ref{fig:GP_num_com}(a). The numerical solution obtained using $81\times81$ grid and $\delta\tau=6.25\times10^{-4}$ is almost indistinguishable from the analytical solution. Investigation on the temporal order of accuracy of the scheme could be found in Table \ref{table:GP_order_temporal}. For this estimation, additional grid translation is withdrawn and a relative finer grid spacing $h=k=0.05$ is employed to keep spatial error low. Three different temporal step sizes $\delta\tau = 0.05$, $0.025$ and $0.0125$ are used and the results obtain are shown in Table \ref{table:GP_order_temporal}. Numerical computations are seen to support the theoretical second order temporal accuracy of the scheme. To further compare our results we compute with higher convection velocity $c_1=10000=c_2$ using spatial grid spacing $h=0.05=k$ and $\delta\tau = 2.5\times10^{-4}$ corresponding to $Pe=200$. A numerical solution along with the exact solution computed with $T=0.012$ can be found in Fig. \ref{fig:GP_num_com}(b). Again a good comparison could be observed.
\subsubsection{Randomized mesh}
Here convection-diffusion of the Gaussian pulse centred at $(0.5, 0.5)$ through a randomized mesh is considered. Starting with an initial uniform mesh a new randomized mesh is generated at every time step by 20\% relative perturbation of the grid spacing in both $x$ and $y$-directions. 
\begin{figure}[htbp]
\begin{minipage}[b]{.55\linewidth}\hspace{-1cm}
\centering\psfig{file=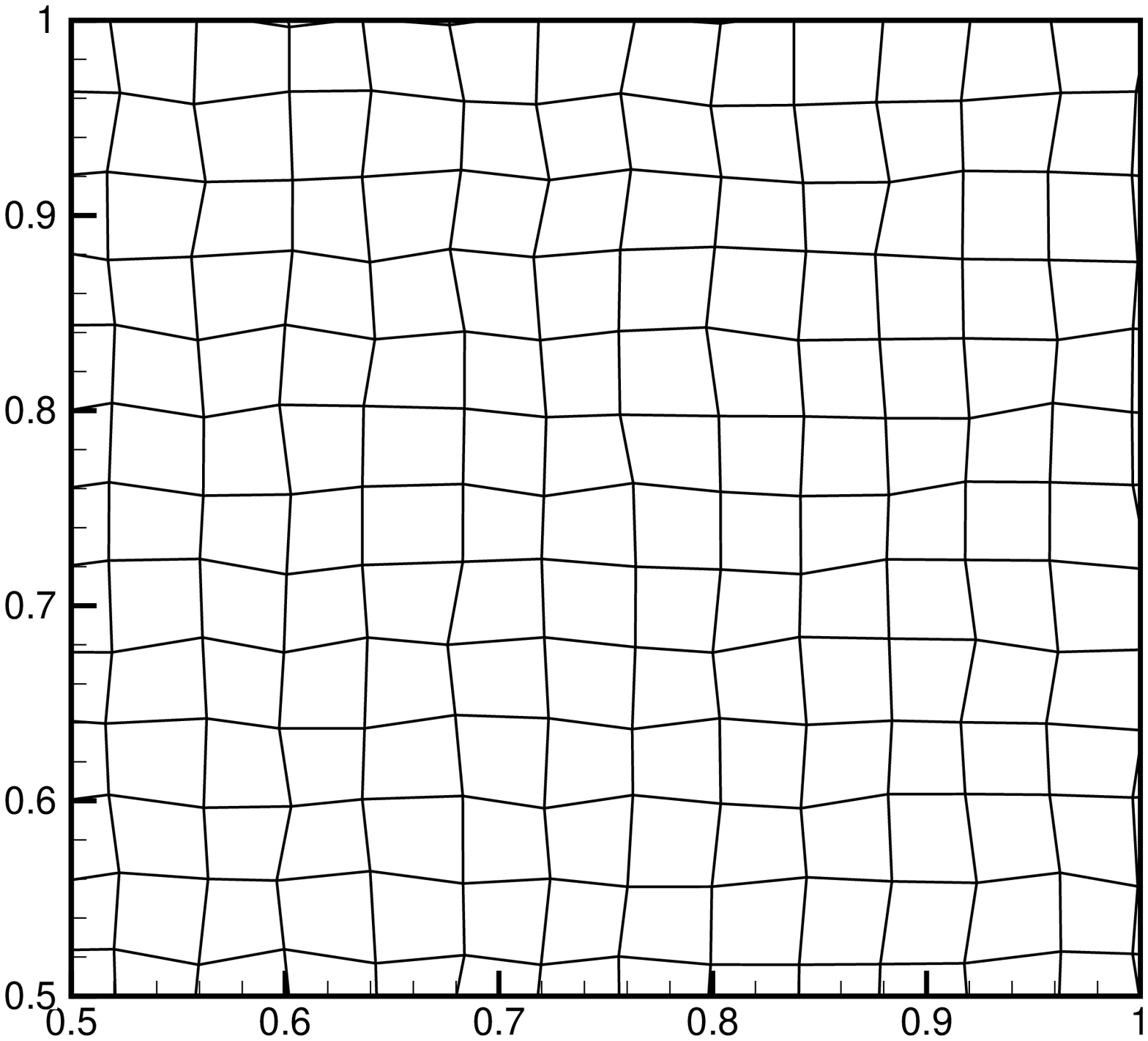,width=0.9\linewidth}
 \\(a)
\end{minipage}
\begin{minipage}[b]{.55\linewidth}\hspace{-1cm}
\centering\psfig{file=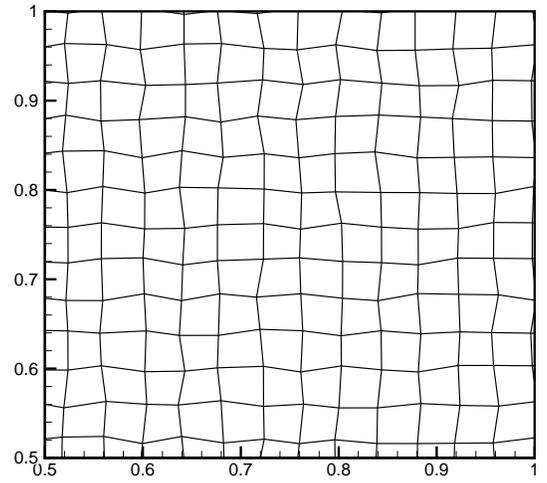,width=0.9\linewidth}
 \\(b)
\end{minipage}
\begin{minipage}[b]{.55\linewidth}\hspace{-1cm}
\centering\psfig{file=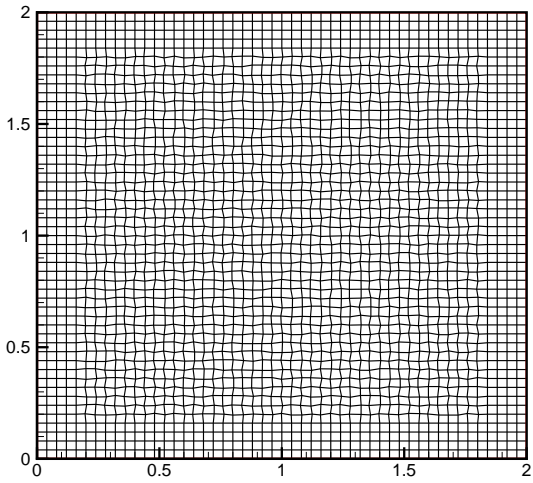,width=0.9\linewidth}
 \\(c)
\end{minipage}
\begin{minipage}[b]{.55\linewidth}\hspace{-1cm}
\centering\psfig{file=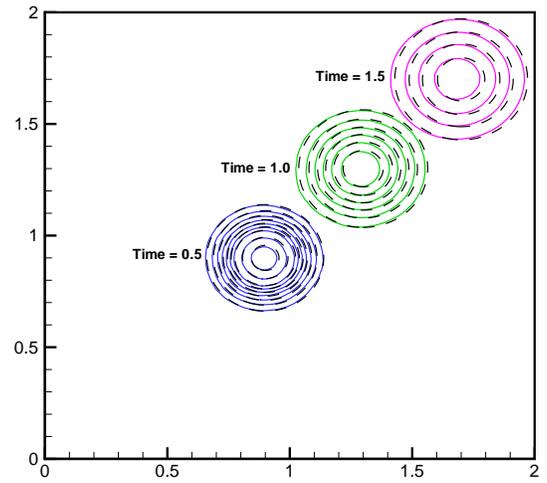,width=0.9\linewidth}
 \\(d)
\end{minipage}
\begin{center}
\caption{{\sl Problem 1: (a) and (b) Randomized mesh at two subsequent time steps. (c) Nominal $51\times51$ mesh. (d) Comparison of  numerical solution (solid lines) and exact solution (black dashed lines) at three different time.} }
    \label{fig:GP_rad}
\end{center}
\end{figure}
The resultant mesh, not only exhibits rapid and persistent non-smooth metric variations but is also different in each time step and might be distinct from the static randomized mesh considered in the work of Visbal and Gaitonde \cite{vis_gai_02}. Randomized variation of nominal mesh at two subsequent time steps are shown in Fig. \ref{fig:GP_rad}(a) and (b). Following \cite{vis_gai_02} metric continuity along all sides of the domain is maintained by leaving a few near boundary points unperturbed. The complete mesh can be found in Fig. \ref{fig:GP_rad}(c). Given the current formulation, where regularity and gradual variation of the grid is implicit, this mesh is an extreme test for correct simulation. We calculate for small Peclet number $Pe=0.25$ by choosing $a=100$, $c_1=80=c_2$ and $\delta\tau=1\times10^{-3}$. Numerical solutions obtained are presented in Fig. \ref{fig:GP_rad}(d) for three different temporal stages along with the exact solution. From this figure, it is clear that the current formulation is quite effective in simulating the motion of the pulse although some phase difference can be noticed.

\subsection{Freestream preservation in deforming domain}
Steger \cite{ste_78} in his study with conservative formulation found that freestream preservation characteristics of the compact methods are similar to that of the standard central difference schemes. Further Visbal and Gaitonde reported that compact schemes account for negligible metric cancellation errors when an identical discretization procedure is adopted for fluxes and the metrics \cite{vis_gai_02}. Nevertheless, working with conservative N-S equations in moving and deforming meshes the GCL identities as delineated in Eq. (\ref{GCL}) are of paramount interest to eliminate metric cancellation errors and thereby ensure freestream preservation. In the context of non-conservative formulation, it is important to document freestream preservation characteristics of compact schemes on highly distorted curvilinear meshes using a rather straightforward approach
to calculating the grid metrics. A comparative study is also envisaged by enforcing the metric identities following the newly developed symmetric-conservative metric evaluation procedure of Abe \emph{et al.} \cite{abe_hag_non_16}. Authors in \cite{abe_hag_non_16} analytically re-expressed coordinate transformation metrics and the procedure was found to ensure automatic satisfaction of GCL identities in 2D and is given by,
\begin{subequations}\label{GCL_comput}
\begin{empheq}[left=\empheqlbrace]{align}
&J=[(x_{\xi}y)_{\eta}-(x_{\eta}y)_{\xi}+(xy_{\eta})_{\xi}-(xy_{\xi})_{\eta}]/2, \label{GCL1_comput}\\
&J_1=[(x_{\tau}y)_{\eta}-(x_{\eta}y)_{\tau}+(xy_{\eta})_{\tau}-(xy_{\tau})_{\eta}]/2,\label{GCL2_comput}\\
&J_2=[(x_{\xi}y)_{\tau}-(x_{\tau}y)_{\xi}+(xy_{\tau})_{\xi}-(xy_{\xi})_{\tau}]/2.\label{GCL3_comput}
\end{empheq}
\end{subequations}

\begin{figure}[!h]
\begin{minipage}[b]{.5\linewidth}\hspace{-1cm}
\centering\psfig{file=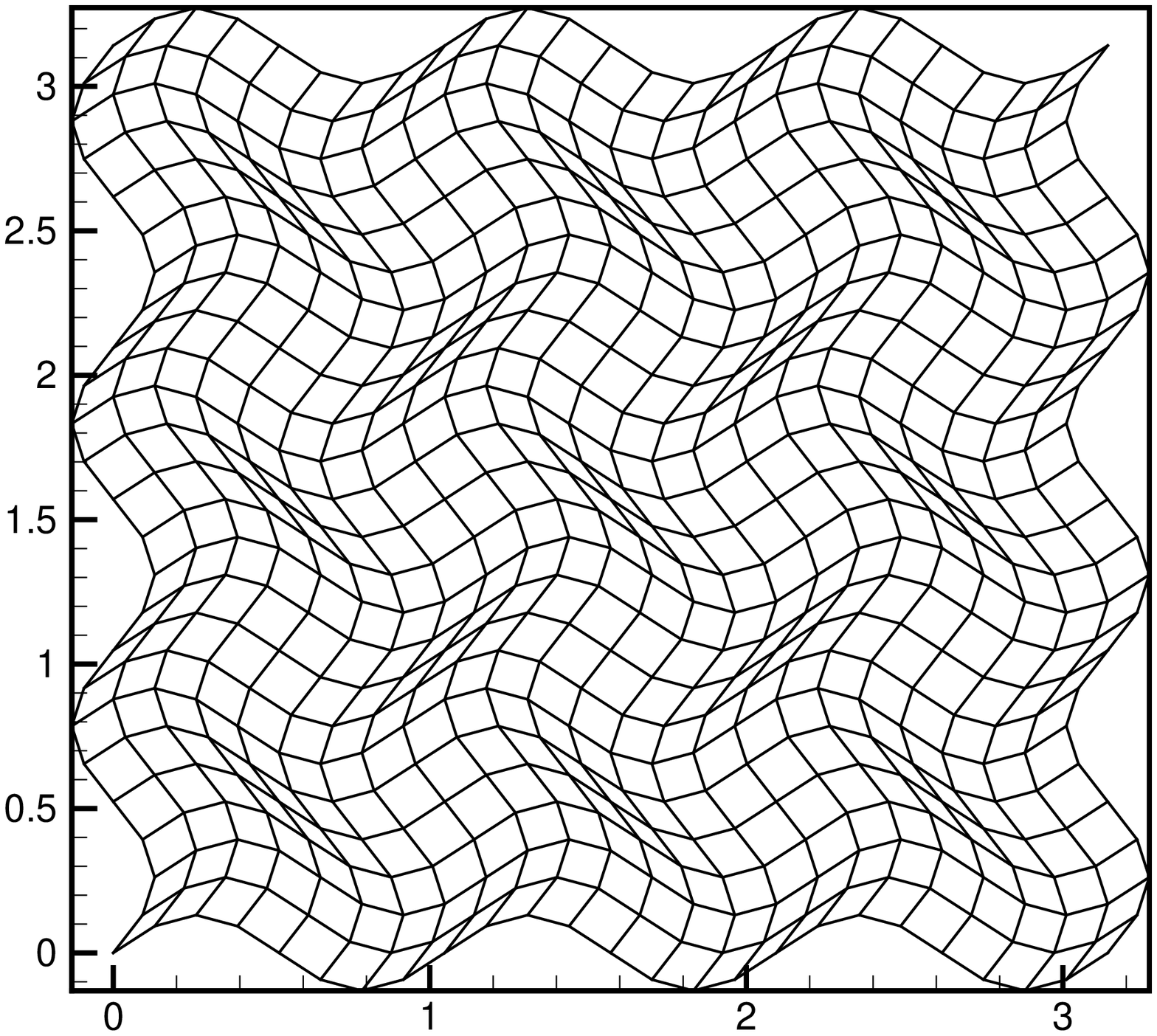,width=0.9\linewidth}
 \\(a)
\end{minipage}
\begin{minipage}[b]{.5\linewidth}\hspace{-1cm}
\centering\psfig{file=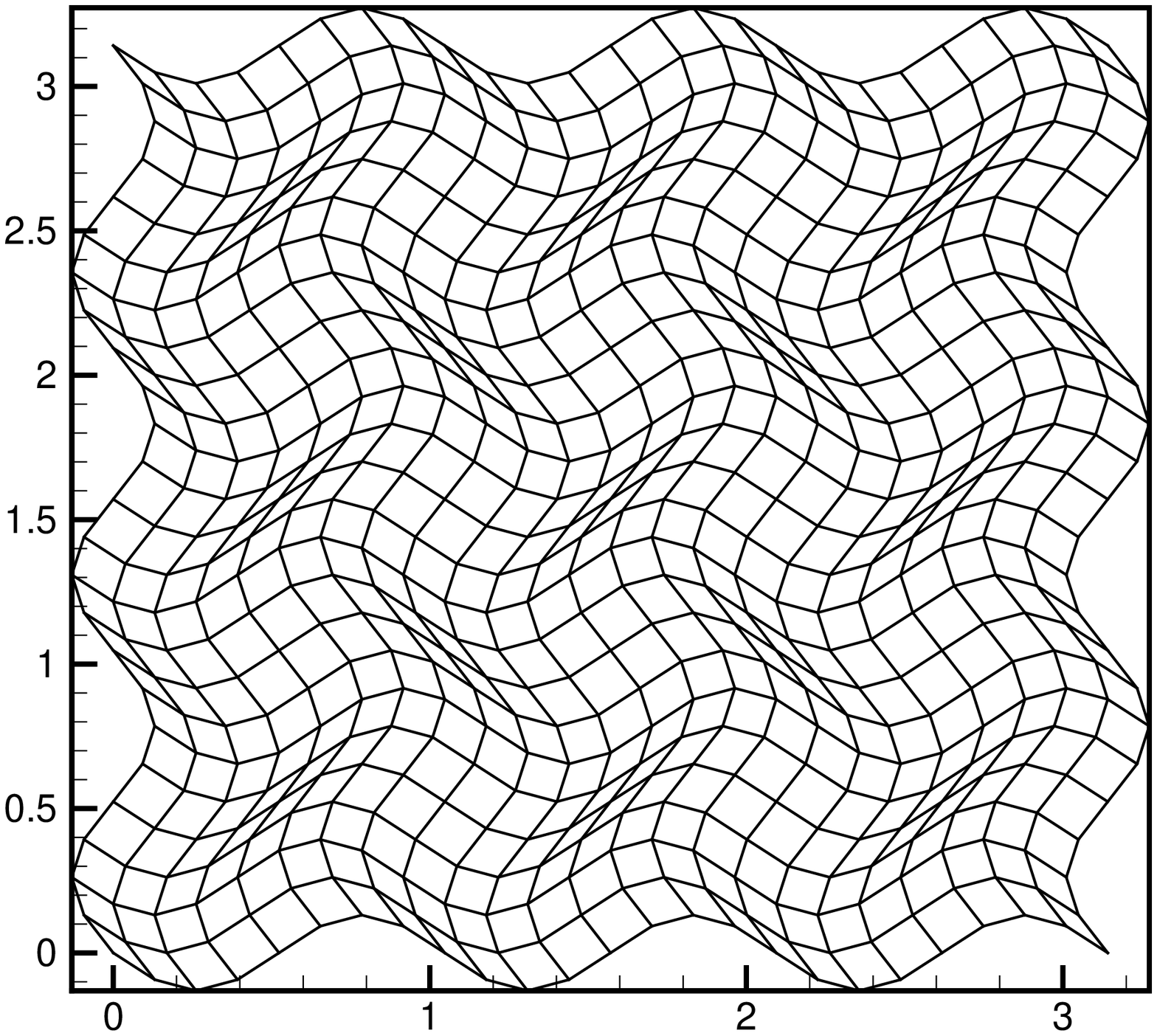,width=0.9\linewidth}
 \\(b)
\end{minipage}
\begin{center}
\caption{{\sl Problem 2: Nominal curvilinear wavy meshes at time (a) $\tau=1.0$ and (b) $\tau=3.0$.} }
    \label{fig:wavy_grid}
\end{center}
\end{figure}

Inviscid uniform flow with $u=1$, $v=0$ and pressure field obtained from pressure Poisson equation is computed using a wavy curvilinear time deforming grid. The grid is generated and subsequently deformed using analytic expression \cite{vis_gai_02}
\begin{subequations}\label{cur_grid}
\begin{empheq}[left=\empheqlbrace]{align}
&x=\xi+A_{\xi}\sin(2\pi\varpi\tau)\left(\frac{\xi_{\max}-\xi_{\min}}{N_{\xi}}\right)\sin\left(n_{\xi}\pi\frac{\eta-\eta_{\min}}{\eta_{\max}-\eta_{\min}}\right), \label{cur_grid1}\\
&y=\eta+A_{\eta}\sin(2\pi\varpi\tau)\left(\frac{\eta_{\max}-\eta_{\min}}{N_{\eta}}\right)\sin\left(n_{\eta}\pi\frac{\xi-\xi_{\min}}{\xi_{\max}-\xi_{\min}}\right),\label{cur_grid2}
\end{empheq}
\end{subequations}
with $A_{\xi}=1=A_{\eta}$, $\varpi=0.25$, $\xi_{\min}=0=\eta_{\min}$, $\xi_{\max}=\pi=\eta_{\max}$. Computations are carried out using two different grids $25\times25$ with $n_{\xi}=4=n_{\eta}$ and $65\times65$ with $n_{\xi}=6=n_{\eta}$ taking $\delta \tau=0.025$ and $0.0025$ respectively upto $T=2.5$ following \cite{vis_gai_02}. Presence of $N_{\xi}$ and $N_{\eta}$ act as stabilizing skewness of the grid with mesh refinement. Nominal grids at two different times are presented in Fig. \ref{fig:wavy_grid}. 
For this test coordinate transformation metrics are approximated using two different approaches. In the first case compact fourth order Pad\'{e} discretizations are used. Additionally, the Jacobian of the transformation $J$, as well as $J_1$ and $J_2$ are arrived at by using GCL identities given by the Eq. (\ref{GCL_comput}) to enforce metric cancellation. In the second approach, straightforward fourth order central differencing is employed. The results obtained are presented in Table \ref{table:freestream_com}.
\begin{table}
\vspace{-0.2cm}
\begin{center}    	
\caption{\sl {Problem 2: Freestream preservation errors for 2D wavy mesh}}
\begin{tabular}{c c c c c c}
  \hline
 	Grid		&\multicolumn{2}{c}{$25\times25$} 				&&\multicolumn{2}{c}{$65\times65$} \\
				 \cline{2-3} 						 			  \cline{5-6}  	
		 		& Error $||v||_{\infty}$& Error $||p||_{\infty}$&& Error $||v||_{\infty}$& Error $||p||_{\infty}$\\
		 		 \hline
  With GCL 		& $4\times10^{-30}$ 	& $4\times10^{-30}$ 	&& $3\times10^{-29}$ 	& $1\times10^{-29}$\\
 Without GCL 	& $5\times10^{-30}$ 	& $6\times10^{-30}$ 	&& $9\times10^{-29}$ 	& $8\times10^{-29}$\\
\hline
\end{tabular}
\label{table:freestream_com}
\end{center}
\end{table}
From the table, it is clear that free stream is preserved for the current formulation even without GCL identities being enforced. This correlates well with our theoretical understanding developed in section 2.3. 

\subsection{N-S equation with analytic solution}
\begin{figure}[htbp]
\begin{minipage}[b]{.55\linewidth}\hspace{-1cm}
\centering\psfig{file=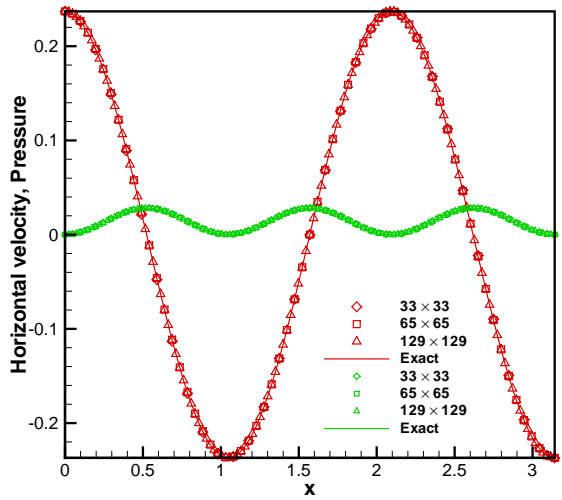,width=0.9\linewidth}
 \\(a)
\end{minipage}
\begin{minipage}[b]{.55\linewidth}\hspace{-1cm}
\centering\psfig{file=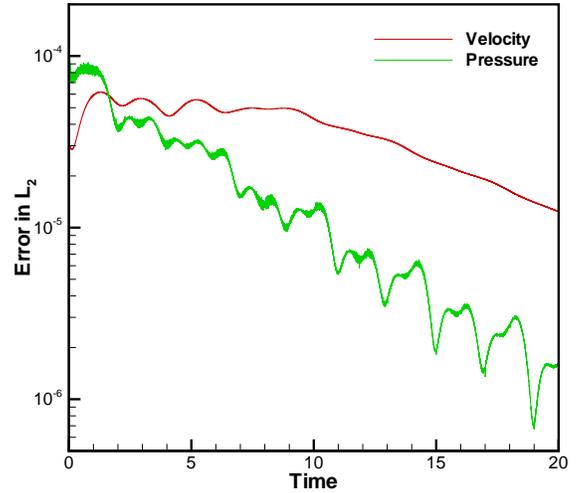,width=0.9\linewidth}
 \\(b)
\end{minipage}
\begin{minipage}[b]{.55\linewidth}\hspace{-1cm}
\centering\psfig{file=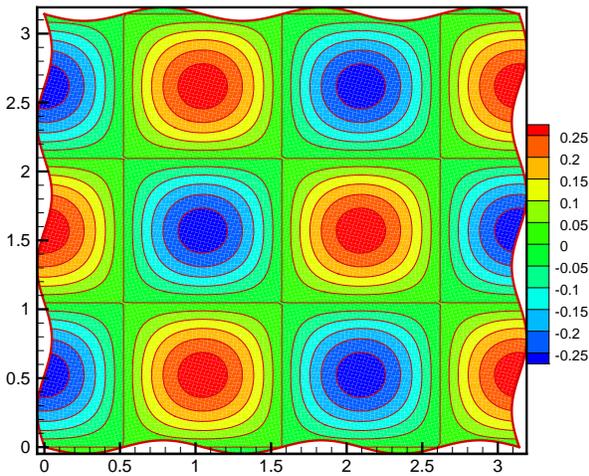,width=0.9\linewidth}
 \\(c)
\end{minipage}
\begin{minipage}[b]{.55\linewidth}\hspace{-1cm}
\centering\psfig{file=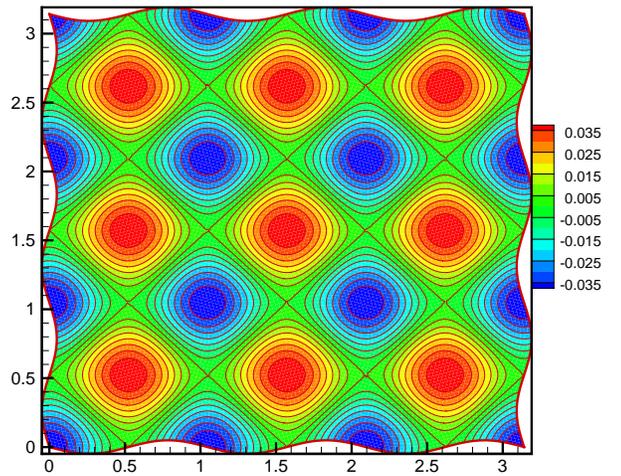,width=0.9\linewidth}
 \\(d)
\end{minipage}
\begin{center}
\caption{{\sl Problem 3: (a) Grid independence of $u_1$ and $p$ at $\tau=8.0$, (b) convergence of $u_1$ and $p$ in $L^2$- norm, (c) $u_1$ contour and (d) $p$ contour at $\tau=7.0$.} }
    \label{fig:FDV_mis}
\end{center}
\end{figure}

Taylor in 1923 established an exact solution of the incompressible N-S equations \cite{tay_23}. Since then it has been used to test numerical schemes applied to fluid flow problems. In this study, we further adapt this problem for a time-varying domain.
We consider the flow decayed by viscosity problem governed by the N-S equations having analytical solutions
\begin{subequations}\label{FDV_sol}
\begin{empheq}[left=\empheqlbrace]{align}
&u_1(x,y,t)=-\cos(Nx)\sin(Ny)e^{-2N^2t/Re}, \label{FDV_sol1}\\
&u_2(x,y,t)=\sin(Nx)\cos(Ny)e^{-2N^2t/Re},\label{FDV_sol2}\\
&p(x,y,t)=-\frac{1}{4}\left(\cos(2Nx)+\cos(2Ny)\right)e^{-4N^2t/Re}.\label{FDV_sol3}
\end{empheq}
\end{subequations}
Here $N$ is a positive integer and controls the numbers of vortices generated in the flow domain. The problem is considered in a time-varying physical domain $\Omega(t)$ that could be readily mapped to a fixed and regular parametric domain $\mathcal{D}=[0,\pi]\times[0,\pi]$ by setting
\begin{subequations}\label{cur_grid}
\begin{empheq}[left=\empheqlbrace]{align}
&x=\xi+A_{\xi}\sin(2\pi\varpi\tau)\sin\left(n_{\xi}\eta\right), \label{cur_grid1}\\
&y=\eta+A_{\eta}\sin(2\pi\varpi\tau)\sin\left(n_{\eta}\xi\right).\label{cur_grid2}
\end{empheq}
\end{subequations}
The above transformation leads to a computationally somewhat more challenging grid for finer meshes \textit{vis-a-vis} the transformation proposed by Visbal and Gaitonde \cite{vis_gai_02} and used in the previous section. This is due to the invariance of maximum deformation with respect to grid spacing leading to increasingly skewed mesh. Here we set $A_{\xi}=\pi/64=A_{\eta}$, $\varpi=0.25$ and $n_{\xi}=6=n_{\eta}$. Further, we take $N=3$, which leads to the creation of $N^2=9$ vortices in the velocity field and $2N^2=18$ vortices in the pressure field. Computations are done for $Re=100$ with three different grids $33\times33$, $65\times65$ and $129\times129$ and the spatial order of accuracy of horizontal velocity component and pressure estimated with $L^2$-norm is verified in Table \ref{table:FDV_order_spatial}. Vertical velocity is not presented as it has features similar to horizontal velocity. Order of accuracy around three can be noticed in  Table \ref{table:FDV_order_spatial} and is found to decrease with a decrease in grid spacing. This might correlate to increased skewness of the refined mesh.
\begin{table}[!h]
\begin{center}
\caption{\sl {Problem 3: $L_2$- norm error and spatial order of convergence with $\delta \tau=h^2=k^2$.}}
%\footnotesize
{\begin{tabular}{ccccccc} \hline \hline
 Time  		&      	&$33\times33$ 		 &Order &$65\times65$  		 &Order &$129\times129$\\
\hline
%$\tau=5.0$  &$u_1$  &6.030$\times10^{-3}$&2.98  &7.660$\times10^{-4}$&3.71  &--$\times10^{-6}$\\
%         	&$p$  	&3.522$\times10^{-3}$&2.80  &5.064$\times10^{-4}$&3.69  &--$\times10^{-6}$\\
%       		&       &                    &      &                    &      &                    \\
%$\tau=8.0$  &$u_1$  &6.712$\times10^{-3}$&3.13  &7.648$\times10^{-4}$&2.95  &--$\times10^{-6}$\\
%         	&$p$	&1.035$\times10^{-3}$&2.44  &1.903$\times10^{-4}$&3.06  &--$\times10^{-6}$\\ 
%       		&       &                    &      &                    &      &                    \\
%$\tau=10.0$ &$u_1$  &6.157$\times10^{-3}$&3.17  &6.843$\times10^{-4}$&2.95  &--$\times10^{-6}$\\
%         	&$p$	&1.316$\times10^{-3}$&2.67  &2.065$\times10^{-4}$&3.06  &--$\times10^{-6}$\\

$\tau=5.0$  &$u_1$  &3.386$\times10^{-3}$&3.15  &3.813$\times10^{-4}$&2.80  &5.493$\times10^{-5}$\\
         	&$p$	&1.712$\times10^{-3}$&3.10  &1.990$\times10^{-4}$&2.69  &3.091$\times10^{-5}$\\
       		&       &                    &      &                    &      &                    \\
$\tau=8.0$  &$u_1$  &3.825$\times10^{-3}$&3.29  &3.920$\times10^{-4}$&3.00  &4.908$\times10^{-5}$\\
         	&$p$	&5.362$\times10^{-4}$&3.04  &6.528$\times10^{-5}$&2.26  &1.360$\times10^{-5}$\\ 
       		&       &                    &      &                    &      &                    \\
$\tau=10.0$ &$u_1$  &3.431$\times10^{-3}$&3.29  &3.486$\times10^{-4}$&3.02  &4.306$\times10^{-5}$\\
         	&$p$	&6.397$\times10^{-4}$&2.99  &8.056$\times10^{-5}$&2.66  &1.272$\times10^{-5}$\\
\hline \hline
\end{tabular}}
\label{table:FDV_order_spatial}
\end{center}
\end{table}
\begin{table}[!h]
\begin{center}
\caption{\sl {Problem 3: $L_2$- norm error and temporal order of convergence.}}
%\footnotesize
{\begin{tabular}{ccccccc} \hline \hline
 Time  		&      	&$\delta\tau=0.04$ 	 &Order &$\delta\tau=0.02$   &Order &$\delta\tau=0.01$ \\
\hline
%$\tau=5.0$  &$u_1$  &7.063$\times10^{-3}$&1.20  &3.066$\times10^{-3}$&1.64  &1.413$\times10^{-3}$\\
%          	&$p$    &3.498$\times10^{-3}$&1.29  &1.426$\times10^{-3}$&1.01  &7.092$\times10^{-4}$\\
%       		&       &                    &      &                    &      &                    \\
%$\tau=8.0$  &$u_1$  &6.671$\times10^{-3}$&1.22  &2.864$\times10^{-3}$&1.15  &1.295$\times10^{-3}$\\
%         	&$p$	&2.389$\times10^{-3}$&0.93  &1.251$\times10^{-3}$&0.93  &6.552$\times10^{-4}$\\               	
%       		&       &                    &      &                    &      &                    \\
%$\tau=10.0$ &$u_1$  &5.754$\times10^{-3}$&1.20  &2.504$\times10^{-3}$&1.11  &1.162$\times10^{-3}$\\
%         	&$p$	&1.642$\times10^{-3}$&0.92  &8.700$\times10^{-4}$&0.86  &4.798$\times10^{-4}$\\

$\tau=5.0$  &$u_1$  &4.684$\times10^{-3}$&1.35  &1.842$\times10^{-3}$&1.26  &7.670$\times10^{-4}$\\
          	&$p$    &1.922$\times10^{-3}$&1.39  &7.354$\times10^{-4}$&1.21  &3.183$\times10^{-4}$\\
       		&       &                    &      &                    &      &                    \\
$\tau=8.0$  &$u_1$  &3.999$\times10^{-3}$&1.36  &1.558$\times10^{-3}$&1.30  &6.330$\times10^{-4}$\\
         	&$p$	&1.120$\times10^{-3}$&1.19  &4.921$\times10^{-4}$&1.12  &2.257$\times10^{-4}$\\               	
       		&       &                    &      &                    &      &                    \\
$\tau=10.0$ &$u_1$  &3.329$\times10^{-3}$&1.35  &1.306$\times10^{-3}$&1.27  &5.430$\times10^{-4}$\\
         	&$p$	&7.492$\times10^{-4}$&1.17  &3.326$\times10^{-4}$&1.02  &1.637$\times10^{-4}$\\     	
\hline \hline
\end{tabular}}
\label{table:FDV_order_temporal}
\end{center}
\end{table}
We investigate the temporal accuracy of the scheme in conjunction with this prototype N-S equation in Table \ref{table:FDV_order_temporal} where computations were carried out using $65\times65$ grid using three different time steps. Our computation report temporal accuracy between one and two for this nonlinear problem. In Fig. \ref{fig:FDV_mis}(a) grid independence of our computation is shown. That the exact solution is approached with grid refinement can be seen in this figure. Convergence of velocity and pressure field in $L_2$- norm is depicted in Fig. \ref{fig:FDV_mis}(b). Although temporal convergence of velocity is linear after some time, pressure is seen to follow an oscillatory path. Finally horizontal velocity contours and pressure contours at time $\tau=7.0$ are presented in Figs. \ref{fig:FDV_mis}(c) and \ref{fig:FDV_mis}(d) respectively. The results indeed correlate well with the exact solution.

\subsection{Driven cavity problem}
The 2D driven cavity problem defined on a unit square is often used as the benchmark test case for incompressible flow solvers is considered to demonstrate the accuracy of the proposed discretization on curvilinear time deforming grid. We carry out numerical simulations on both Cartesian uniform and time varying stretched curvilinear grids. Results obtained are essentially indistinguishable and are further compared with benchmark solutions to amply demonstrate the accuracy of the procedure advocated here. Such an indirect approach to validate the accuracy of a discretization strategy could be traced in the works of Wesseling \emph{et al.} \cite{wes_seg_kas_99} and  Ge and Sotiropoulos \cite{ge_sot_07} where the curvilinear grids are time independent and in the work of Chen and Xie \cite{che_xie_16} where the strategy was furthered to a time varying grid. 

\begin{figure}[!h]
\begin{minipage}[b]{.5\linewidth}\hspace{-1cm}
\centering\psfig{file=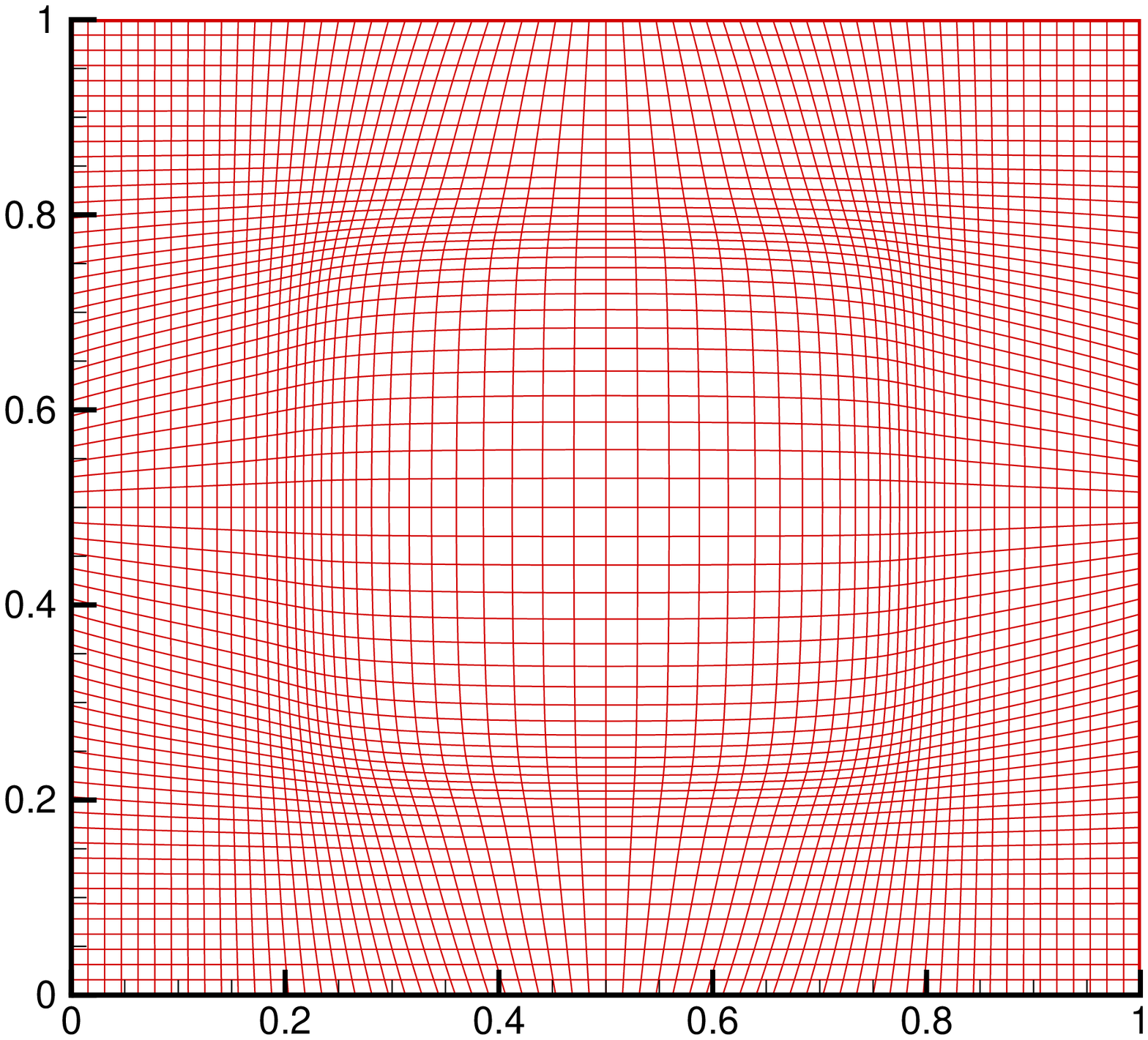,width=0.9\linewidth}
 \\(a)
\end{minipage}
\begin{minipage}[b]{.5\linewidth}\hspace{-1cm}
\centering\psfig{file=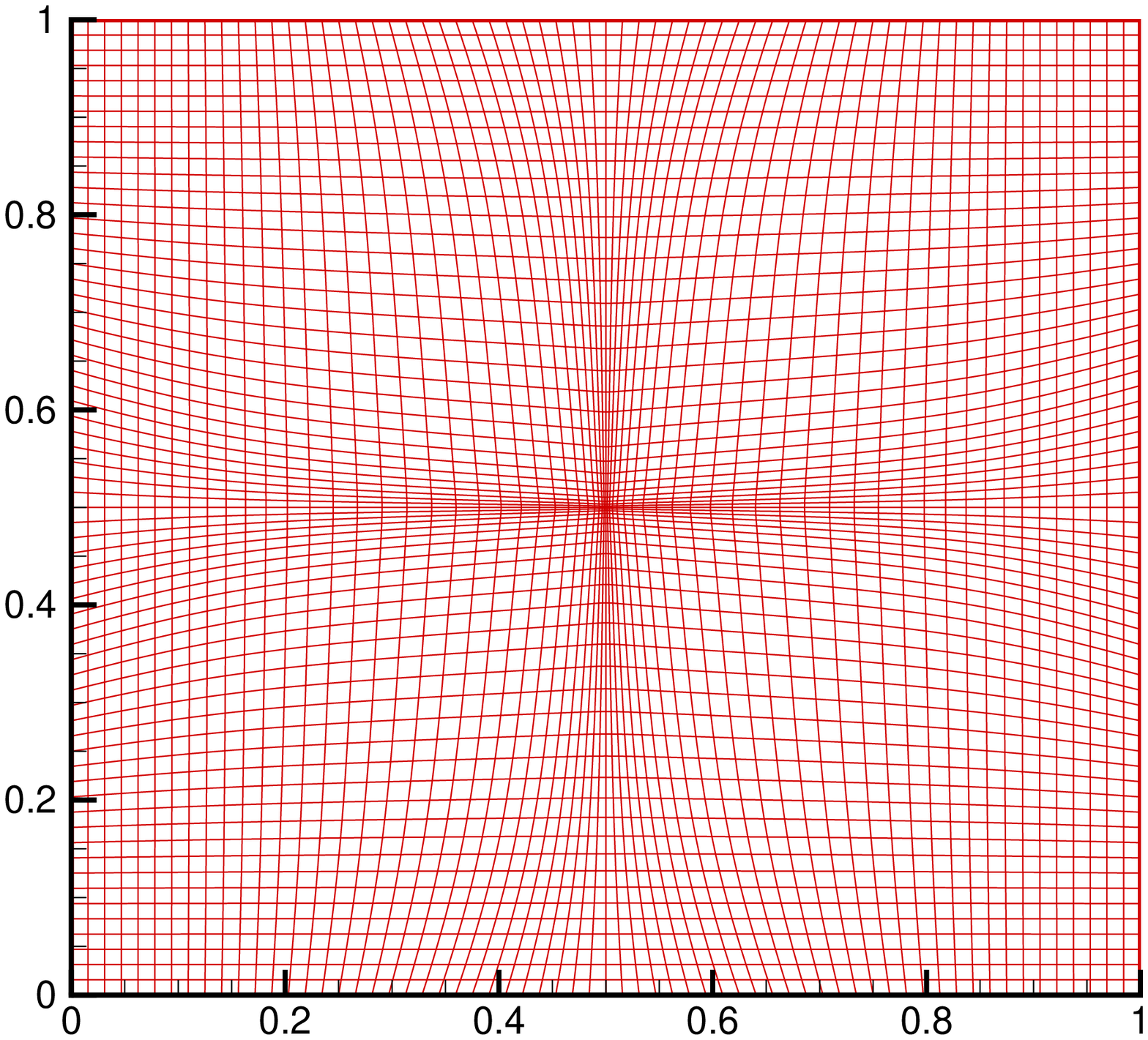,width=0.9\linewidth}
 \\(b)
\end{minipage}
\begin{center}
\caption{{\sl Problem 4: Deforming $65\times65$ mesh at time (a) $\tau=0.25$ and (b) $\tau=0.75$.} }
    \label{fig:ldc_grid}
\end{center}
\end{figure}

Employing primitive variable form of the Navier-Stokes equation steady-state solution is arrived at for diverse $Re$. The time deforming curvilinear grid over a rectangular domain $\Omega=[\xi_{min}, \xi_{max}]\times[\eta_{min}, \eta_{max}]$ is constructed combining polynomial stretching with trigonometric temporal variation and is given by  
\begin{subequations}\label{ldc_grid}
\begin{empheq}[left=\empheqlbrace]{align}
x=\left(\frac{\xi-\xi_{\scriptscriptstyle min}}{\xi_{\scriptscriptstyle max}-\xi_{\scriptscriptstyle min}}\right)+&A_{\xi}\left(\frac{\xi-\xi_{\scriptscriptstyle min}}{\xi_{\scriptscriptstyle max}-\xi_{\scriptscriptstyle min}}\right)^6\left(\frac{\xi-0.5(\xi_{\scriptscriptstyle min}+\xi_{\scriptscriptstyle max})}{\xi_{\scriptscriptstyle max}-\xi_{\scriptscriptstyle min}}\right)\left(\frac{\xi_{\scriptscriptstyle max}-\xi}{\xi_{\scriptscriptstyle max}-\xi_{\scriptscriptstyle min}}\right)^6\nonumber\\
&\left(\frac{\eta-\eta_{\scriptscriptstyle min}}{\eta_{\scriptscriptstyle max}-\eta_{\scriptscriptstyle min}}\right)\left(\frac{\eta_{max}-\eta}{\eta_{\scriptscriptstyle max}-\eta_{\scriptscriptstyle min}}\right)\sin(2\pi\tau), \label{ldc_grid1}\\
y=\left(\frac{\eta-\eta_{\scriptscriptstyle min}}{\eta_{\scriptscriptstyle max}-\eta_{\scriptscriptstyle min}}\right)+&A_{\eta}\left(\frac{\eta-\eta_{\scriptscriptstyle min}}{\eta_{\scriptscriptstyle max}-\eta_{\scriptscriptstyle min}}\right)^6\left(\frac{\eta-0.5(\eta_{\scriptscriptstyle min}+\eta_{\scriptscriptstyle max})}{\eta_{\scriptscriptstyle max}-\eta_{\scriptscriptstyle min}}\right)\left(\frac{\eta_{\scriptscriptstyle max}-\eta}{\eta_{\scriptscriptstyle max}-\eta_{\scriptscriptstyle min}}\right)^6\nonumber\\
&\left(\frac{\xi-\xi_{\scriptscriptstyle min}}{\xi_{\scriptscriptstyle max}-\xi_{\scriptscriptstyle min}}\right)\left(\frac{\xi_{\scriptscriptstyle max}-\xi}{\xi_{\scriptscriptstyle max}-\xi_{\scriptscriptstyle min}}\right)\sin(2\pi\tau).\label{ldc_grid2}
\end{empheq}
\end{subequations}
In our computations unit square is contemplated taking $\xi_{min}=0=\eta_{min}$ and $\xi_{max}=1=\eta_{max}$ with $A_{\xi}=15000=A_{\eta}$. Evolution of a $65\times65$ grid at time $0.25$ and $0.75$ corresponding to extreme stretching and clustering are presented in Figs. \ref{fig:ldc_grid}(a) and \ref{fig:ldc_grid}(b) respectively. As the grid deform continuously with a unit period, recurrence of such extremal twist along with loss of grid orthogonality could be noticed in the interior of the cavity at each nondimensional time.

\begin{figure}[htbp]
\begin{minipage}[b]{.55\linewidth}\hspace{-1cm}
\centering\psfig{file=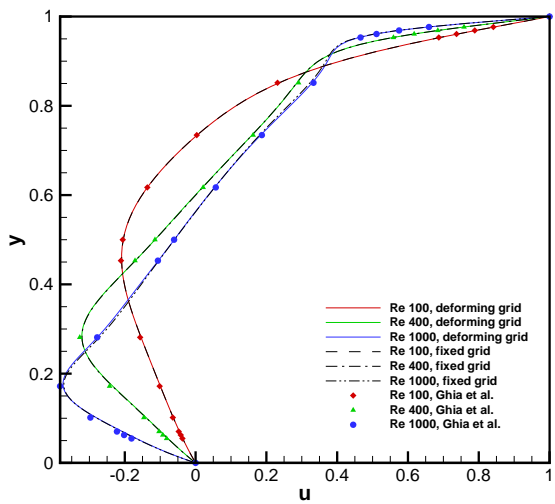,width=0.9\linewidth}
 \\(a)
\end{minipage}
\begin{minipage}[b]{.55\linewidth}\hspace{-1cm}
\centering\psfig{file=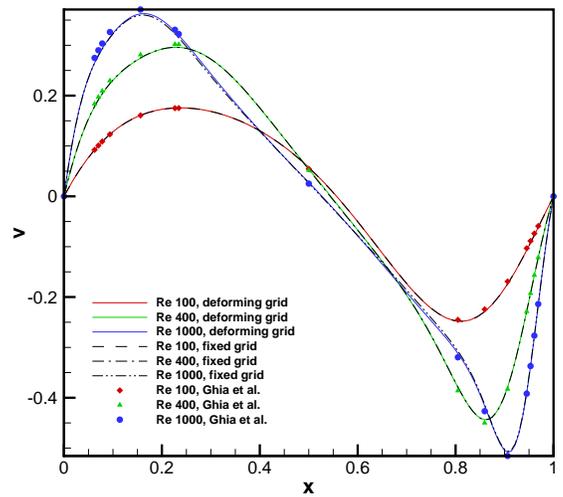,width=0.9\linewidth}
 \\(b)
\end{minipage}
\begin{minipage}[b]{.55\linewidth}\hspace{-1cm}
\centering\psfig{file=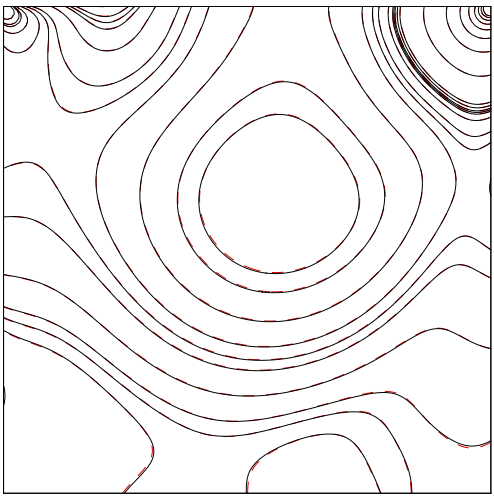,width=0.9\linewidth}
 \\(c)
\end{minipage}
\begin{minipage}[b]{.55\linewidth}\hspace{-1cm}
\centering\psfig{file=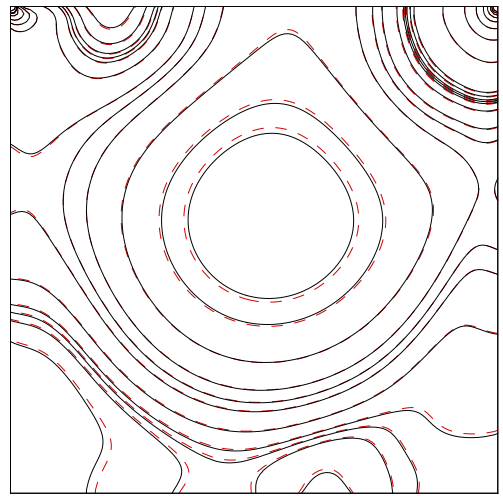,width=0.9\linewidth}
 \\(d)
\end{minipage}
\begin{center}
\caption{{\sl Problem 4: Comparison of velocity profile with benchmark for (a) $u$-component along vertical centreline, (b) $v$-component along horizontal centreline. Pressure contours computed with deforming grid (solid line) and static grid (dashed line) for (c) $Re=400$, (d) $Re=1000$.} }
    \label{fig:LDC_com1}
\end{center}
\end{figure}

We carried out computations on static uniform Cartesian mesh apart from the continually distorting curvilinear grid for $Re=$ 100, 400 and 1000. Time marched steady state solutions are compared to test the ability of the scheme to yield accurate solutions as the underlying grid alternate. Figs. \ref{fig:LDC_com1}(a) and \ref{fig:LDC_com1}(b) illustrate the horizontal and vertical velocity profile along the vertical and horizontal center lines respectively. The velocity field obtained on these two different grid systems are indistinguishable from each other. Further, the profiles compare well to the benchmark result of Ghia \emph{et al.} \cite{ghi_ghi_shi_82}. What is more encouraging, however, is that the equi-pressure contours for $Re=400$ and $1000$ as seen in Figs. \ref{fig:LDC_com1}(c) and \ref{fig:LDC_com1}(d) respectively, computed on deforming and static $129\times129$ grids are very close to each other.

\begin{figure}[htbp]
\begin{minipage}[b]{.55\linewidth}\hspace{-1cm}
\centering\psfig{file=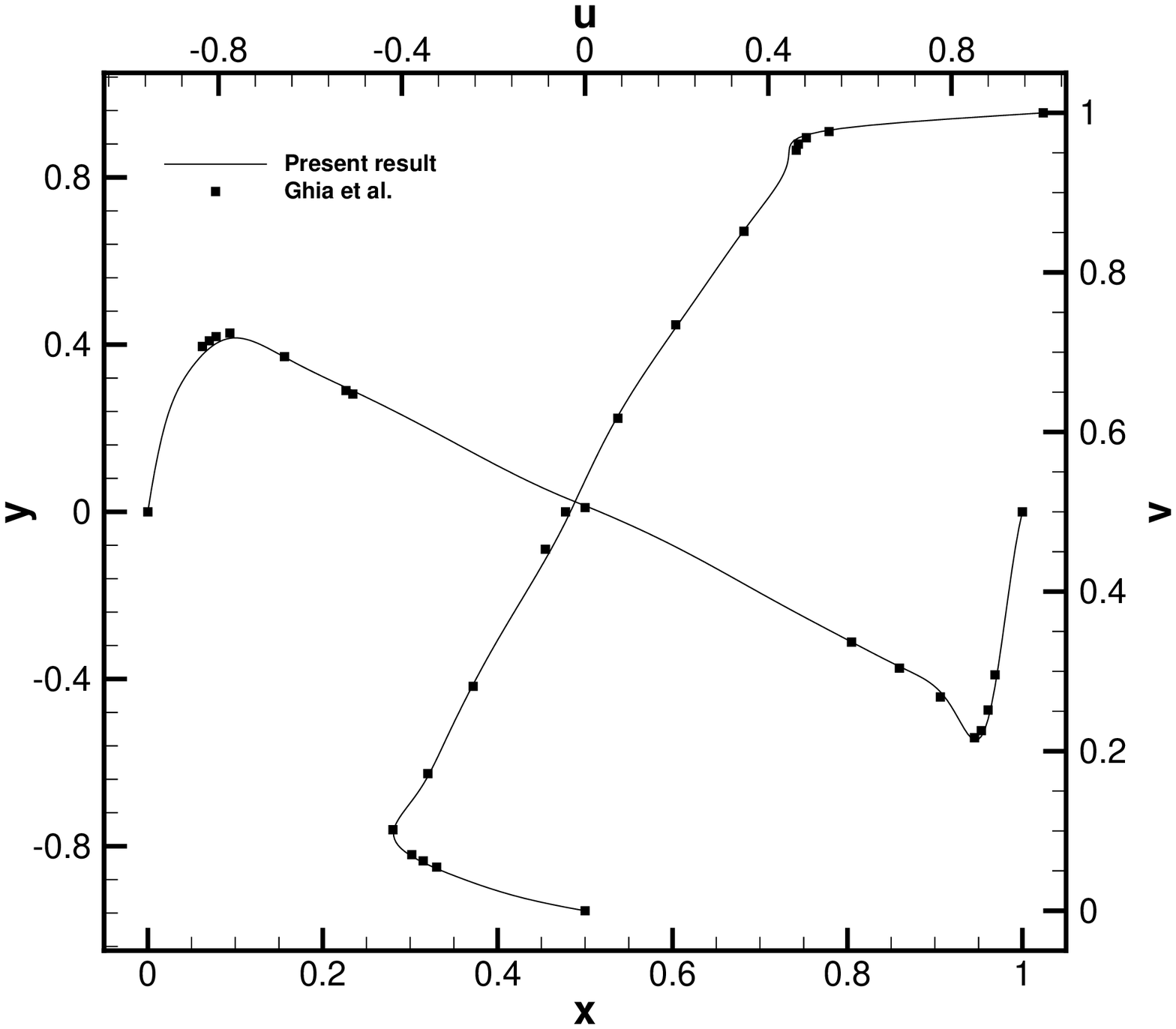,width=0.9\linewidth}
 \\(a)
\end{minipage}
\begin{minipage}[b]{.55\linewidth}\hspace{-1cm}
\centering\psfig{file=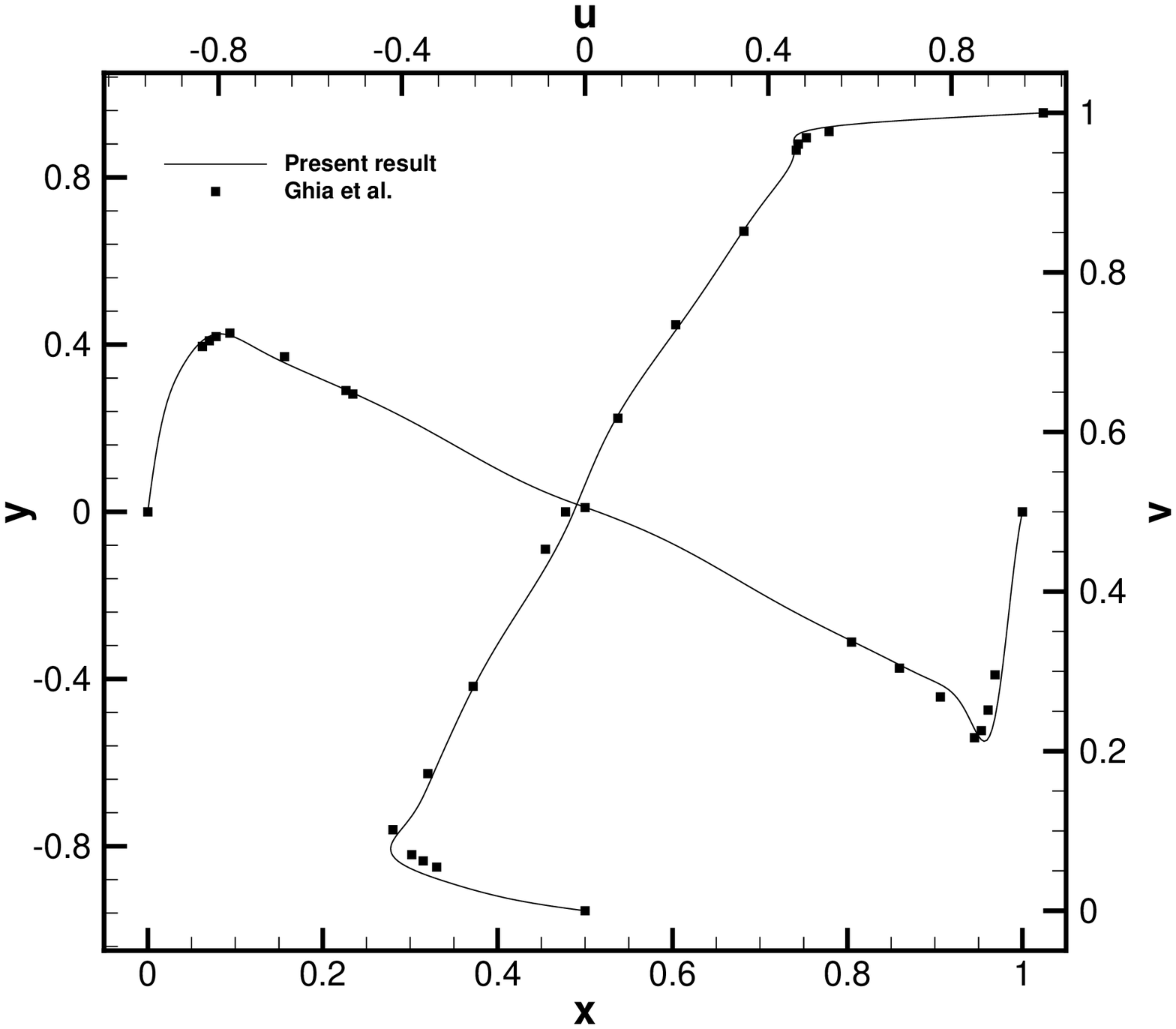,width=0.9\linewidth}
 \\(b)
\end{minipage}
\begin{minipage}[b]{.55\linewidth}\hspace{-1cm}
\centering\psfig{file=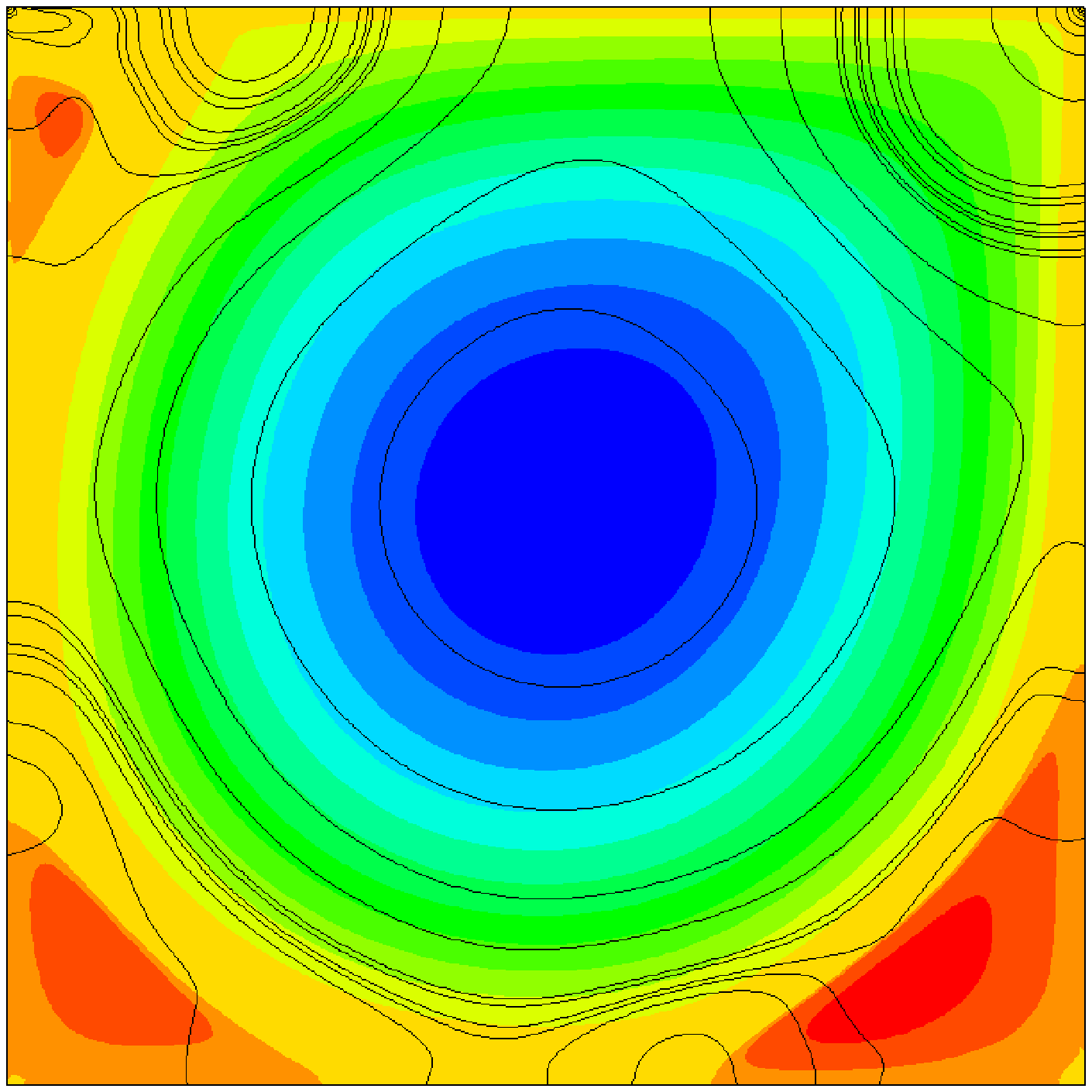,width=0.9\linewidth}
 \\(c)
\end{minipage}
\begin{minipage}[b]{.55\linewidth}\hspace{-1cm}
\centering\psfig{file=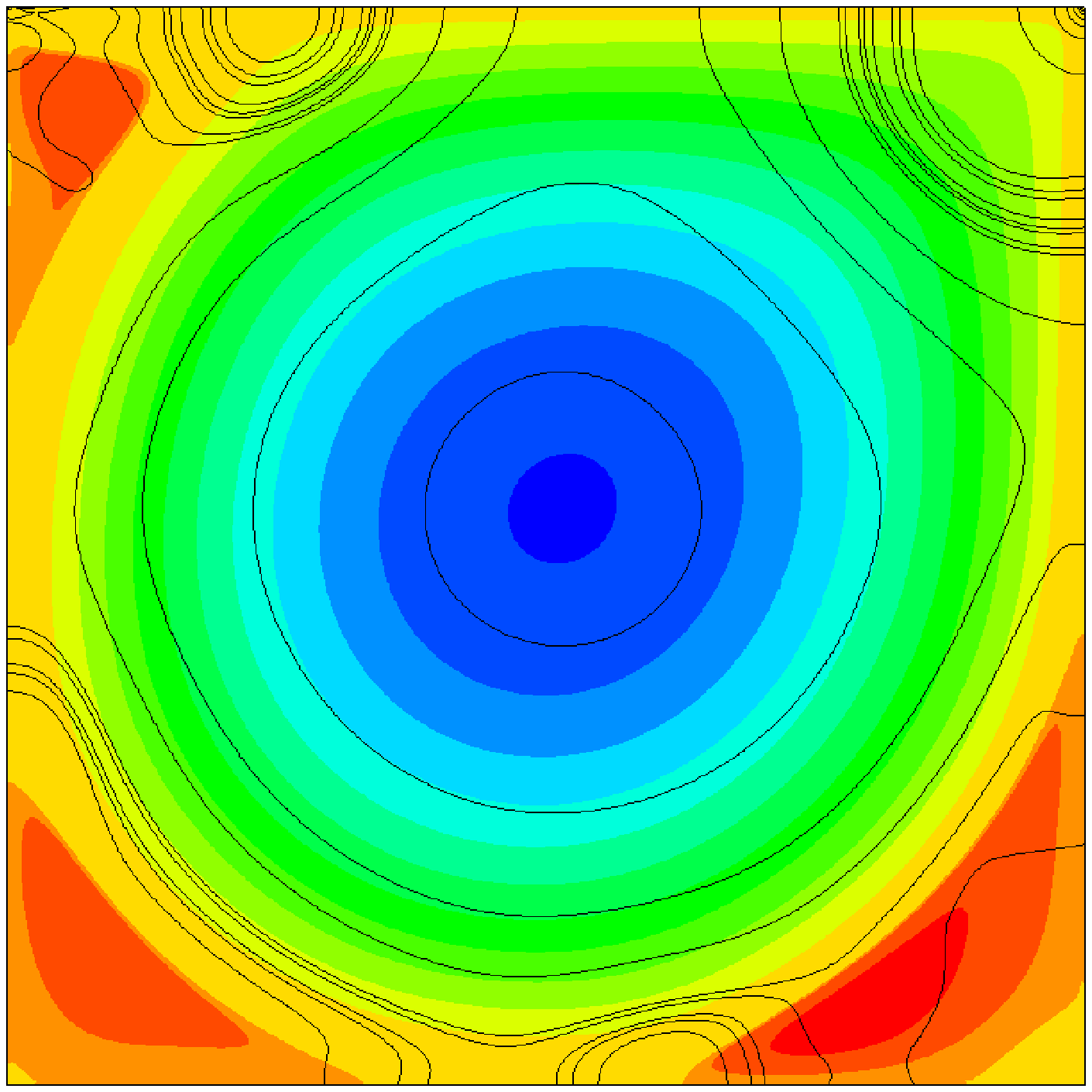,width=0.9\linewidth}
 \\(d)
\end{minipage}
\vspace{-1cm}
\begin{center}
\caption{{\sl Problem 4: Profiles of horizontal and vertical velocities along vertical and horizontal lines respectively through geometric center of cavity for (a) $Re=3200$, (b) $Re=5000$. Plots of pressure contour and
post processed streamlines for (c) $Re=3200$, (d) $Re=5000$.} }
    \label{fig:LDC_com2}
\end{center}
\end{figure}

We also present results for $Re=3200$ and $5000$ computed using $257\times257$ and $321\times321$ time wavering grids respectively in Fig. \ref{fig:LDC_com2}. All simulations are carried out with $\delta \tau=0.0025$. Figs. \ref{fig:LDC_com2}(a) and \ref{fig:LDC_com2}(b), which illustrate the velocity profile along the vertical and horizontal center lines for $Re=3200$ and $5000$ respectively, the velocity field obtained on deforming grid again match well with the values provided in \cite{ghi_ghi_shi_82}. Fig. \ref{fig:LDC_com2}(c) depicts a snapshot of the flow field for $Re=3200$. Based on the contour plots of the pressure and the streamlines, the separated flow regions are correctly predicted and compares well with the results available in the literature. Similarly, in Fig. \ref{fig:LDC_com2}(d) computed flow field and pressure contours at steady state for $Re=5000$ are correctly captured. 

\subsection{Deforming cavity problem}
\begin{figure}[htbp]
\begin{minipage}[b]{.5\linewidth}\hspace{-1cm}
\centering\psfig{file=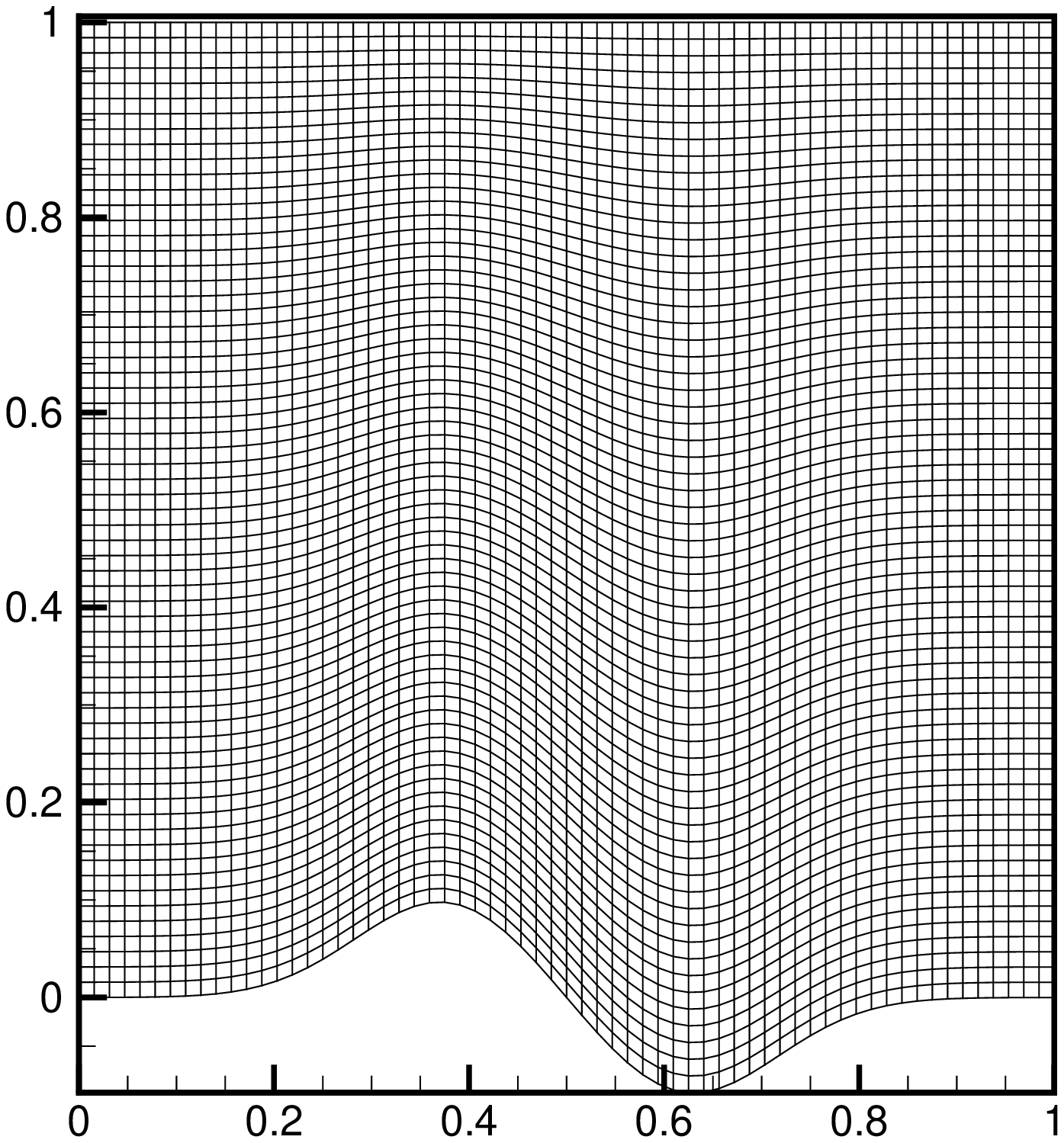,width=0.9\linewidth}(a)
\end{minipage}
\begin{minipage}[b]{.5\linewidth}\hspace{-1cm}
\centering\psfig{file=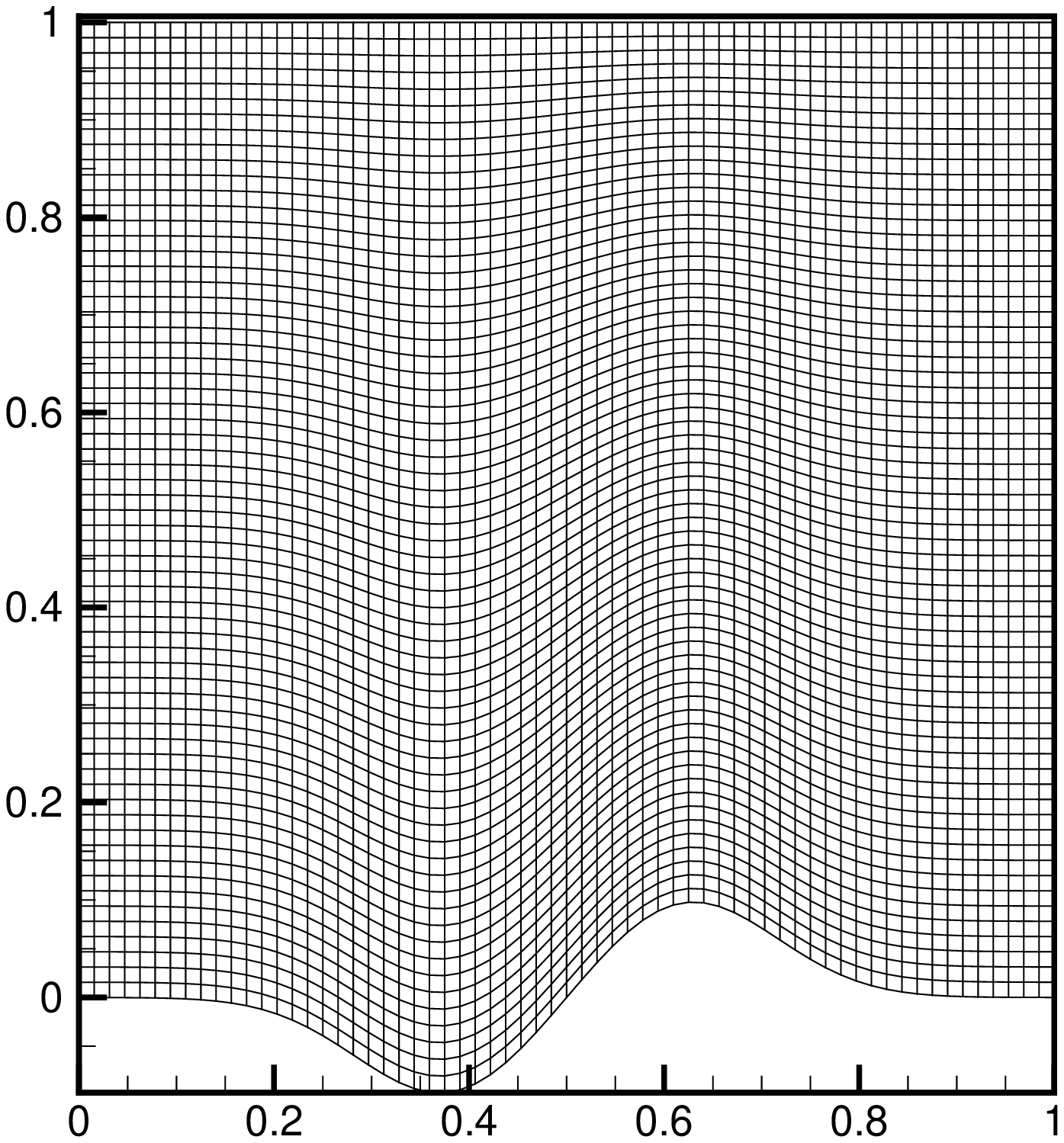,width=0.9\linewidth}(b)
\end{minipage}
\begin{center}
\caption{{\sl Problem 5: Configuration of the deformable cavity and $65\times65$ mesh at time (a) $T/4$ and (b) $3T/4$.} }
    \label{fig:ldc_grid_def}
\end{center}
\end{figure}

We next simulate incompressible flow inside a closed but deforming boundary. This specially designed case is used to investigate our approach on the flows with the deformation of the boundary. The difference between this case and the previous lid-driven cavity problem is that the lower wall of the cavity is deforming and the domain is meshed using the prescribed function. Here we adopt the work of Chen and Xie \cite{che_xie_16} and stipulate contortion of the lower boundary of the cavity as
\begin{equation}\label{ldc_def}
y(x,t)=A\left(e^{a(x-c_1)^2}-e^{a(x-c_2)^2}\right)\sin(2\pi f t),
\end{equation}
where following \cite{che_xie_16} we take amplitude and frequency of wall deformation as $A=0.1$ and $f=0.2$ respectively. Further, positions of maximum deformations are taken symmetrically at $c_1=0.375$ and $c_2=0.625$ with constant for range of deformation $a=60$. The curvilinear coordinates and the geometrical configuration of the physical domain at nondimensional time $T/4$ and $3T/4$ where $T$ is the period of domain alteration is presented in Figs. \ref{fig:ldc_grid_def}(a) and \ref{fig:ldc_grid_def}(b) respectively.

\begin{figure}[htbp]
\begin{minipage}[b]{.5\linewidth}\hspace{-1cm}
\centering\psfig{file=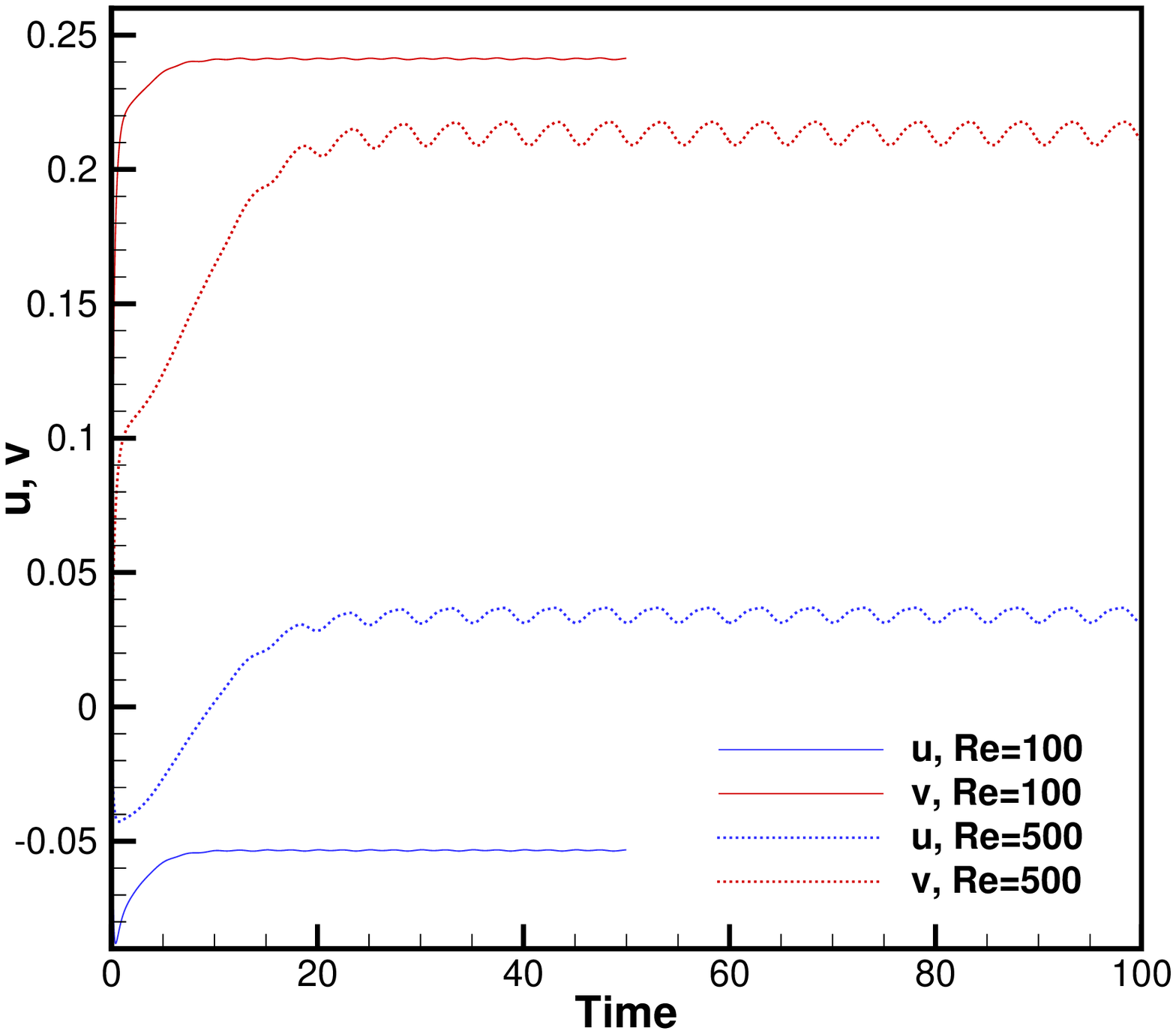,width=0.9\linewidth}(a)
\end{minipage}
\begin{minipage}[b]{.5\linewidth}\hspace{-1cm}
\centering\psfig{file=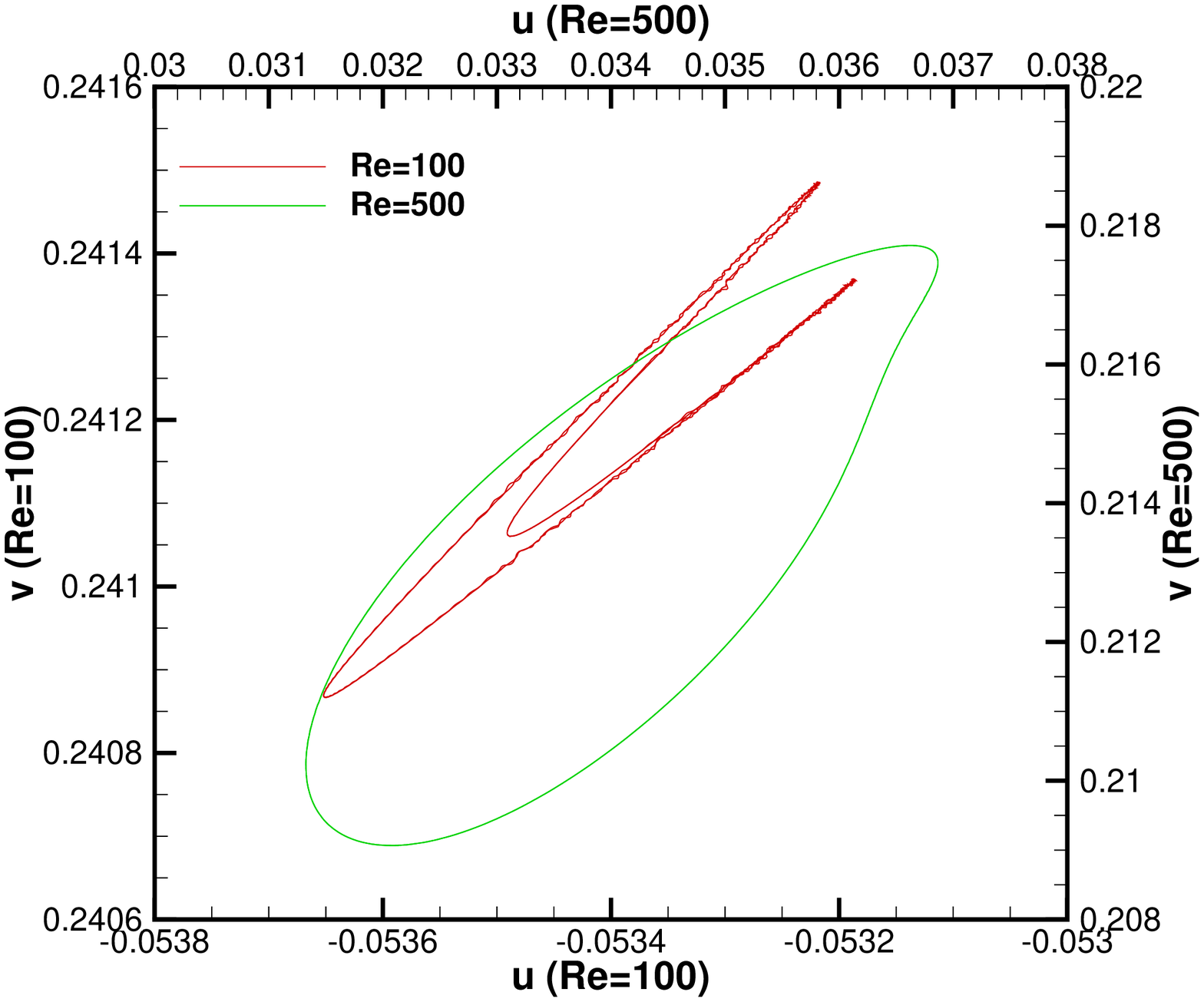,width=0.9\linewidth}(b)
\end{minipage}
\begin{center}
\caption{{\sl Problem 5: (a) Time evolution of velocity components and (b) phase portrait at monitoring point $(2/16, 13/16)$.} }
    \label{fig:ldc_griddef_timeevo}
\end{center}
\end{figure}

In the parametric domain, the flow field is computed for $Re=100$ and $500$ on $65\times65$ and $129\times129$ grids respectively with $\delta\tau=0.0025$. For both these $Re$ values steady-state is reached for lid-driven cavity problem in square cavity discussed in the previous section. But in this section with the domain constantly deforming following a cyclical function in time it is impossible to arrive at steady motion \cite{che_xie_16}. We thus inquire about the periodic solution and monitor the time evolution of velocity at a chosen point $(2/16, 13/16)$. The choice of monitoring point away from the deforming lower boundary is partly influenced by the work of Bruneau and Saad \cite{bru_saa_06} which was used to analyze the periodic solution for classical driven cavity problem. Time development of horizontal and vertical velocity components at the monitoring point is presented in Fig. \ref{fig:ldc_griddef_timeevo}(a). It is seen that the flow away from the lower boundary after starting from rest settles to periodic variation. This is emphasized for $Re=500$ whereas for $Re=100$ variations are much subdued and is captured by phase portrait in Fig. \ref{fig:ldc_griddef_timeevo}(b). 

\begin{figure}[htbp]
\begin{minipage}[b]{.33\linewidth}\hspace{-1cm}
\centering\psfig{file=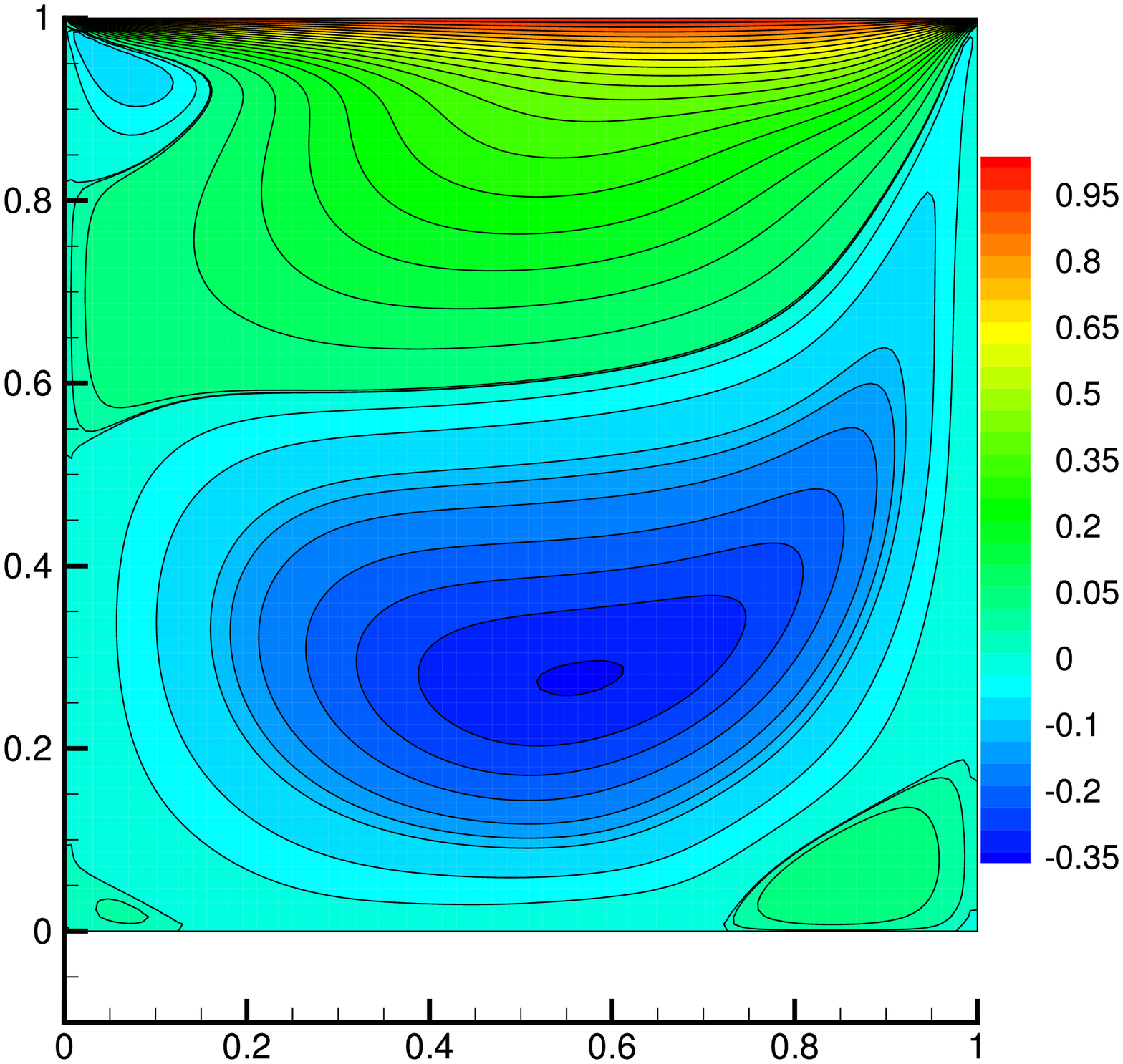,width=0.8\linewidth}
 (a)
\end{minipage}
\begin{minipage}[b]{.33\linewidth}\hspace{-1cm}
\centering\psfig{file=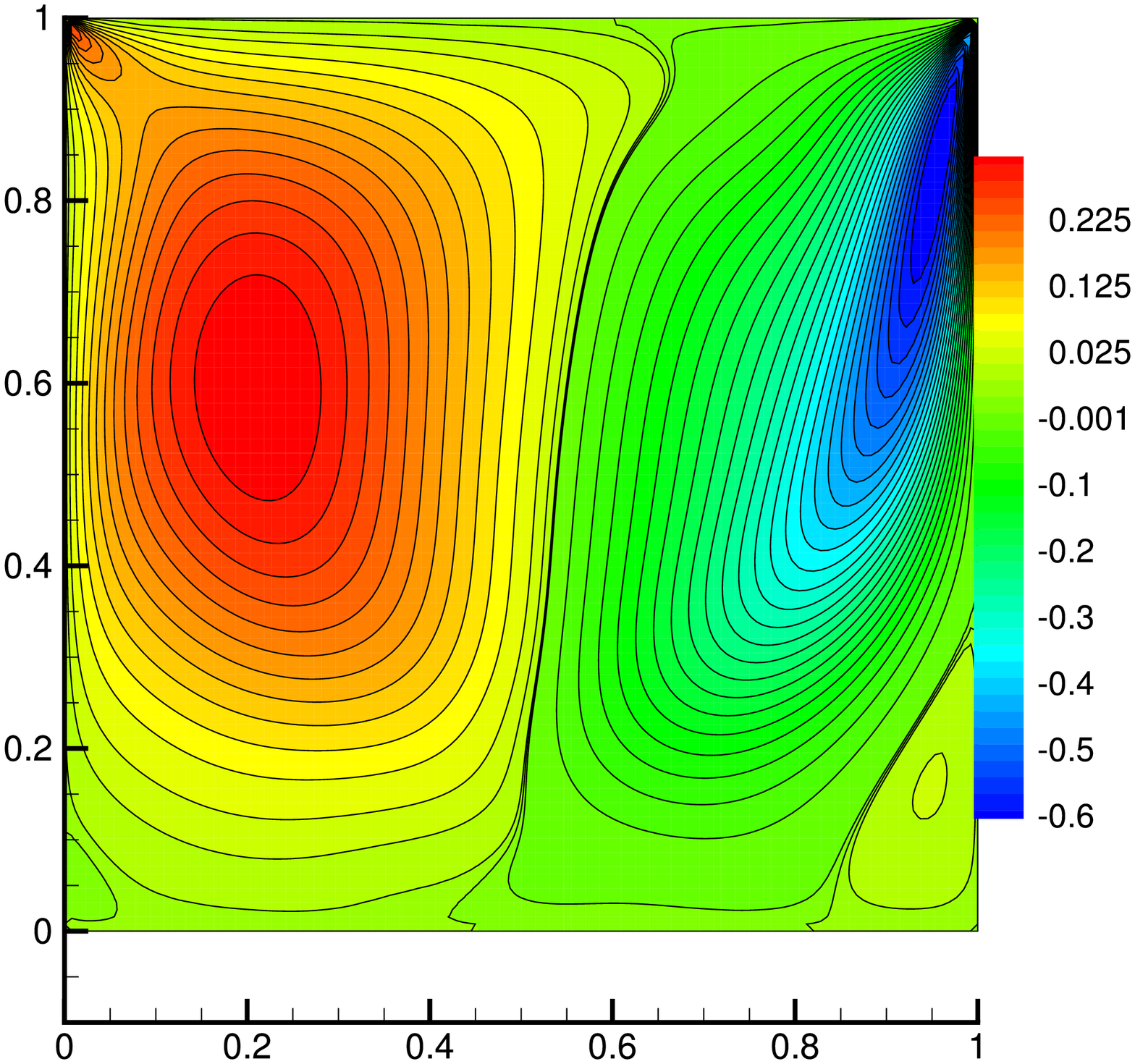,width=0.8\linewidth}
 (b)
\end{minipage}
\begin{minipage}[b]{.33\linewidth}\hspace{-1cm}
\centering\psfig{file=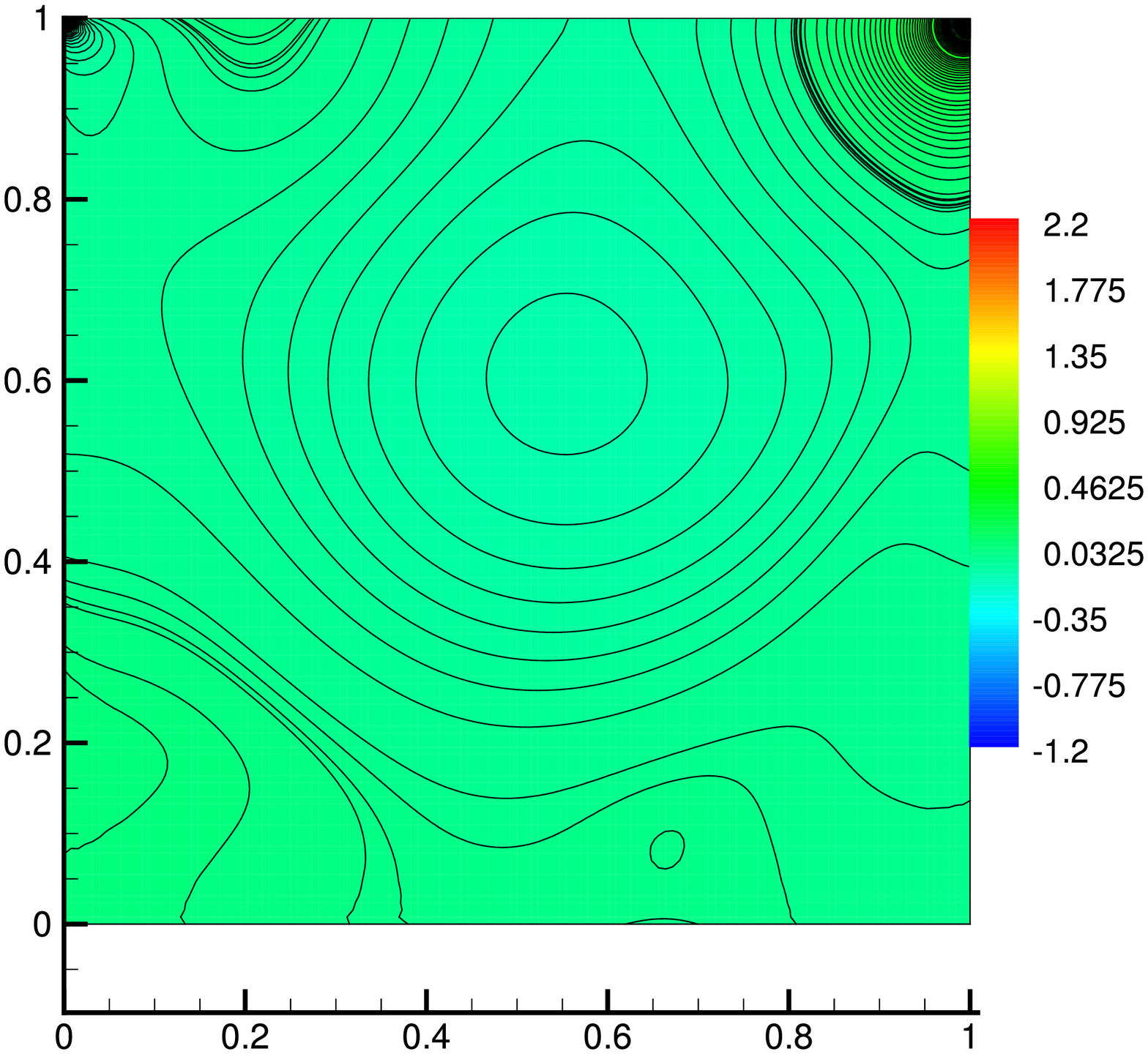,width=0.8\linewidth}
 (c)
\end{minipage}
\begin{minipage}[b]{.33\linewidth}\hspace{-1cm}
\centering\psfig{file=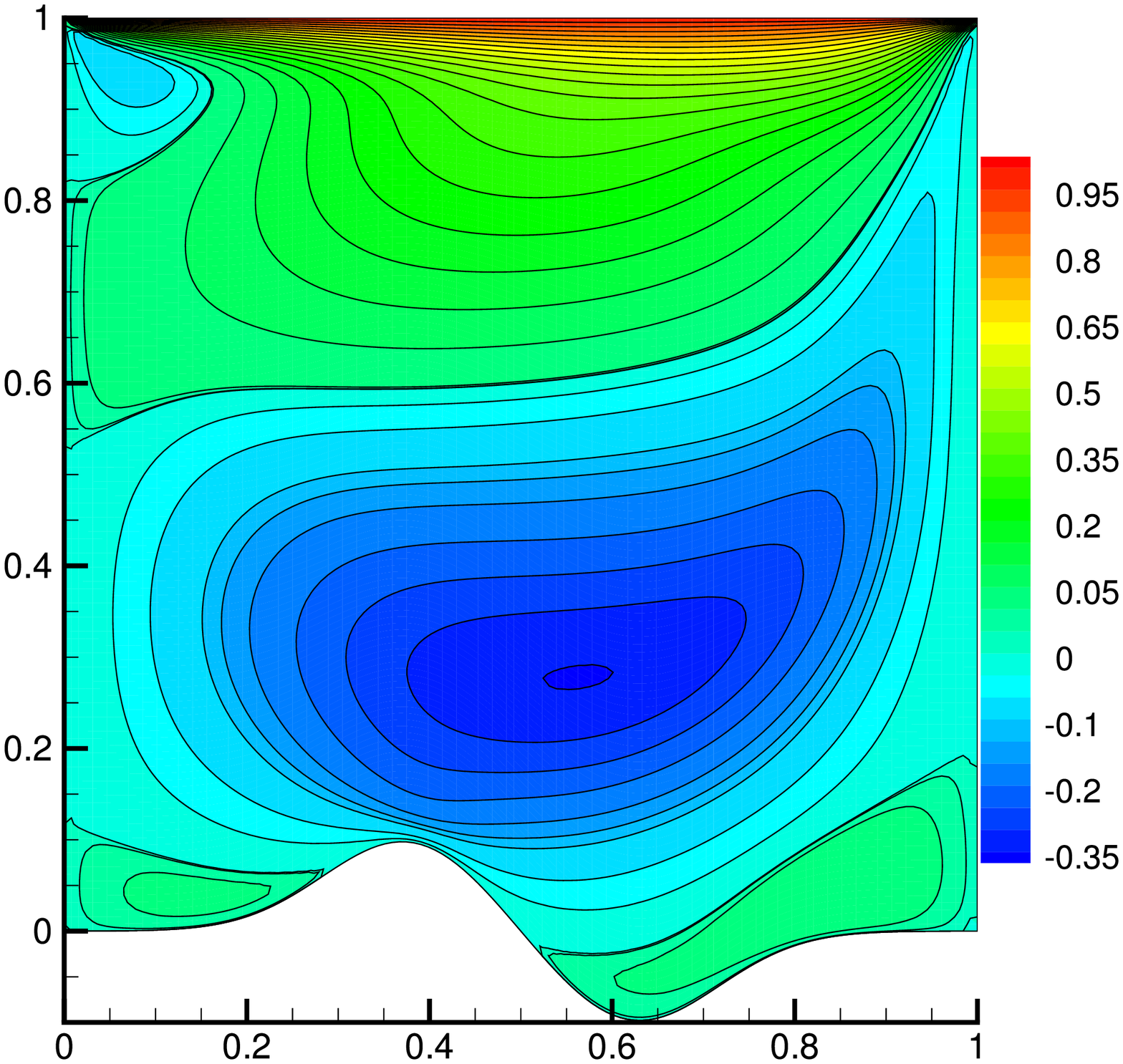,width=0.8\linewidth}
 (d)
\end{minipage}
\begin{minipage}[b]{.33\linewidth}\hspace{-1cm}
\centering\psfig{file=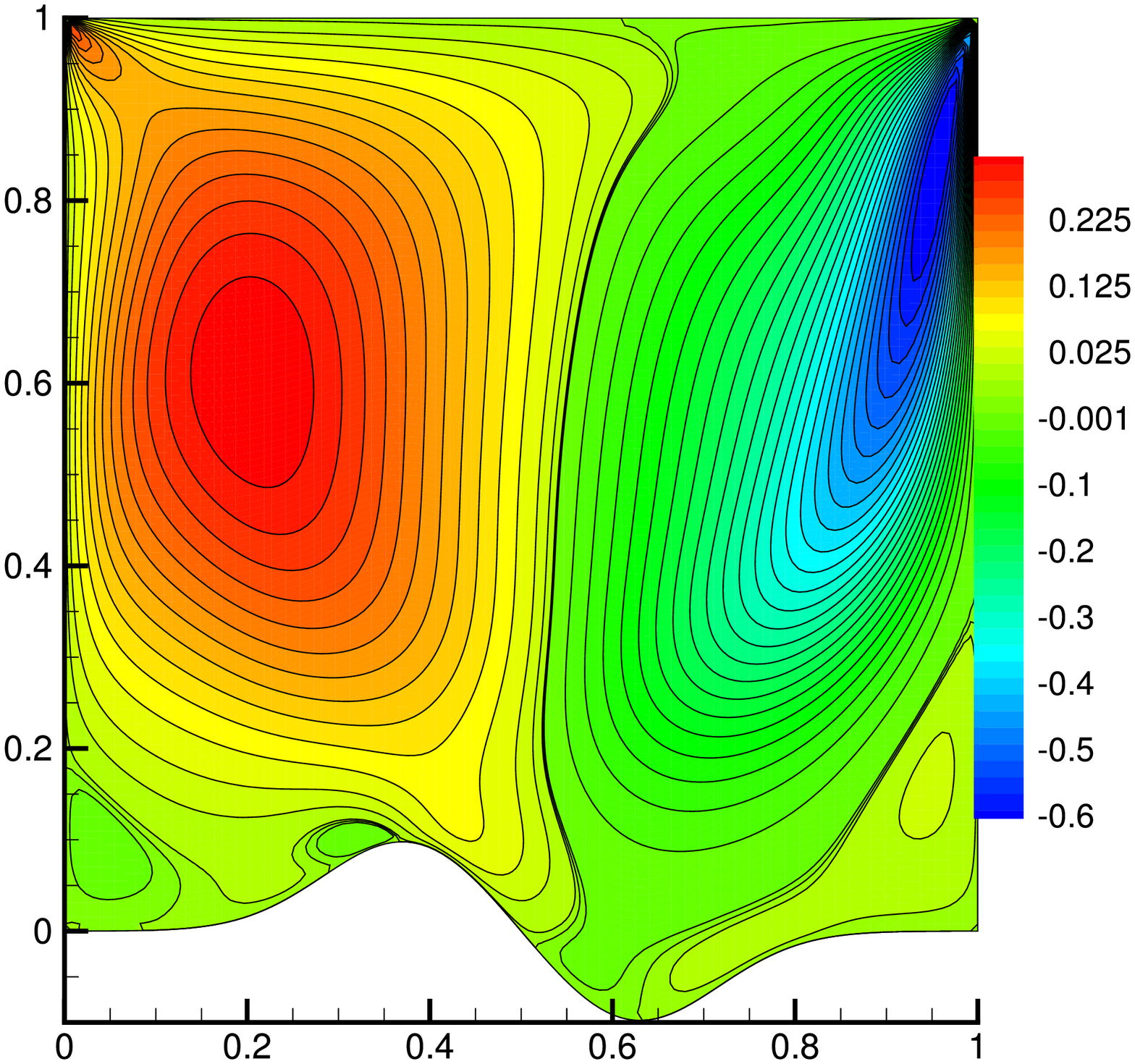,width=0.8\linewidth}
 (e)
\end{minipage}
\begin{minipage}[b]{.33\linewidth}\hspace{-1cm}
\centering\psfig{file=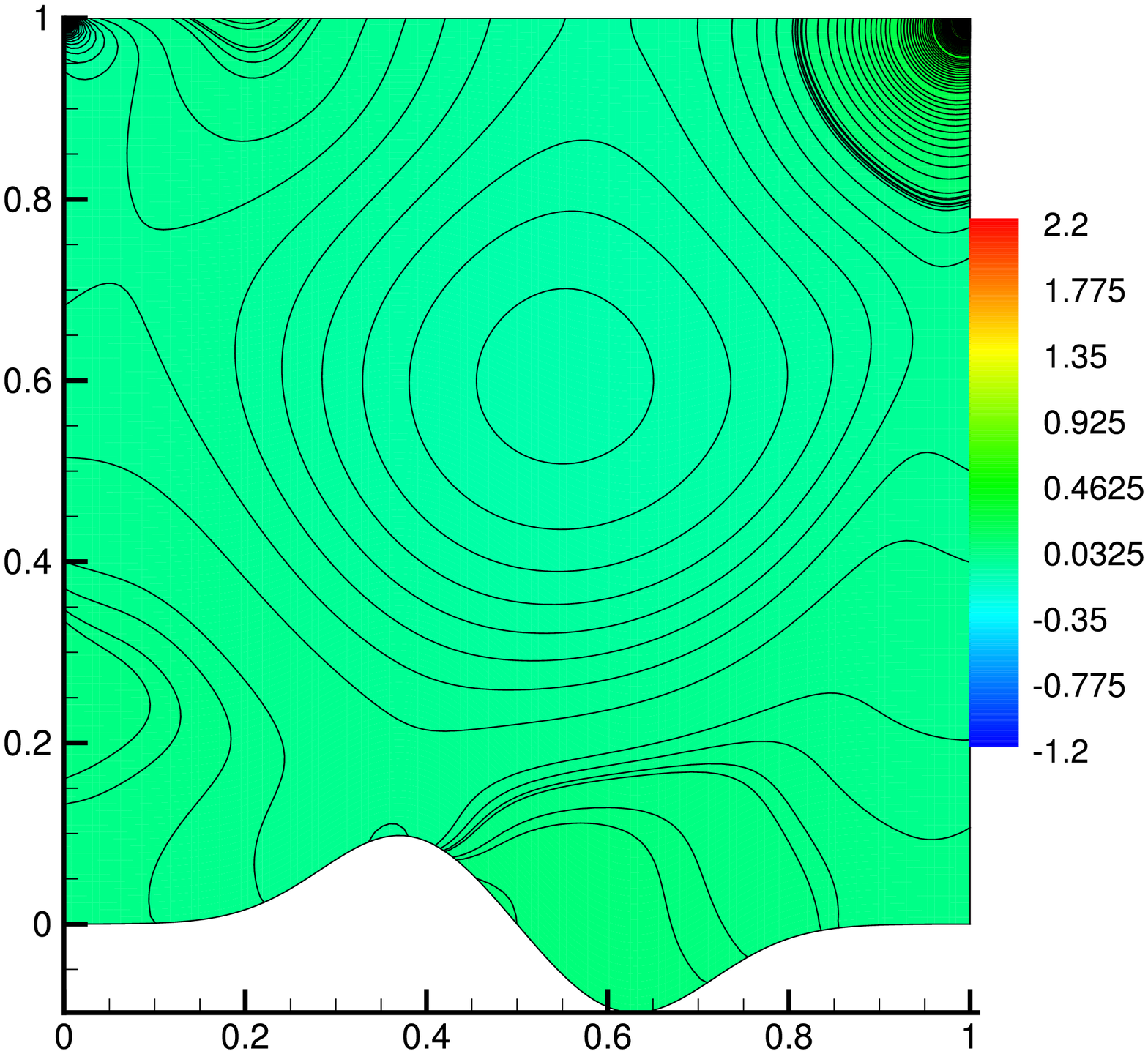,width=0.8\linewidth}
 (f)
\end{minipage}
\begin{minipage}[b]{.33\linewidth}\hspace{-1cm}
\centering\psfig{file=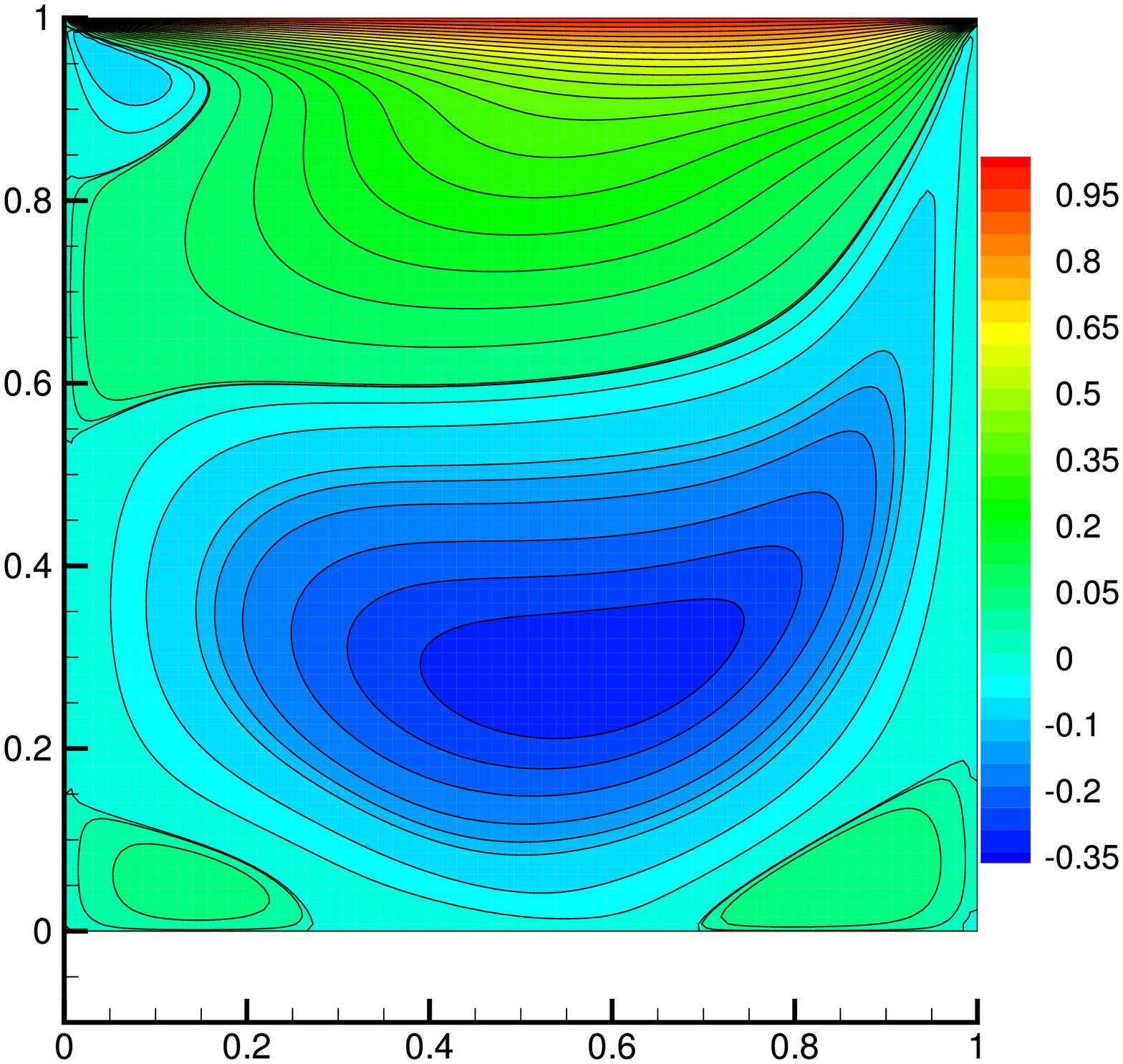,width=0.8\linewidth}
 (g)
\end{minipage}
\begin{minipage}[b]{.33\linewidth}\hspace{-1cm}
\centering\psfig{file=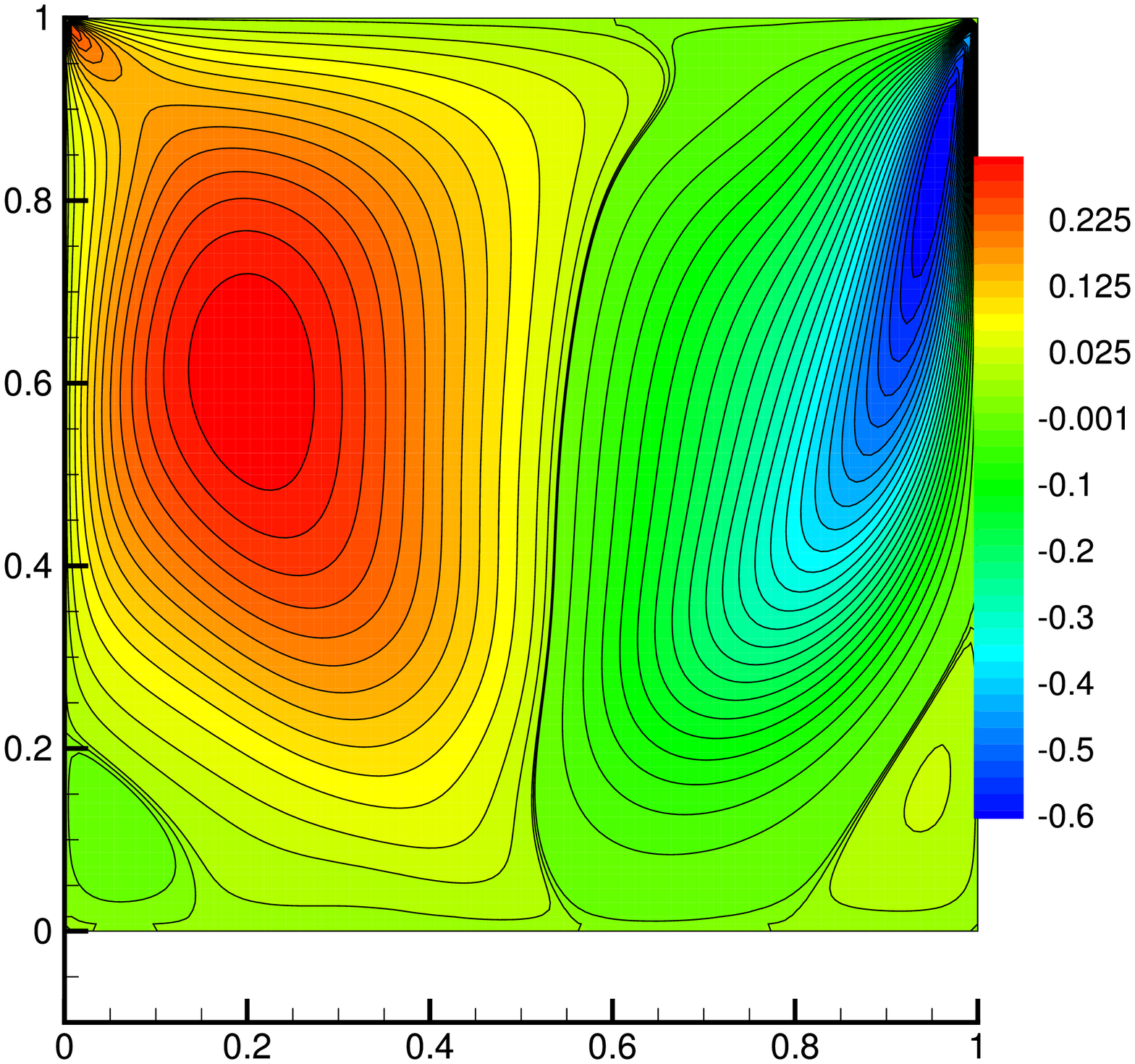,width=0.8\linewidth}
 (h)
\end{minipage}
\begin{minipage}[b]{.33\linewidth}\hspace{-1cm}
\centering\psfig{file=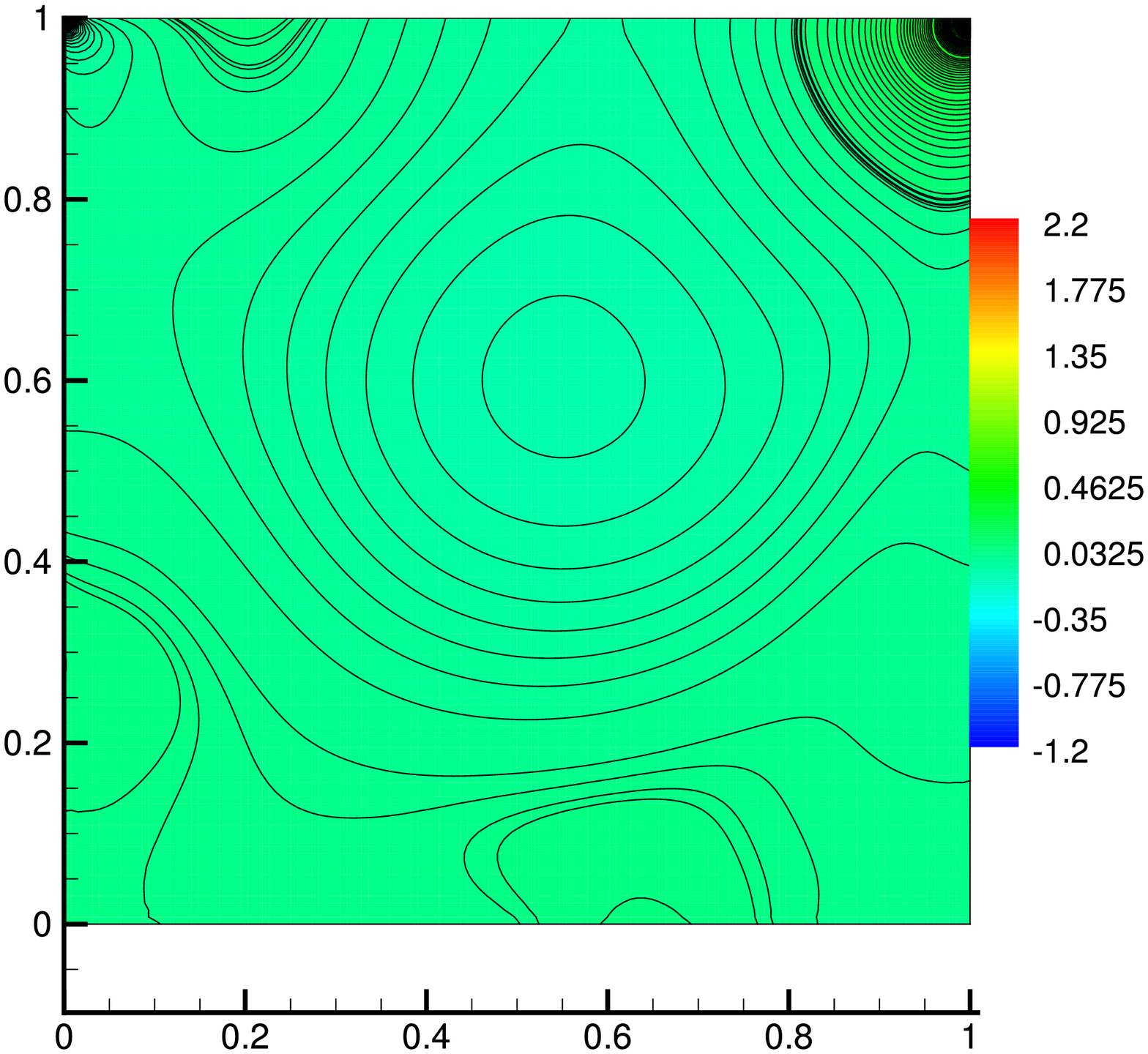,width=0.8\linewidth}
 (i)
\end{minipage}
\begin{minipage}[b]{.33\linewidth}\hspace{-1cm}
\centering\psfig{file=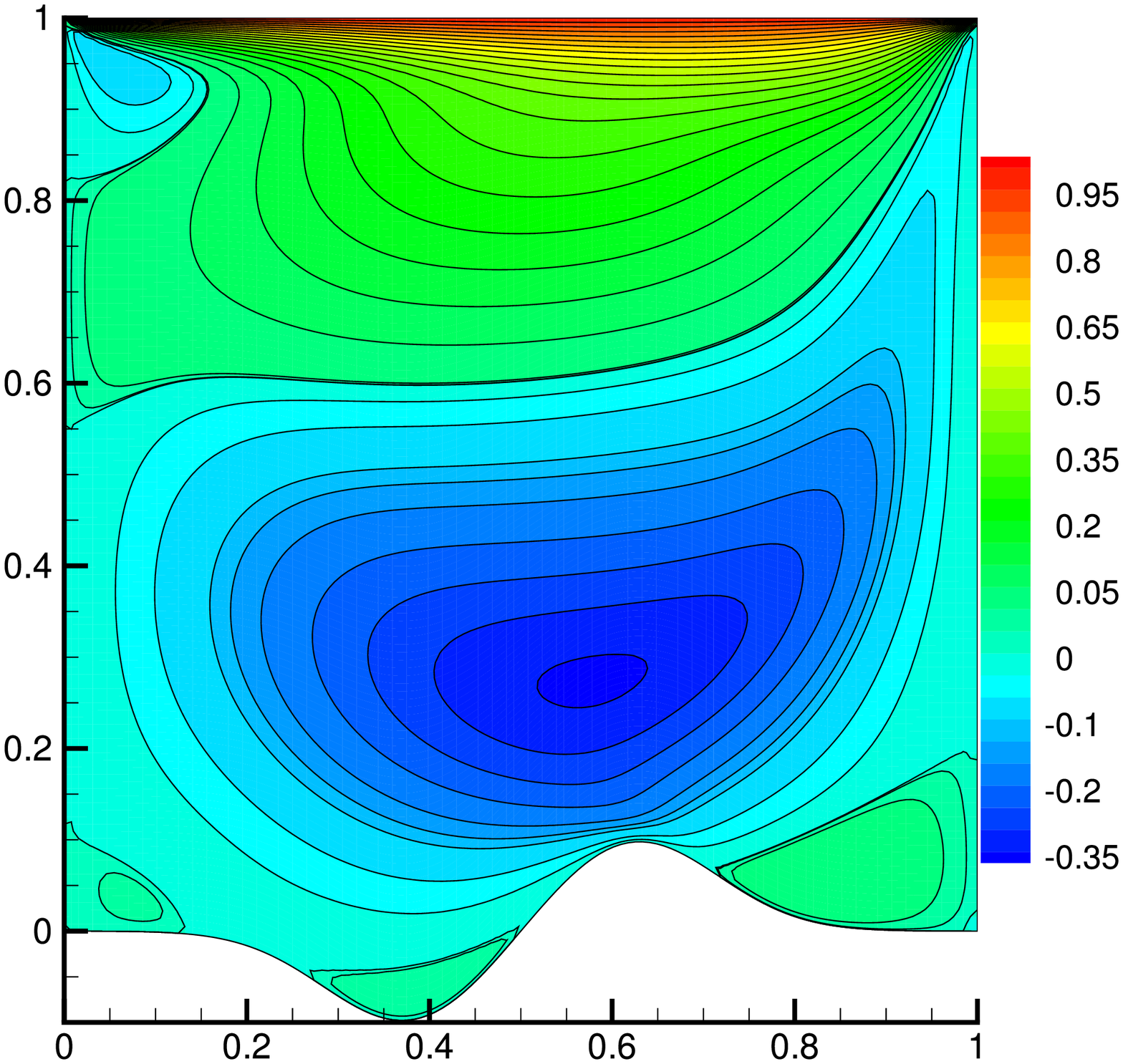,width=0.8\linewidth}
 (j)
\end{minipage}
\begin{minipage}[b]{.33\linewidth}\hspace{-1cm}
\centering\psfig{file=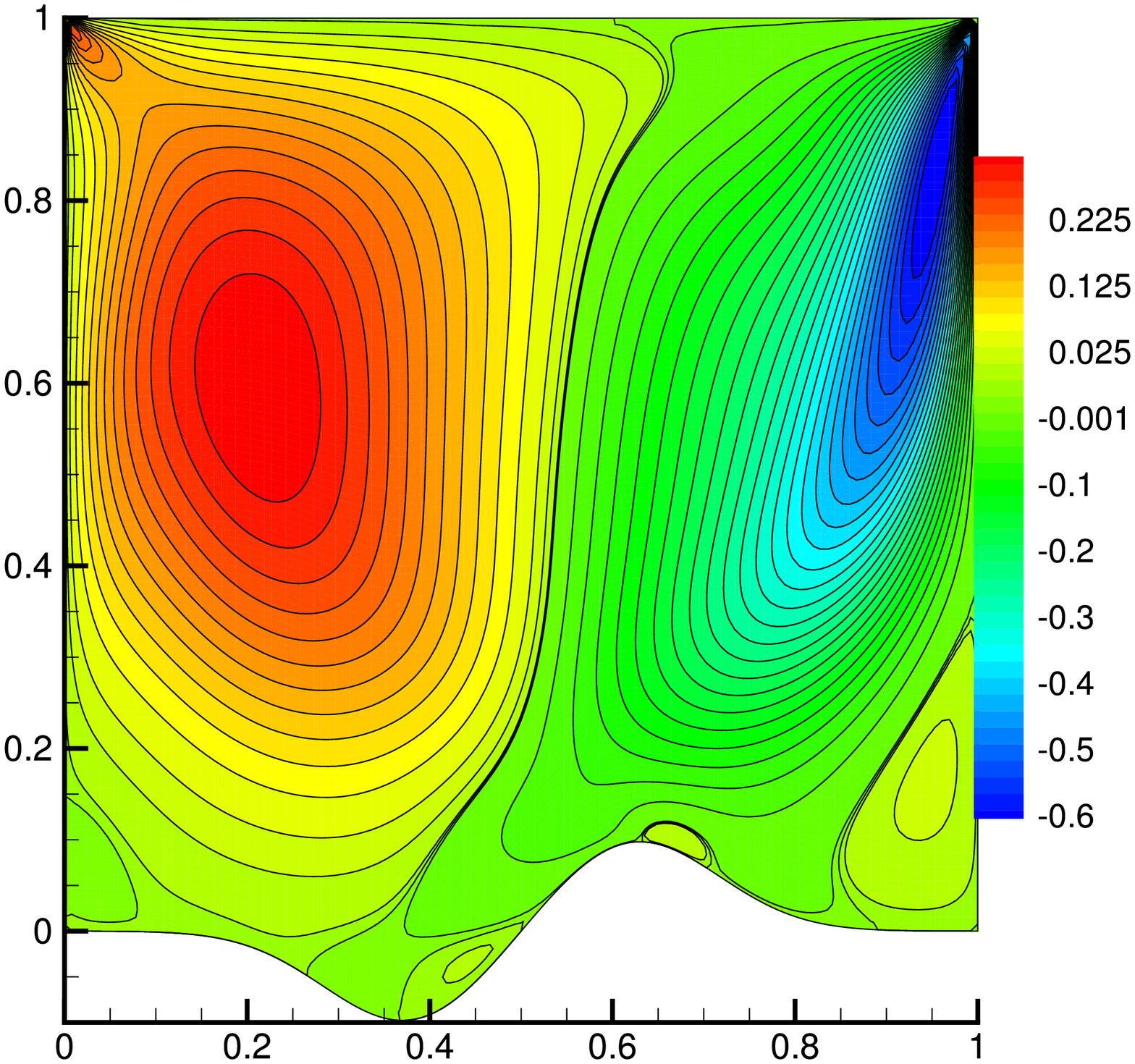,width=0.8\linewidth}
 (k)
\end{minipage}
\begin{minipage}[b]{.33\linewidth}\hspace{-1cm}
\centering\psfig{file=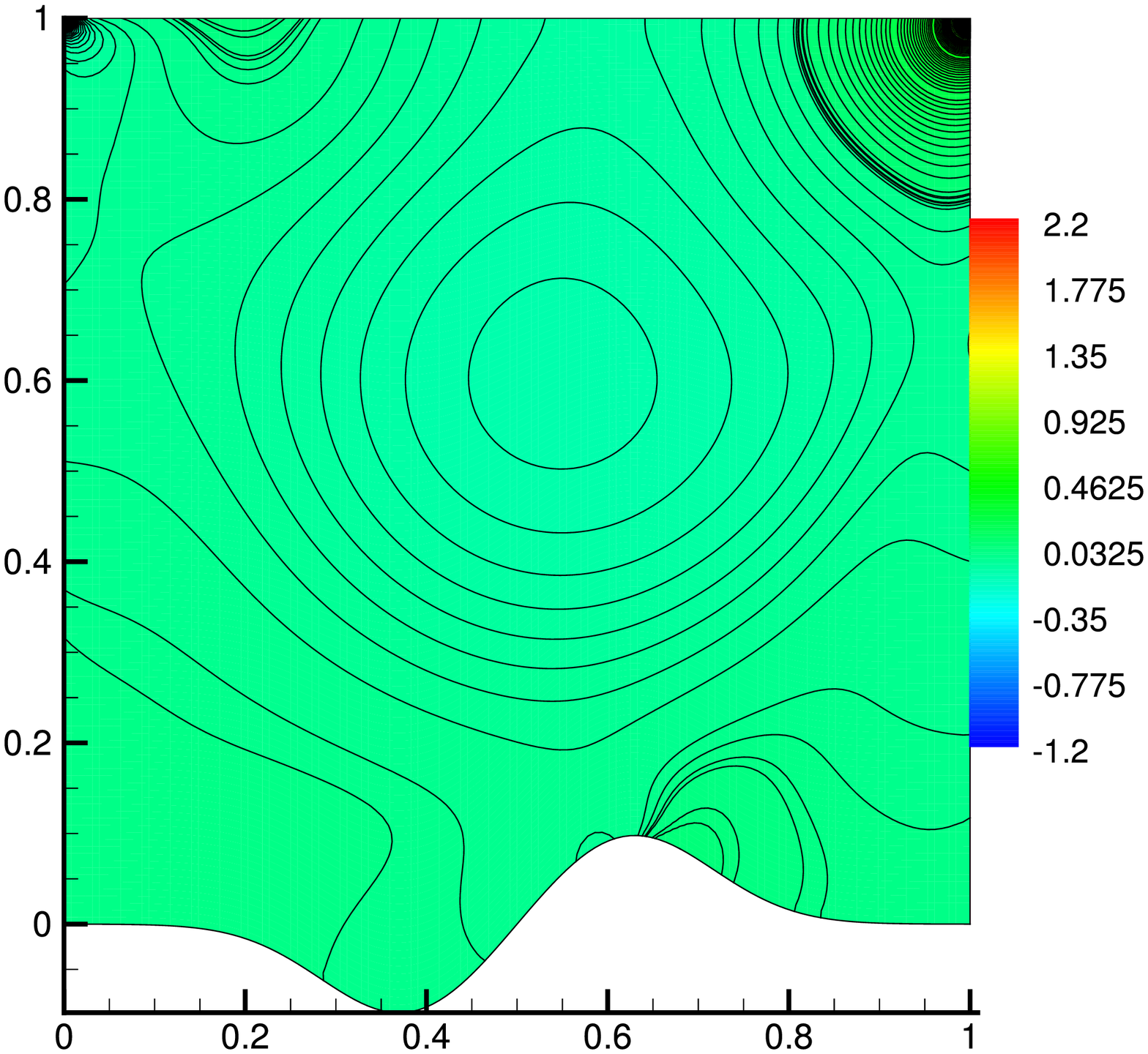,width=0.8\linewidth}
 (l)
\end{minipage}
\begin{center}
\caption{{\sl Problem 5: Evolution of the horizontal velocity (left), vertical velocity (middle) and pressure (right) during one period for $Re=500$. From top to bottom rows correspond to times $0$, $T/4$, $2T/4$ and $3T/4$.} }
    \label{fig:ldc_griddef_contour}
\end{center}
\end{figure}

Fig. \ref{fig:ldc_griddef_contour} shows the instantaneous velocity and pressure field for $Re=500$ at different time phases within one period. From the snapshots, it is clear that although overall periodic, in the region far from the deforming surface the flow is similar to the steady case perceptible in classical cavity problem. From the plot of velocity components, it is amply clear that with changing bottom topography of the cavity secondary vortices on each side of the domain undergoes continuous change in shape and size. Nevertheless, it can be inferred that these vortices remain stable and attached to the lower boundary. The periodic appearance of a third secondary vortex is also noticed from Figs. \ref{fig:ldc_griddef_contour}(j) and \ref{fig:ldc_griddef_contour}(k).
\subsection{Oscillating cylinder in cross flow}
\begin{figure}[htbp]
\begin{minipage}[b]{.5\linewidth}\hspace{-1cm}
\centering\psfig{file=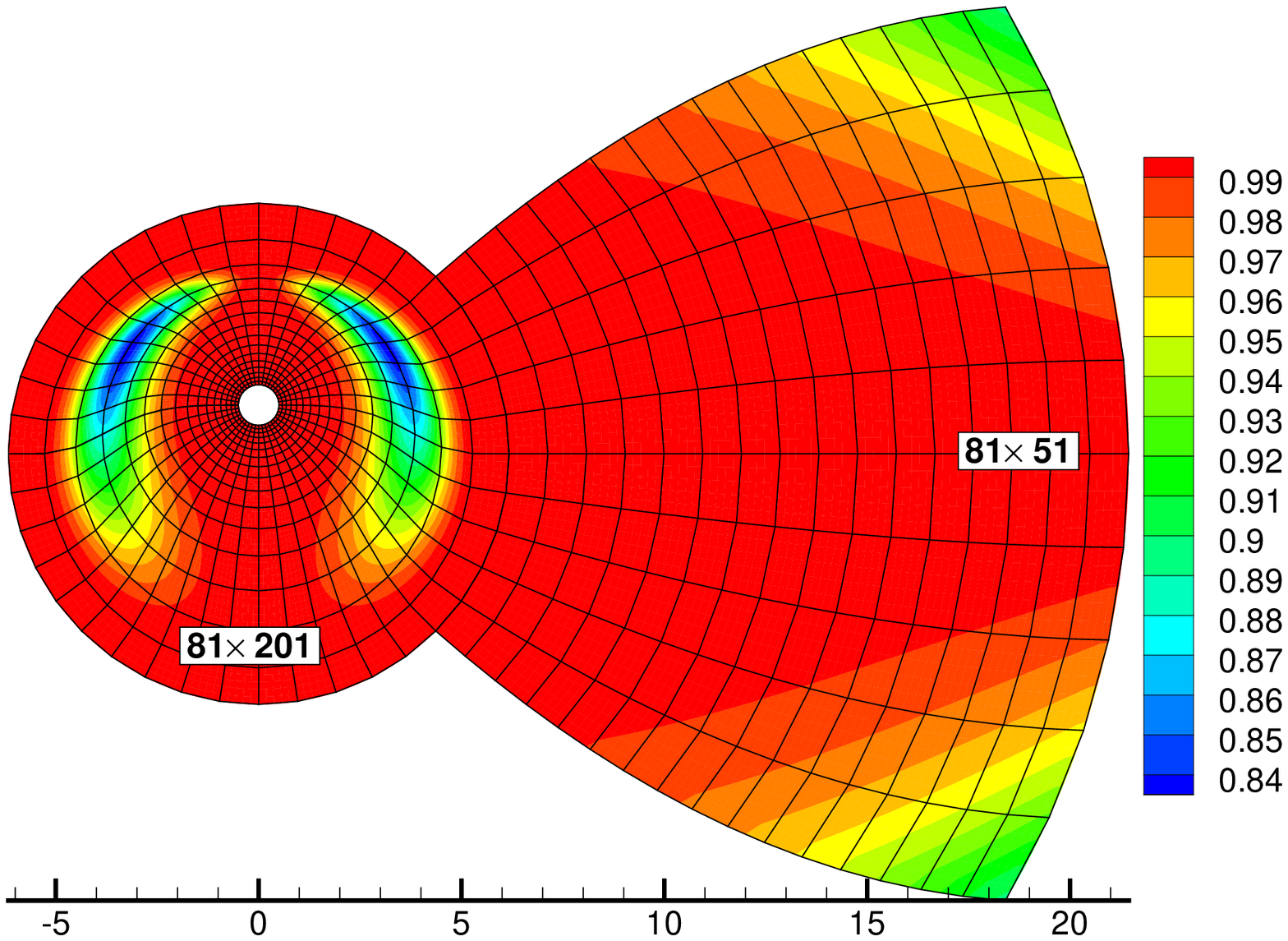,width=0.9\linewidth}(b)
\end{minipage}
\begin{minipage}[b]{.5\linewidth}\hspace{-1cm}
\centering\psfig{file=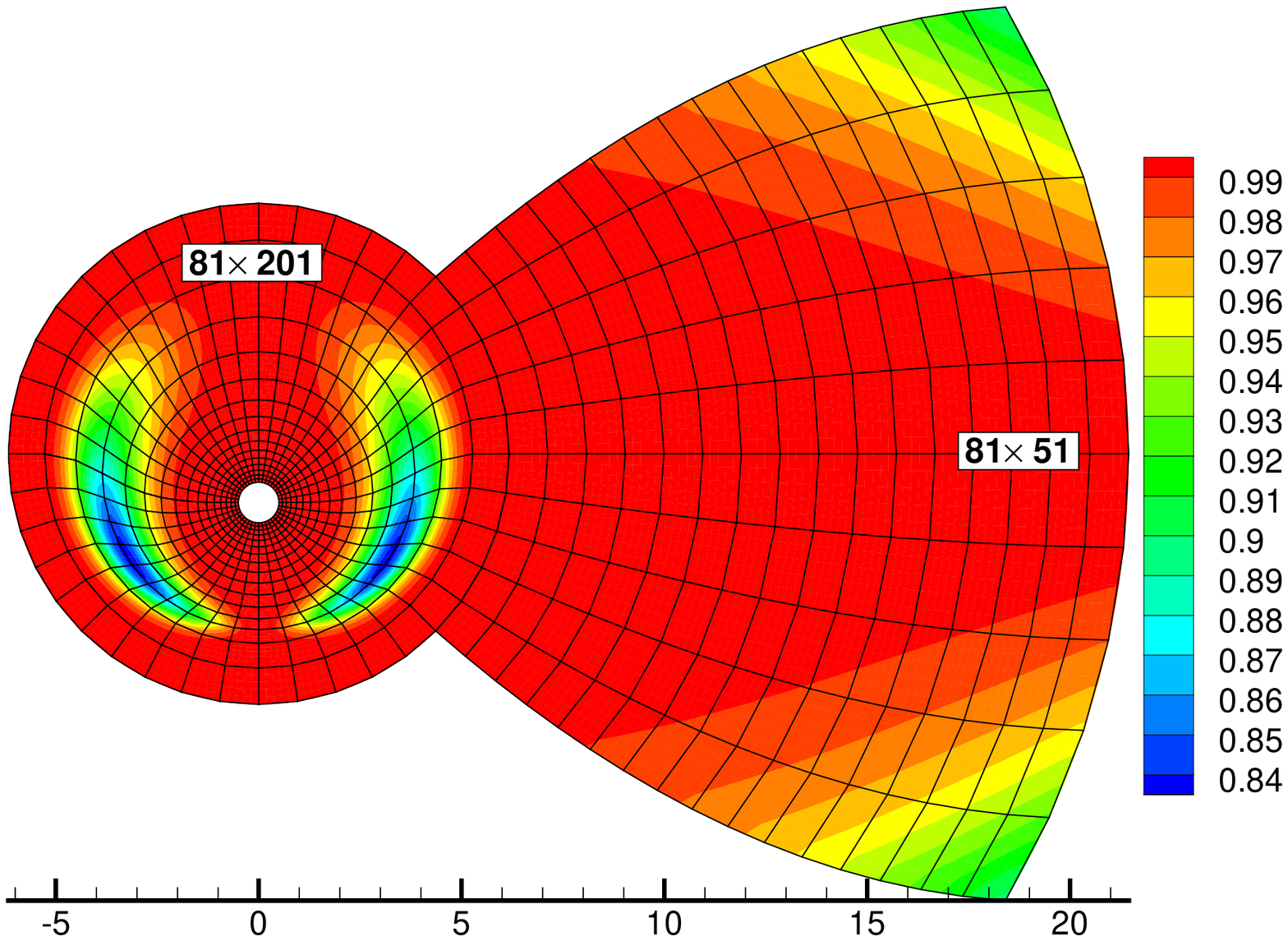,width=0.9\linewidth}(b)
\end{minipage}
\begin{center}
\caption{{\sl Problem 6: Grid along with skew metrics for cylinder positioned (a) up and (b) down 1.2 times diameter from the mean position. Every fifth grid line is shown.} }
    \label{fig:Cyl_cross_grid}
\end{center}
\end{figure}

In this example, numerical simulation of the flow past a circular cylinder that is able to oscillate transversely in a uniform background flow is presented. In a seminal work, Williamson and Roshko in 1988 \cite{wil_ros_88} documented the influence of bluff body oscillation on vortex shedding. They experimentally established several distinct wake structures which they classified in terms of the number of vortices shed per oscillation cycle. Afterward, this involved fluid--structure interaction problem has drawn much attention and is considered a typical test case for numerical techniques. 
\begin{figure}[htbp]
\begin{minipage}[b]{.9\linewidth}\hspace{-1cm}
\centering\psfig{file=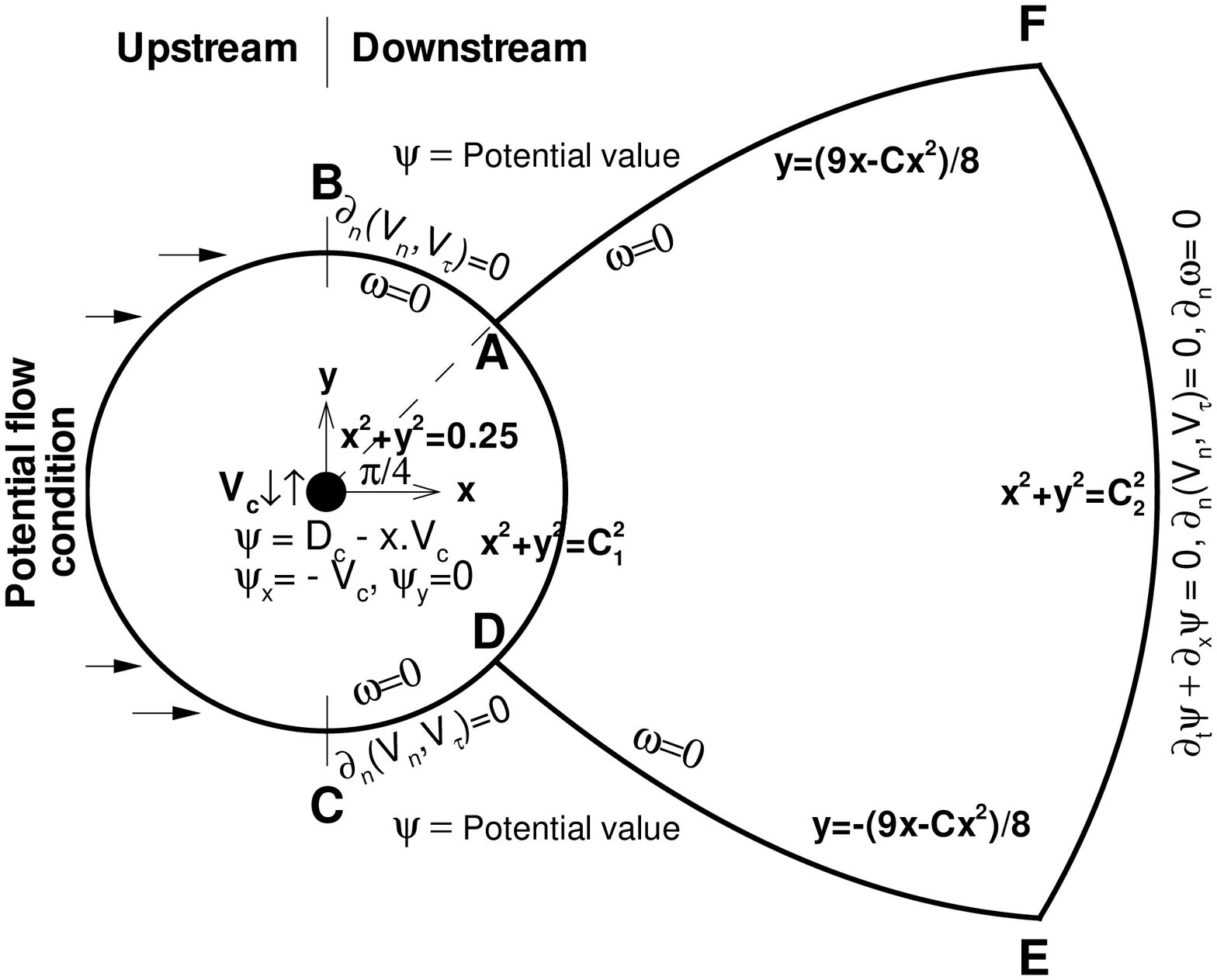,width=0.5\linewidth}
\end{minipage}
\begin{center}
\caption{{\sl Problem 6: Schematic diagram exhibiting domain and boundary conditions for cylinder oscillating in cross-flow.} }
    \label{fig:Cyl_cross_sch}
\end{center}
\end{figure}
In the present study main motivation is to validate some important phenomena hitherto detailed in the literature. In this process, one of the crucial steps is to continually and adequately transform the body fitted grid. As discussed earlier recent advances in IDW interpolation help us generate a good quality time deforming grid. The process starts with an initial grid generated by encompassing the cylinder of unit diameter at origin details of which could be found in the author's earlier work \cite{sen_16}. As shown in schematic diagram Fig. \ref{fig:Cyl_cross_sch} the flow domain consists of two separate blocks joined at the abutting boundary DA. Such a multi-block grid helps us not only to cluster an adequate number of mesh points in the vicinity of the cylinder surface but also to capture flow wake up to a reasonable distance while being economical with the distribution of nodal points. With the computational domain held fixed, adequate capture of transverse oscillation of the cylinder makes it essential to acceptably adjust the grid at each time step. Since the motion of the cylinder is in the transverse direction we resist calculating rotation quaternions and spontaneously engage IDW interpolation methodology. Two representative grids with cylinder displaced 1.2 times its diameter from the mean position at the origin are depicted in Fig. \ref{fig:Cyl_cross_grid}. Here every fifth grid line is shown. For the generated meshes skew metric is an important characteristic \cite{knu_03} and in the current context ensures the quality of grid deformation. Having implemented IDW interpolation we are able to maintain skew value close to unity in large regions of the domain and the same deviate marginally only at a few pockets as seen in Fig. \ref{fig:Cyl_cross_grid}. 

Having discussed grid generation we turn our attention to briefly describe problem formulation and associated boundary conditions. For this test case $\psi-\omega$ formulation is tabbed and the flow commences impulsively. Here a rigid circular cylinder of diameter $D$ is assumed to be elastically mounted. The cylinder is restrained to executes harmonic motion with amplitude $A_c$ and frequency $f_c$ transverse to the incoming flow having free-stream velocity $U$. Computations are done for $Re=200$ and $392$ where $Re=UD/\nu$, $\nu$ being kinematic viscosity. These choices are influenced by the availability of results in the literature. 

It is well known that the wake region of the flow past a stationary circular cylinder exhibit second Hopf bifurcation at some critical $Re$ value ($Re_{cr}$). In the literature, this threshold value is scattered in the range $43\le Re_{cr} \le 50$ details of which could be found elsewhere \cite{kal_sen_12}. Flows above $Re_{cr}$ are characterized by von K\'arm\'an vortex street and exhibits certain periodicities identified by Strouhal number ($St$) and defined as $St=Df_n/U$ where $f_n$ is the frequency of periodic variation of lift exerted on a stationary cylinder. Following Williamson and Roshko \cite{wil_ros_88} we take $St=0.200$. 

In our computations, with cylinder performing forced oscillations, excitation amplitude and frequency are characterized by adimensional amplitude $A_r=A_c/D$ and frequency ratio $f_r=f_c/f_n$. Further oscillatory velocity transverse to the flow is imposed as $V_c=2\pi f_cA_c\cos(2\pi f_c\tau)$. With $(u,  v)=(\psi_y, -\psi_x)=(0, V_c)$ it is easy to see that on the cylinder surface $\psi=D_c-xV_c$ where $D_c$ is the time varying transverse displacement of the cylinder. Again with vorticities produced at no slip boundaries correct and accurate imposition of vorticity boundary conditions is of critical importance. In our previous work \cite{sen_she_17} we have advocated the global integral vorticity condition which is primitive. This procedure was successfully realized as varied explicit boundary vorticity conditions. Further, we were also able to establish a connection between vorticity integral philosophy and a few well known vorticity boundary conditions. Insight acquired lead us to use Briley's vorticity boundary approximation \cite{bri_71} on the surface of the cylinder as an optimum option for this test case.

Referred to Fig. \ref{fig:Cyl_cross_sch} we impose uniform flow conditions in the upstream leading to $\psi=Uy$ and $\omega=0$. Similar expressions are also imposed on the boundaries \textbf{BA}, \textbf{AF}, \textbf{CD}, and \textbf{DE}. Here it is important to mention that although $(u=U, v=0)$ is imposed on boundaries \textbf{BC}, \textbf{DE}, and \textbf{AF} we work with zero Neumann condition $\partial_n(V_n, V_{\tau})=(0, 0)$ for normal ($V_n$) and tangential ($V_{\tau}$) velocity components on \textbf{BA} and \textbf{CD} as discussed in earlier work \cite{kal_sen_12}. Along the downstream boundary \textbf{EF} we impose convective boundary condition for $\psi$ and Neumann conditions for $\omega$, $V_n$, and $V_{\tau}$ as shown in the schematic diagram Fig. \ref{fig:Cyl_cross_sch}. Approximations for gradients of $\omega$ at all boundary arcs are obtained by using one sided approximations following the work of Tian \emph{et al.} \cite{tia_lia_yu_11} except at \textbf{EF} where the normal gradient of $\omega$ is set to zero. 

The drag force $(F_D)$ often normalized by the dynamic free stream pressure $\frac{1}{2}\rho U^2$ referred to as the drag coefficient $(C_D)$ might be represented as 
\begin{equation}\label{osc1}
C_D=\frac{F_D}{\frac{1}{2}\rho U^2D}.
\end{equation}
Following Fornberg \cite{for_80}, $F_D$ evaluated by a line integral around the body,
\begin{equation}\label{osc2}
C_D=\frac{1}{Re}\int^1_{-1}(\omega_{\xi}-\pi\omega)\sin(\pi\eta)d\eta.
\end{equation}
In a similar fashion lift coefficient $(C_L)$ could be written as
\begin{equation}\label{osc3}
C_L=\frac{1}{Re}\int^1_{-1}(\omega_{\xi}-\pi\omega)\cos(\pi\eta)d\eta.
\end{equation}

\begin{figure}[htbp]
\begin{minipage}[b]{.55\linewidth}\hspace{-1cm}
\centering\psfig{file=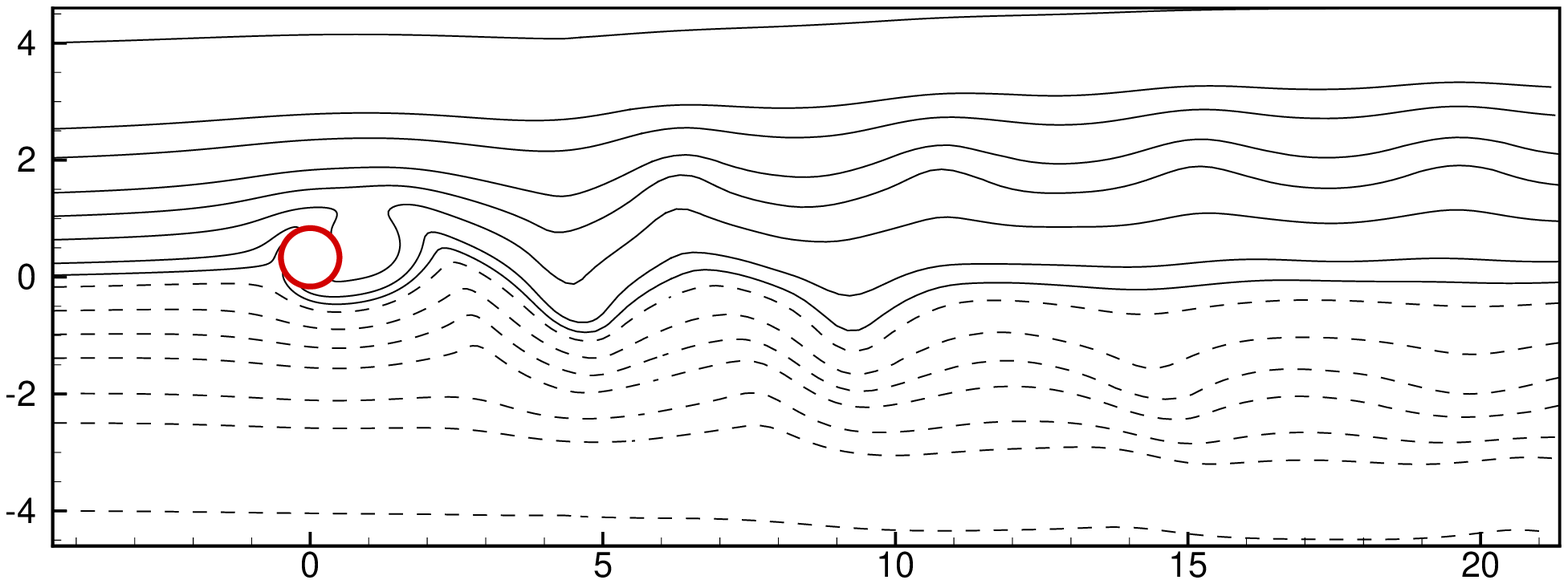,width=0.9\linewidth}
 \\(a)
\end{minipage}
\begin{minipage}[b]{.55\linewidth}\hspace{-1cm}
\centering\psfig{file=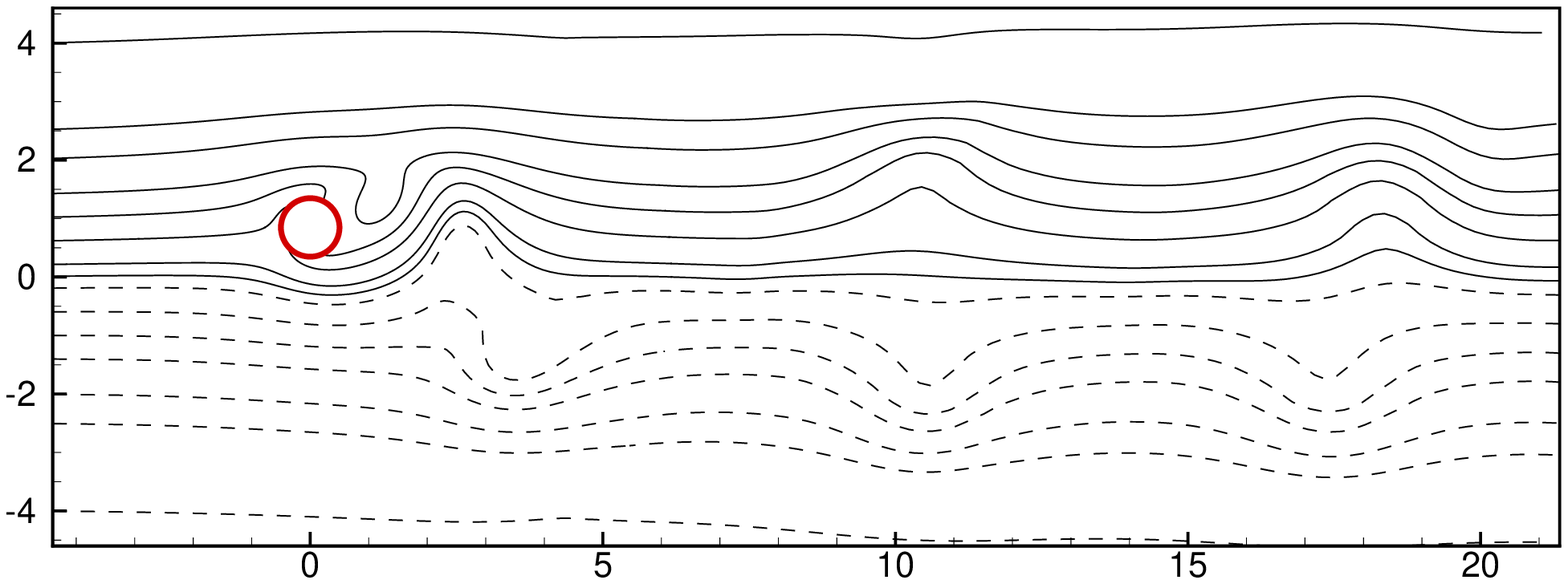,width=0.9\linewidth}
 \\(b)
\end{minipage}
\begin{minipage}[b]{.55\linewidth}\hspace{-1cm}
\centering\psfig{file=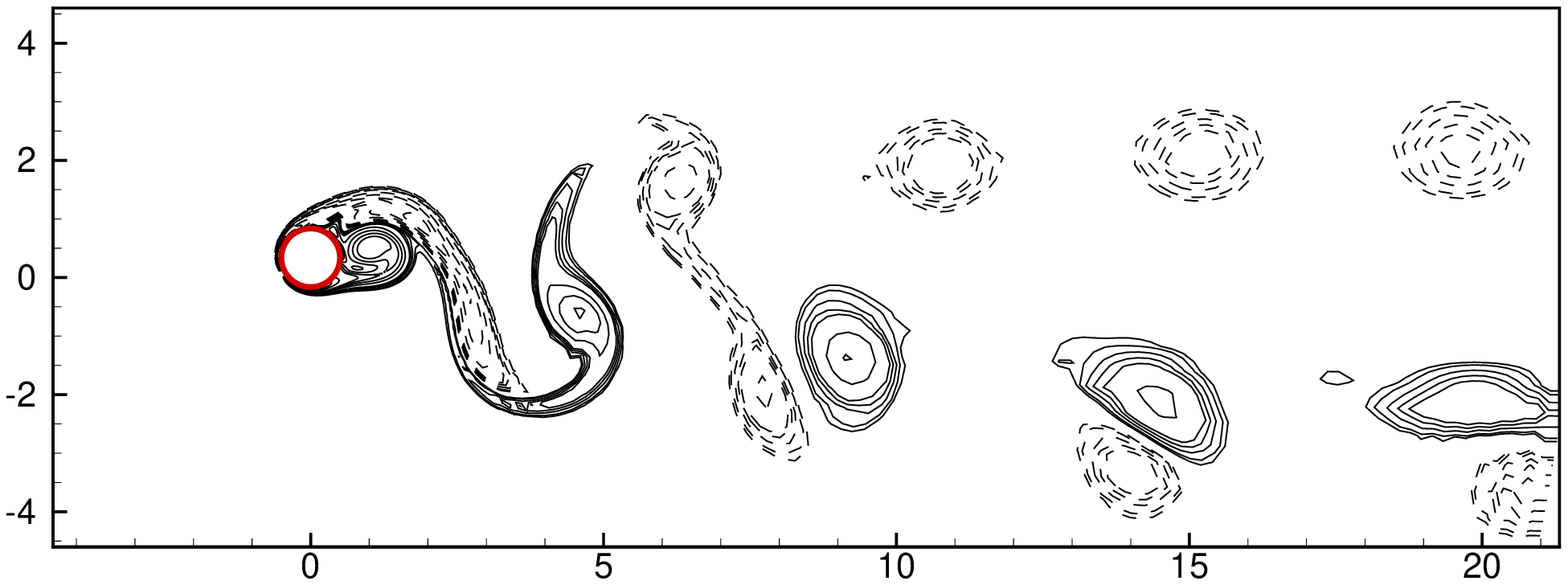,width=0.9\linewidth}
 \\(c)
\end{minipage}
\begin{minipage}[b]{.55\linewidth}\hspace{-1cm}
\centering\psfig{file=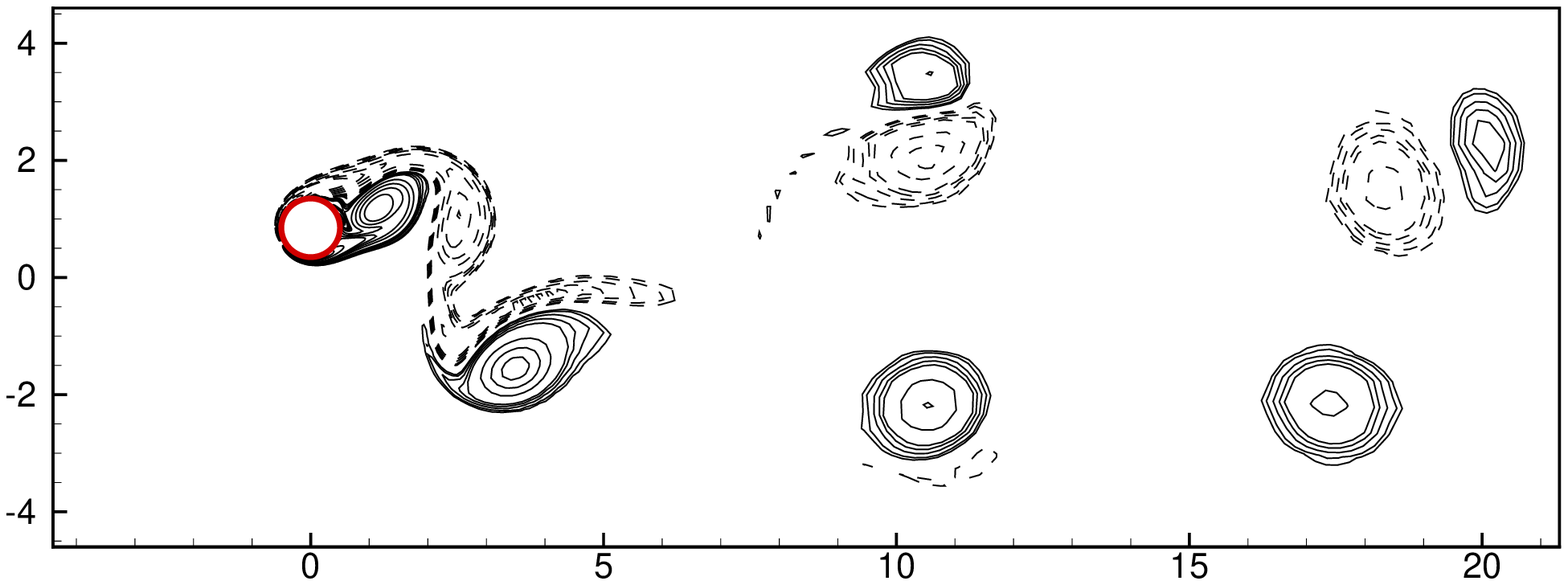,width=0.9\linewidth}
 \\(d)
\end{minipage}
\begin{minipage}[b]{.55\linewidth}\hspace{-1cm}
\centering\psfig{file=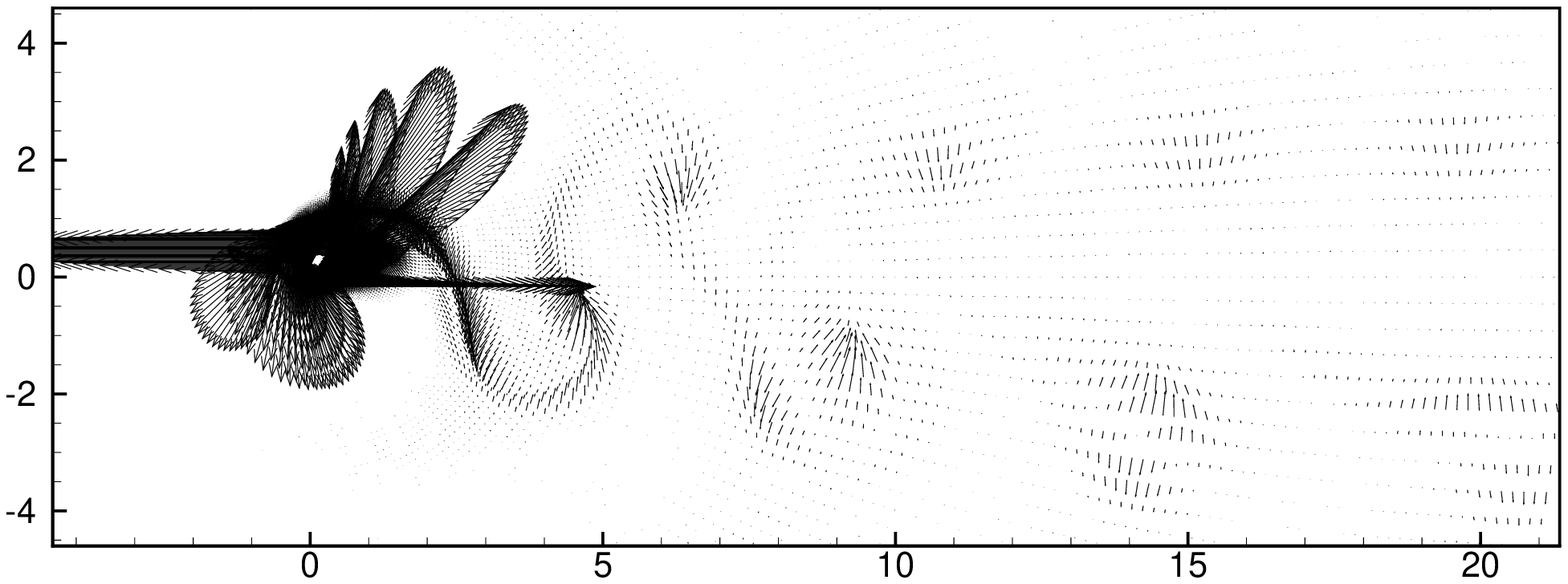,width=0.9\linewidth}
 \\(e)
\end{minipage}
\begin{minipage}[b]{.55\linewidth}\hspace{-1cm}
\centering\psfig{file=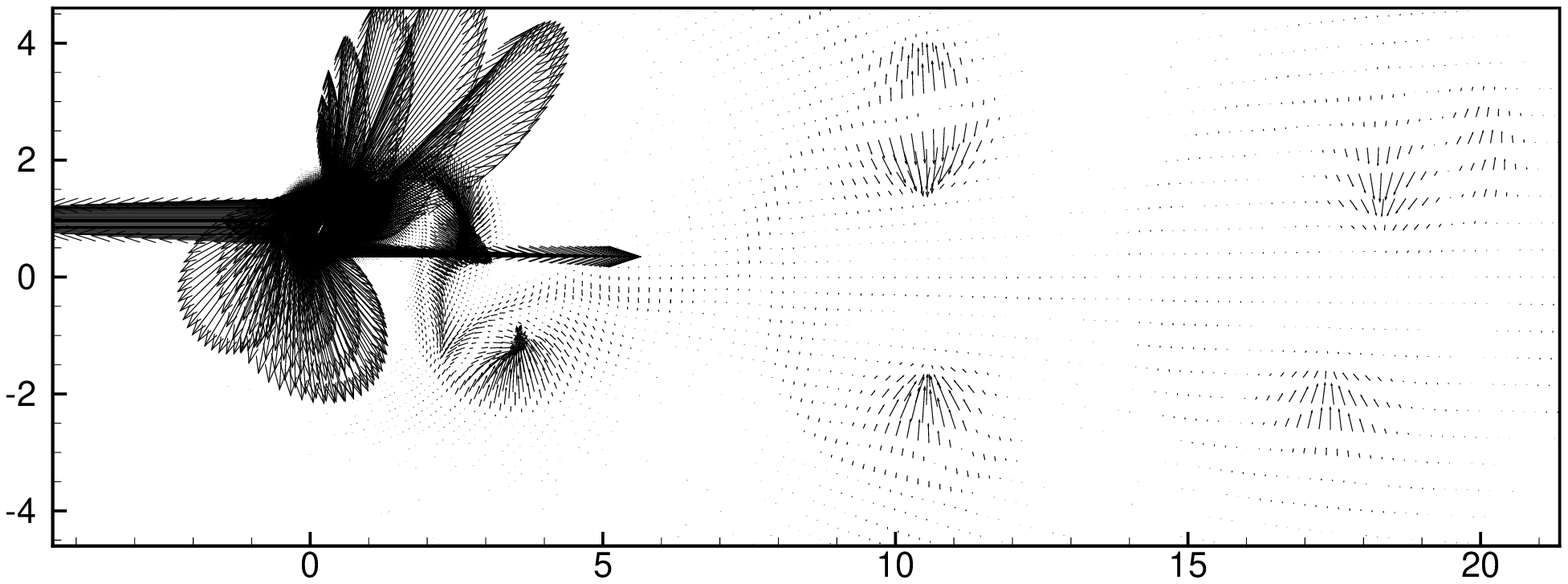,width=0.9\linewidth}
 \\(f)
\end{minipage}
\begin{minipage}[b]{.55\linewidth}\hspace{-1cm}
\centering\psfig{file=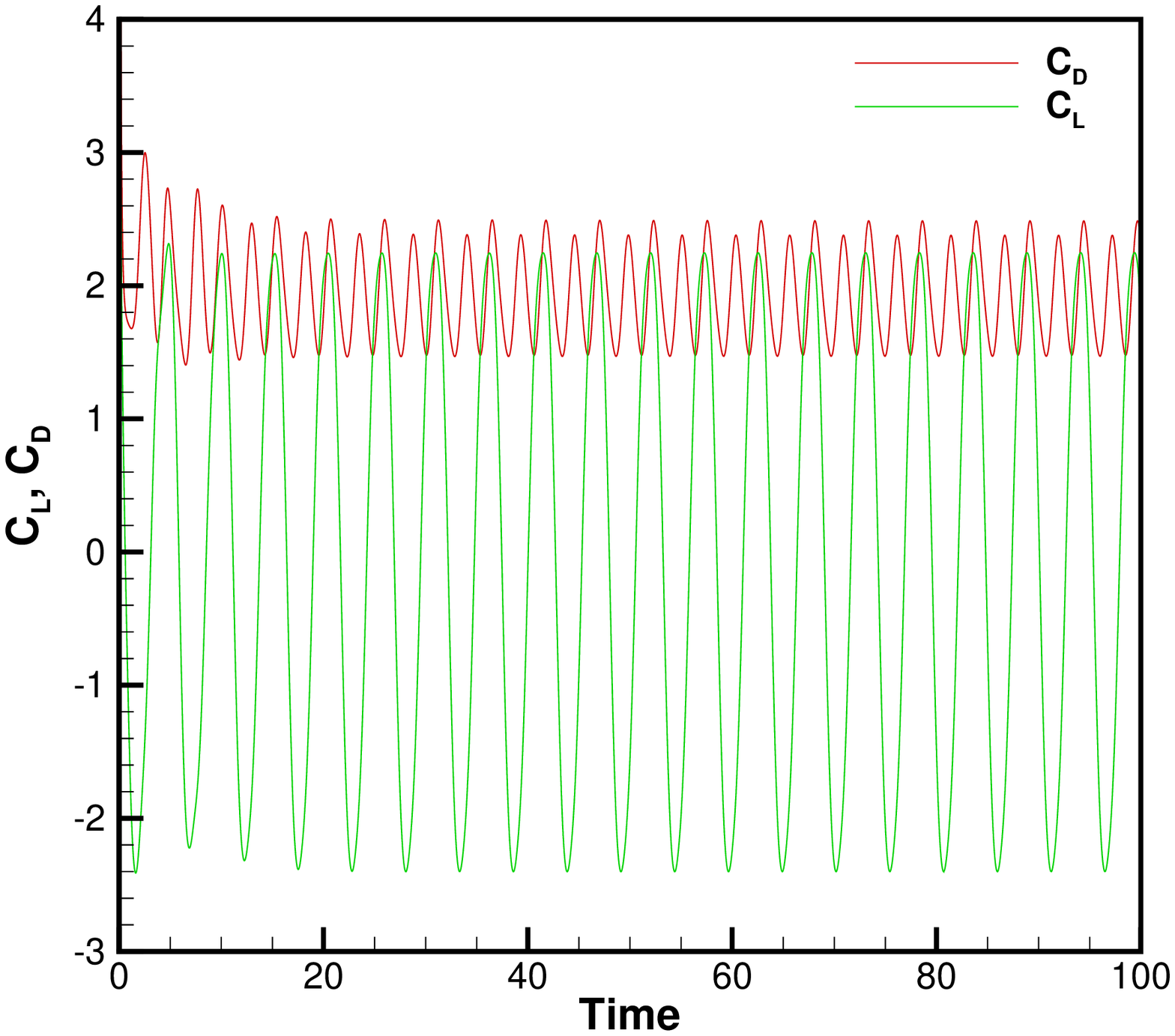,width=0.8\linewidth}
 \\(g)
\end{minipage}
\begin{minipage}[b]{.55\linewidth}\hspace{-1cm}
\centering\psfig{file=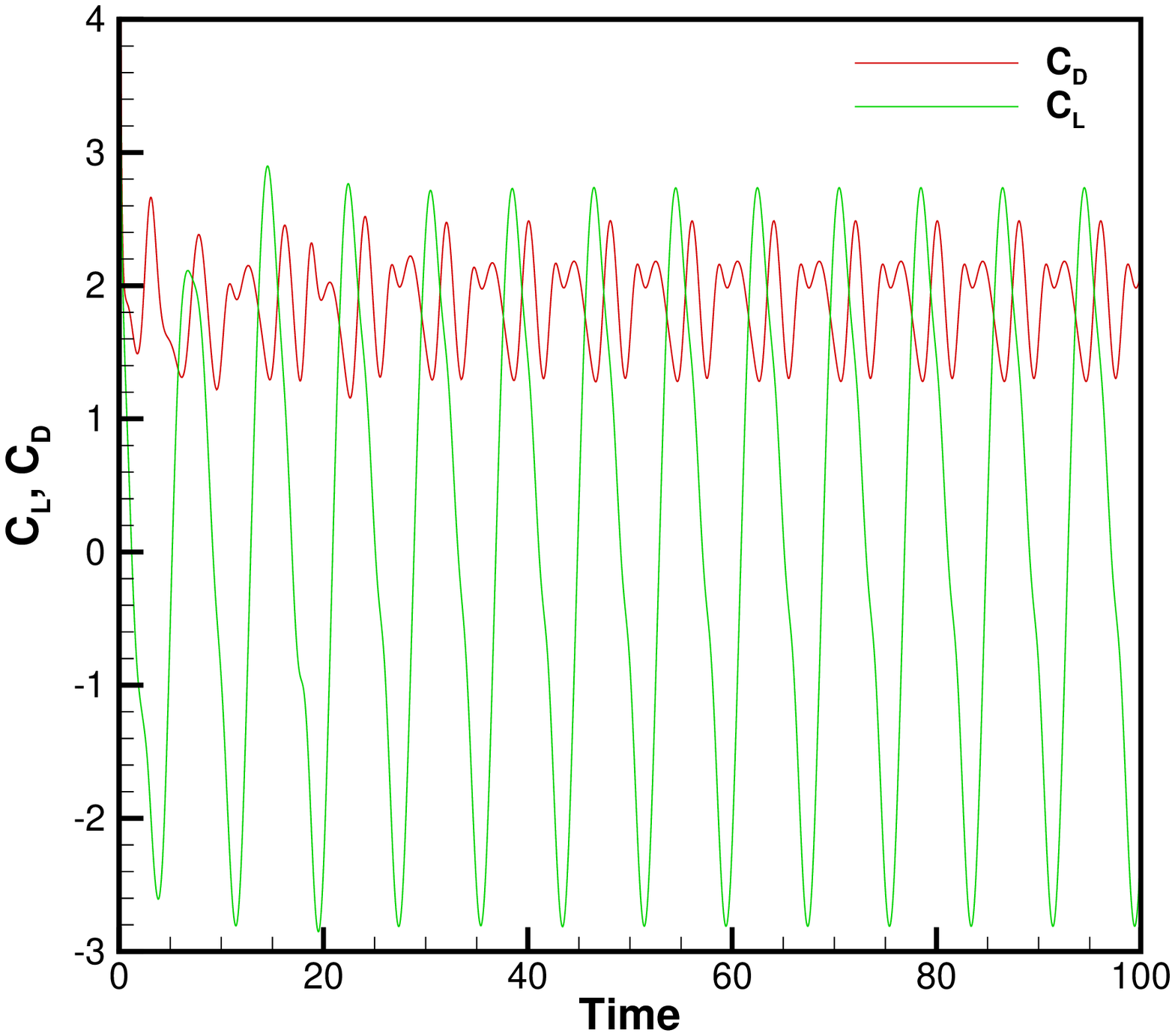,width=0.8\linewidth}
 \\(h)
\end{minipage}
\begin{center}
\caption{{\sl Problem 6: Streamlines (top row), Vorticity contours (second row), Lamb vector plot (third row) and Variation of drag and lift (bottom row) for $Re=200$ at (i) $A_r=0.7$, $f_r=0.95$ (left) and (ii) $A_r=1.2$, $f_r=0.625$ (right). Negative contours appear in dashed lines.} }
    \label{fig:osc_R200}
\end{center}
\end{figure}

Four distinct cases with varied flow parameters are studied here. In the first two cases, we work with $Re=200$ and choose $(A_r, f_r)$ as $(0.7, 0.95)$ and $(1.2, 0.625)$ respectively. Both these cases are partially motivated by the availability of data in the literature \cite{leo_ste_tho_06} and the pressing nature of the specifications. At $A_r=0.7$ we have a comparatively lesser displacement of the cylinder from its mean position at origin but with $f_r=0.95$ we are close to the region of primary synchronization. At such a choice $P+S$ wake mode is reported by Leontini \emph{et al.} \cite{leo_ste_tho_06} in their extensive 2D simulations of an oscillating circular cylinder. From Fig. \ref{fig:osc_R200}(c) it is clear that even with a modest computational setup our scheme is able to accurately capture $P+S$ mode consisting of a pair and single vortex per oscillation cycle. The corresponding streamlines are shown in Fig. \ref{fig:osc_R200}(a). Lamb vector ${\bm \omega}\times\bm u$ is analogous to the Coriolis force observed in a rotating frame of reference. It is of crucial importance in vorticity dynamics and drives the particles to move around the vorticity direction \cite{wu_ma_zho_07}. The Lamb vector plot of the flow field presented in Fig. \ref{fig:osc_R200}(c) identifies asymmetric Karman vortex shedding where the cylinder sheds a paired counter rotating and a single vortex in each cycle. Variation of drag and lift coefficients are presented in Fig \ref{fig:osc_R200}(g). It is seen that after impulsive start the flow settles down to a periodic state at around nondimensional time 30. In the second case with $(A_r, f_r)=(1.2, 0.625)$ the amplitude of oscillation is increased and the cylinder oscillates away from primary synchronization. According to Leontini \emph{et al.} \cite{leo_ste_tho_06}, this choice of parameters leads to a region well over the boundary defining the transition from the $2S$ to $P+S$ wake modes for $Re=200$. Our computations corroborate the same as seen in Fig. \ref{fig:osc_R200}(d). This is further accentuated in the Lamb vector plot in Fig. \ref{fig:osc_R200}(f). Fig. \ref{fig:osc_R200}(h) demonstrates that with increased amplitude of oscillation variation of lift is more pronounced. With cylinder oscillating further away from primary synchronization variation of drag attains three points of maxima instead of the customary two during each time period.

\begin{figure}[htbp]
\begin{minipage}[b]{.55\linewidth}\hspace{-1cm}
\centering\psfig{file=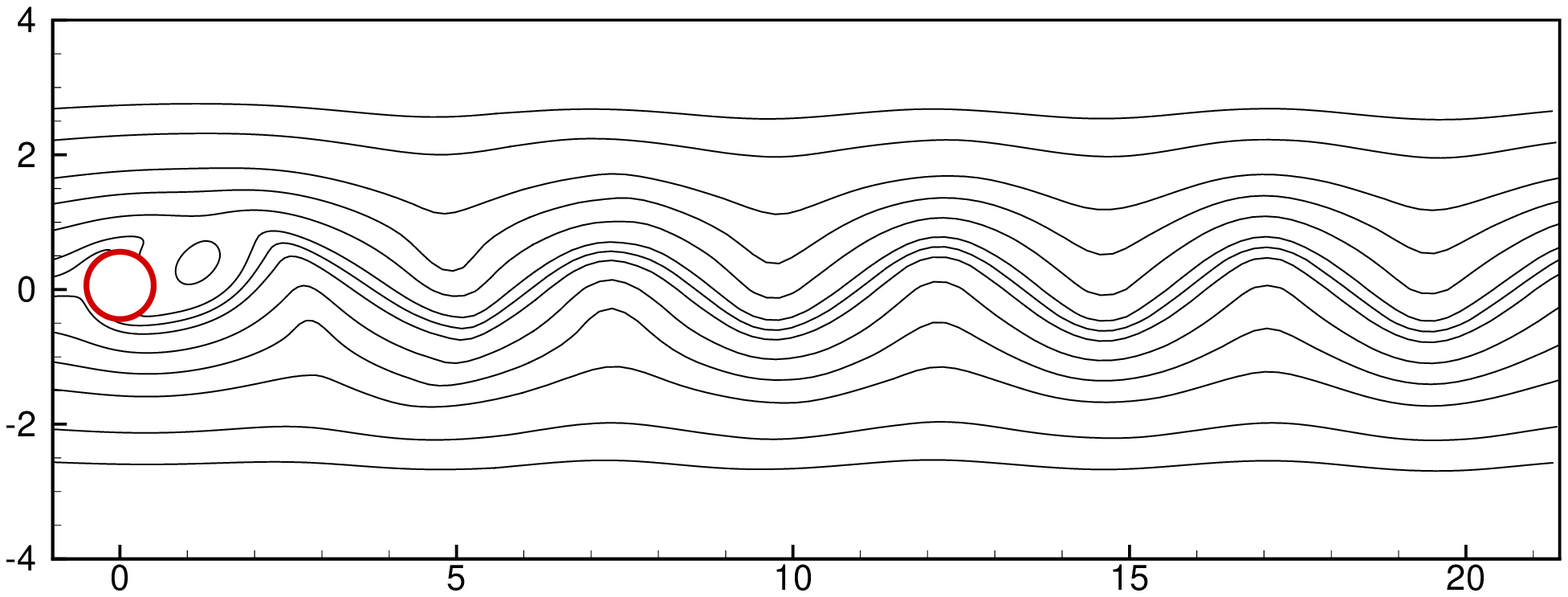,width=0.9\linewidth}
 \\(a)
\end{minipage}
\begin{minipage}[b]{.55\linewidth}\hspace{-1cm}
\centering\psfig{file=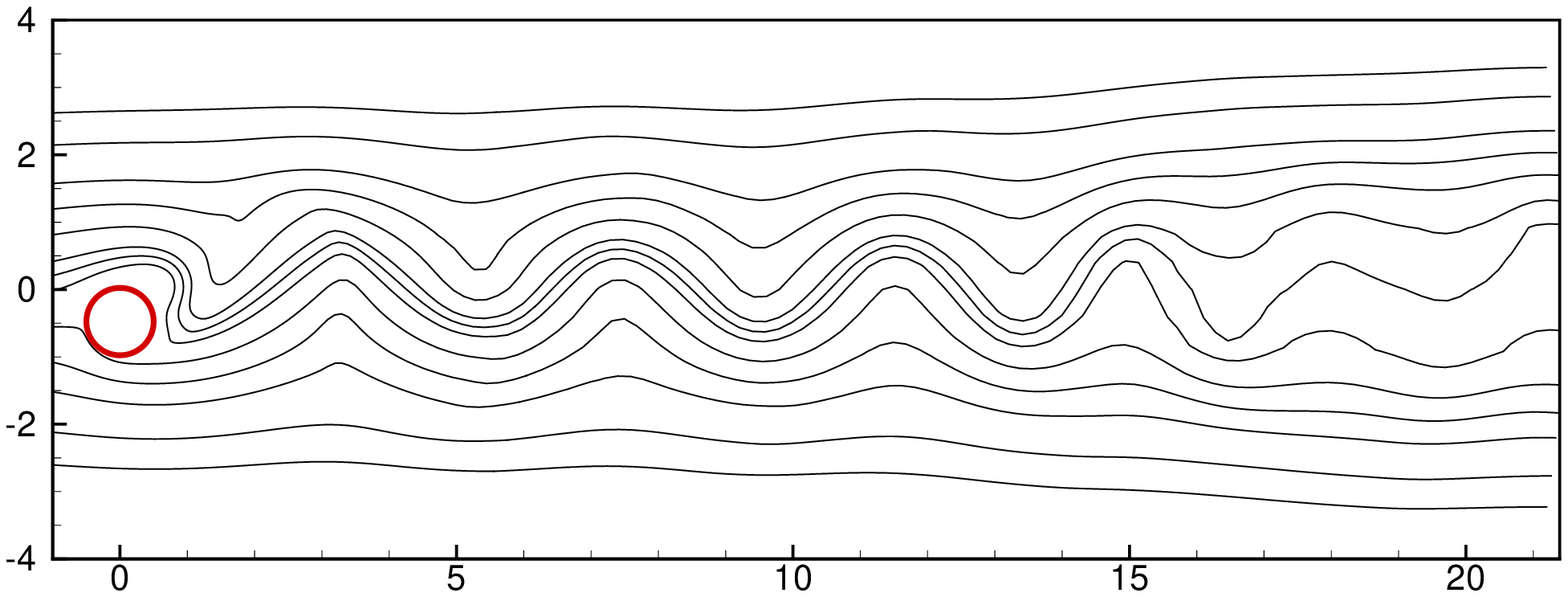,width=0.9\linewidth}
 \\(b)
\end{minipage}
\begin{minipage}[b]{.55\linewidth}\hspace{-1cm}
\centering\psfig{file=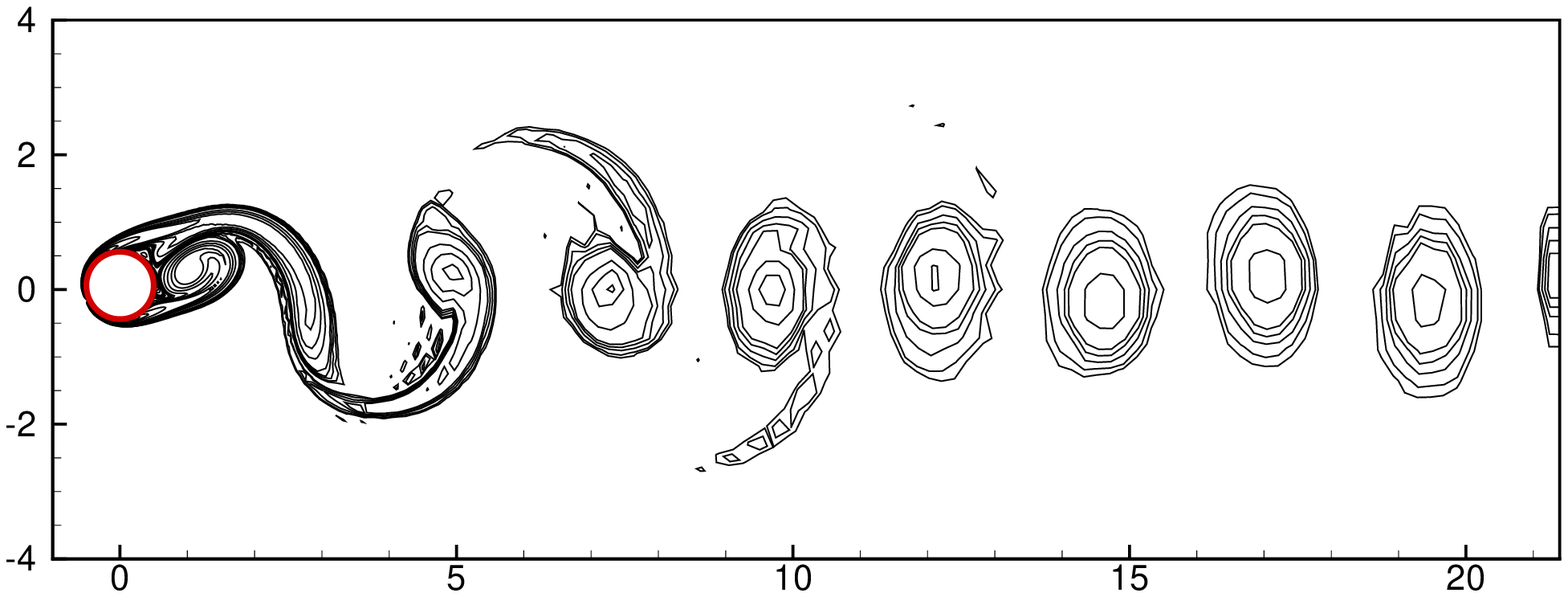,width=0.9\linewidth}
 \\(c)
\end{minipage}
\begin{minipage}[b]{.55\linewidth}\hspace{-1cm}
\centering\psfig{file=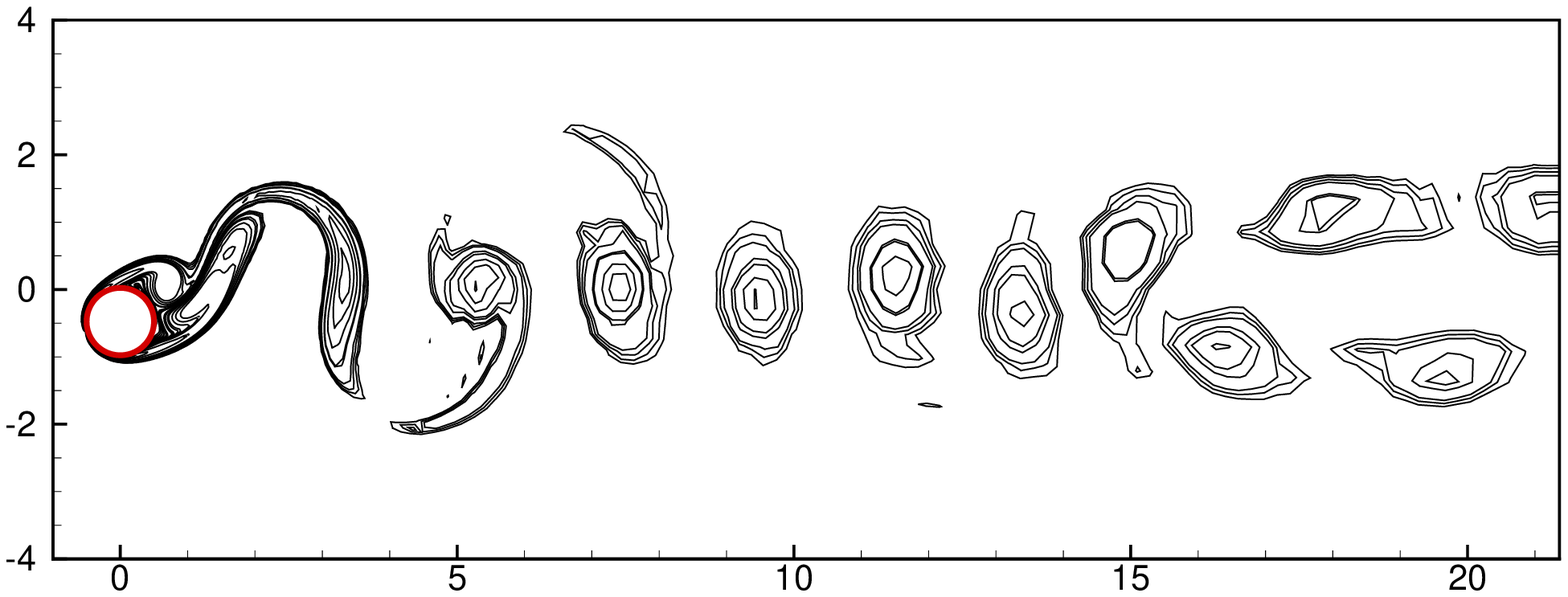,width=0.9\linewidth}
 \\(d)
\end{minipage}
\begin{minipage}[b]{.55\linewidth}\hspace{-1cm}
\centering\psfig{file=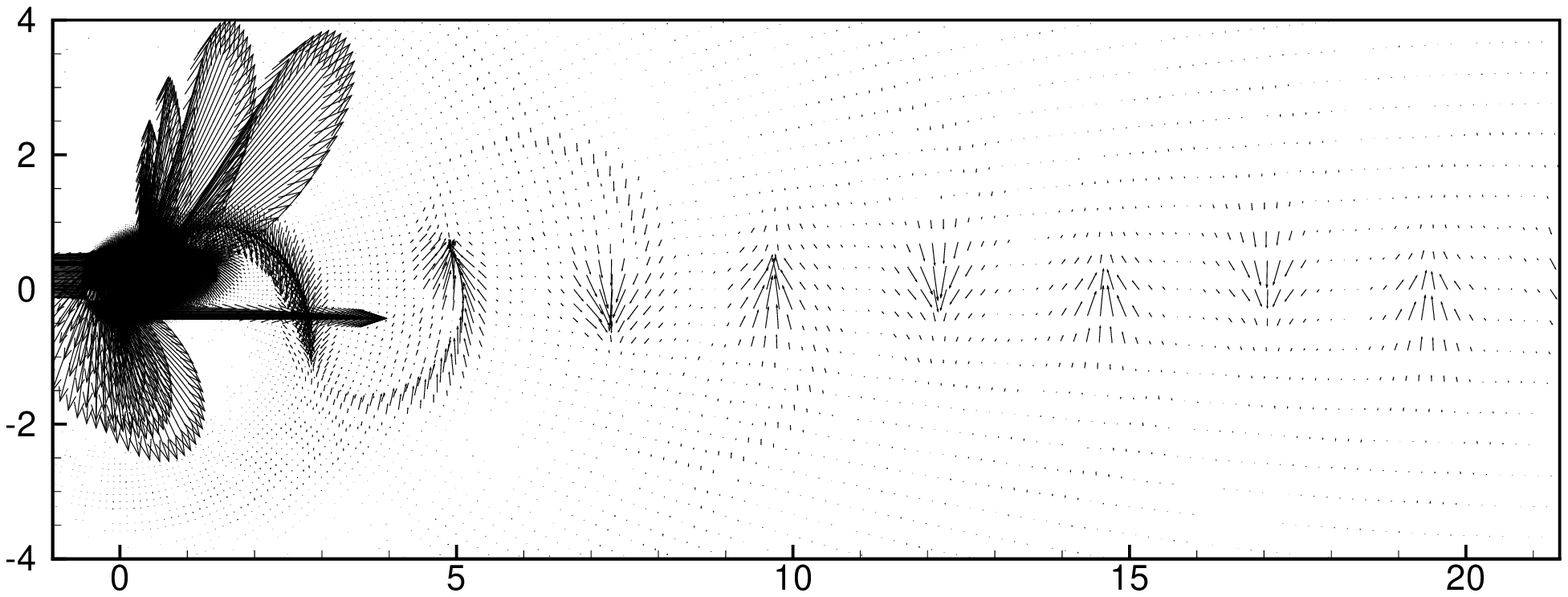,width=0.9\linewidth}
 \\(e)
\end{minipage}
\begin{minipage}[b]{.55\linewidth}\hspace{-1cm}
\centering\psfig{file=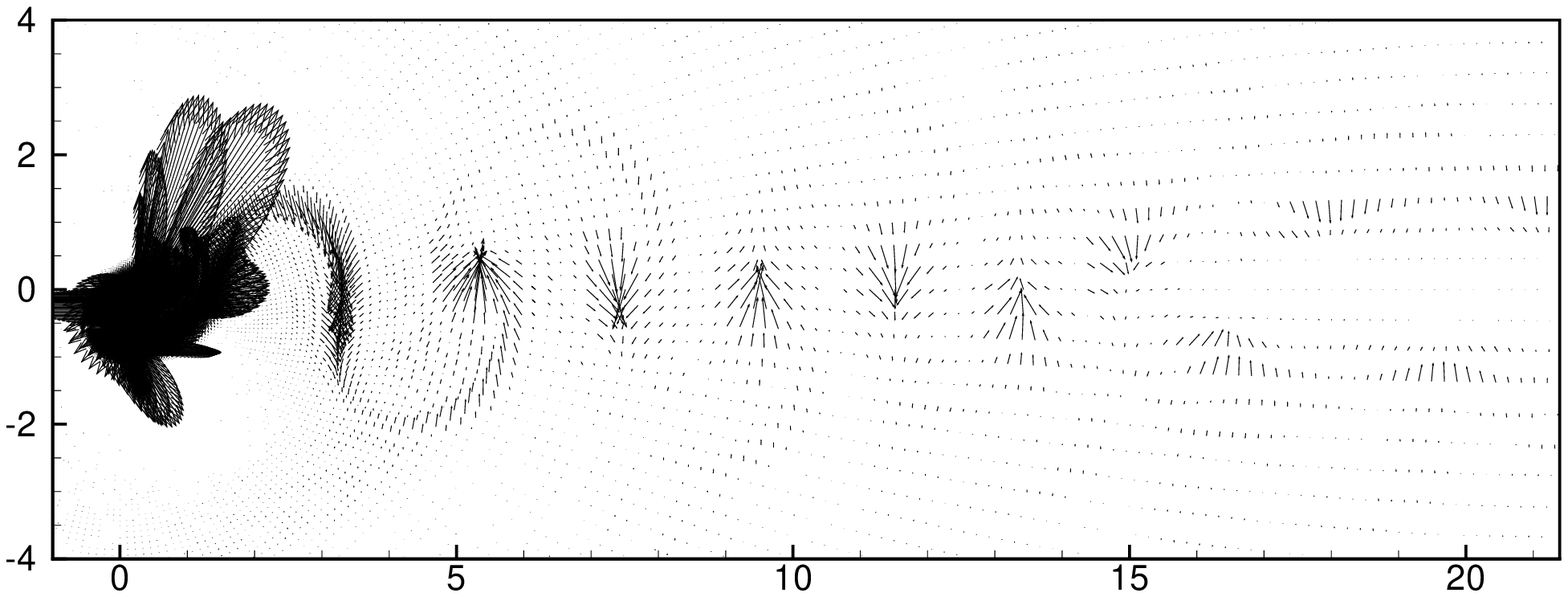,width=0.9\linewidth}
 \\(f)
\end{minipage}
\begin{minipage}[b]{.55\linewidth}\hspace{-1cm}
\centering\psfig{file=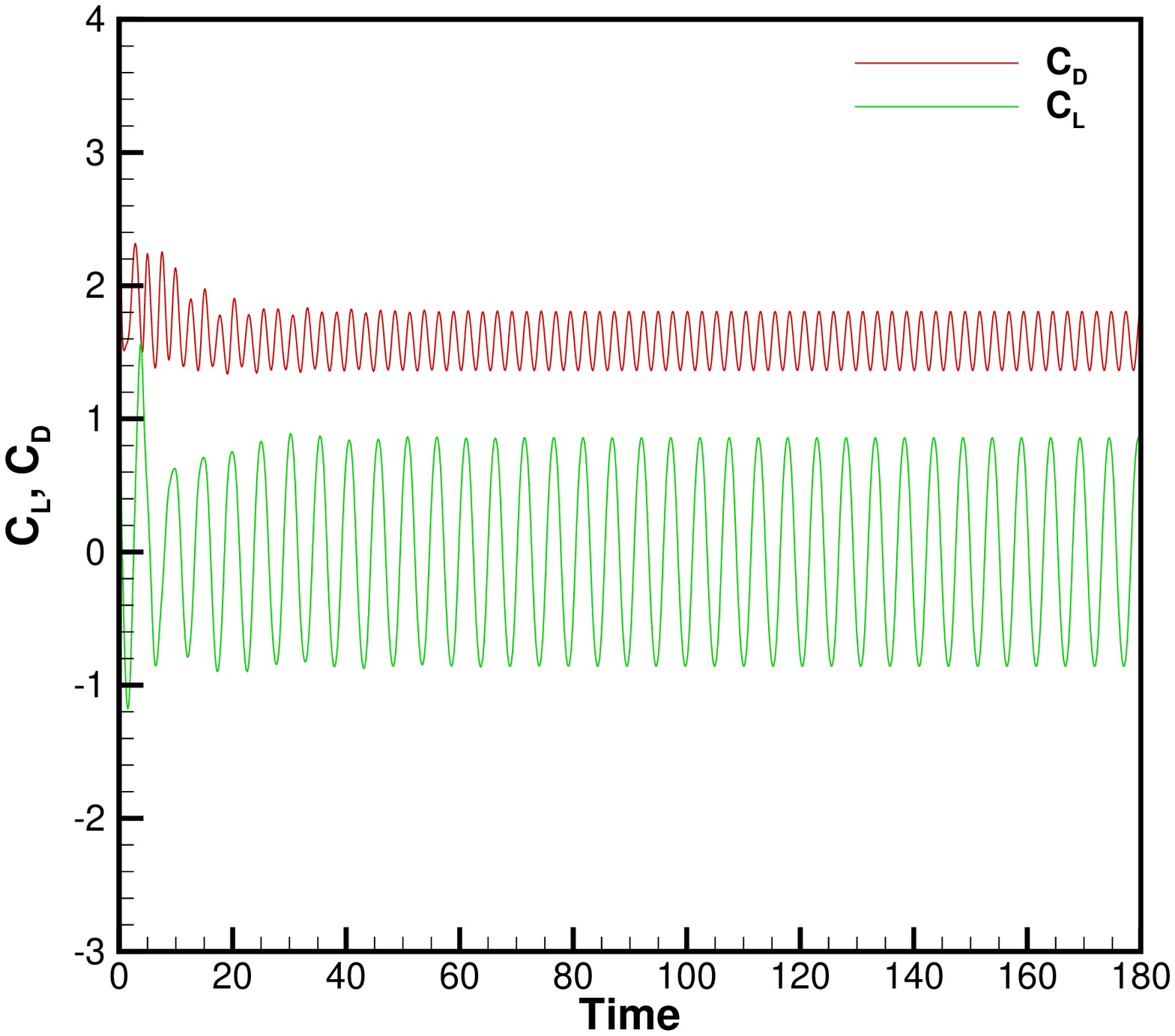,width=0.8\linewidth}
 \\(g)
\end{minipage}
\begin{minipage}[b]{.55\linewidth}\hspace{-1cm}
\centering\psfig{file=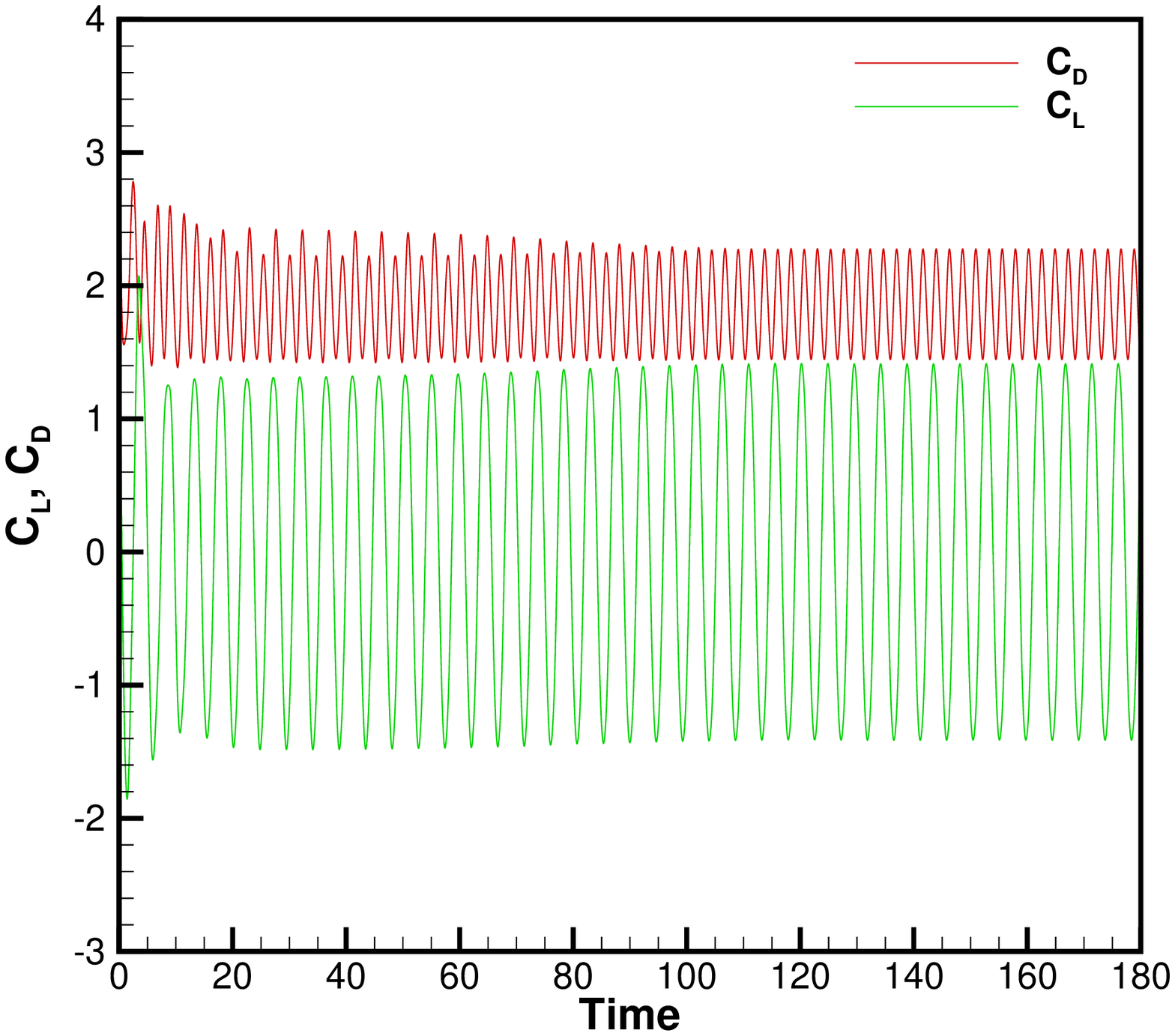,width=0.8\linewidth}
 \\(h)
\end{minipage}
\begin{center}
\caption{{\sl Problem 6: Streamlines (top row), Vorticity contours (second row), Lamb vector plot (third row) and Variation of drag and lift (bottom row) for $Re=392$ at (i) $A_r=0.4$, $f_r=0.971$ (left) and (ii) $A_r=0.5$, $f_r=1.075$ (right).} }
    \label{fig:osc_R392}
\end{center}
\end{figure}

In the second set of studies we work with $Re=392$ and frequency ratios are chosen close to unity (0.971 and 1.075). Amplitude ratios are fixed at 0.4 and 0.5 for the first and second cases respectively. These instances are partly motivated by the experimental work of Williamson and Roshko \cite{wil_ros_88} where authors have documented $2S$ mode of vortex synchronization. Our computation captures periodic flow with $2S$ vortex shedding mode as seen in Figs. \ref{fig:osc_R392}(a), \ref{fig:osc_R392}(c) and \ref{fig:osc_R392}(e). In the fundamental lock-in region, the flow indeed settles to a periodic state at a nondimensional time close to 40. But with a frequency ratio slightly above unity, the flow requires much larger time and gradually attains the periodic state at nondimensional time 100 as seen for the second case in Fig. \ref{fig:osc_R392}(h). Note that lift force also increases with a higher amplitude ratio at 0.5. It is interesting to note that our computation report change in the instantaneous wake structure in the distance downstream. The shed vortices are found to organize themselves into two rows instead of one and are amply clear from the Lamb vector plot in Fig.  \ref{fig:osc_R392}(f). Such a phenomena is reported in the computational study of Leontini \emph{et al.} \cite{leo_ste_tho_06} as well for $Re=200$ at $f_r=1.01$, $A_r=0.5$. Thus it could be inferred our newly developed method is able to capture delicate flow phenomena akin to the oscillating cylinder.
 
\subsection{Pitching and heaving airfoil}
\begin{figure}[htbp]
\begin{minipage}[b]{.5\linewidth}\hspace{-2cm}
\centering\psfig{file=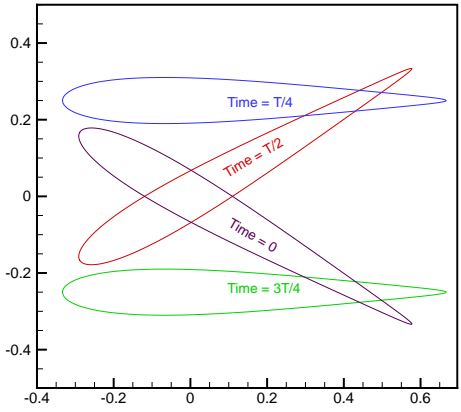,width=0.85\linewidth}(a)
\end{minipage}
\begin{minipage}[b]{.5\linewidth}\hspace{-2cm}
\centering\psfig{file=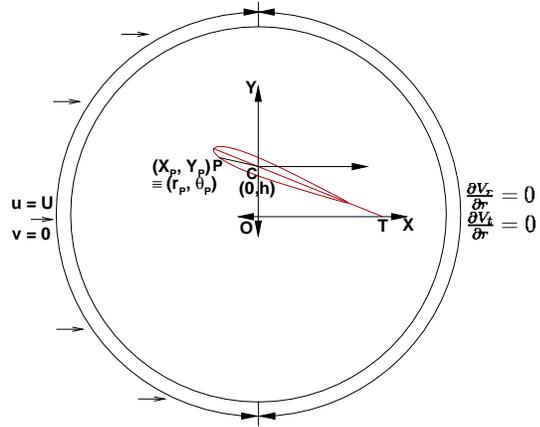,width=0.85\linewidth}(b)
\end{minipage}
\begin{minipage}[b]{.5\linewidth}\hspace{-2cm}
\centering\psfig{file=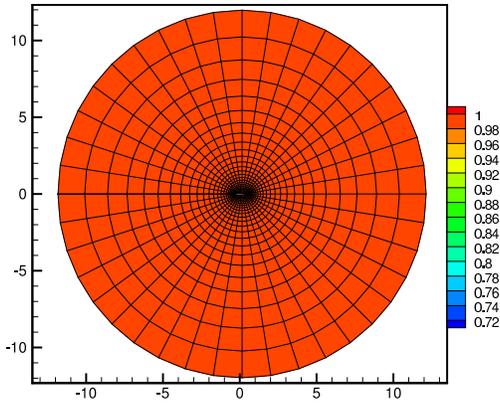,width=0.85\linewidth}(c)
\end{minipage}
\begin{minipage}[b]{.5\linewidth}\hspace{-2cm}
\centering\psfig{file=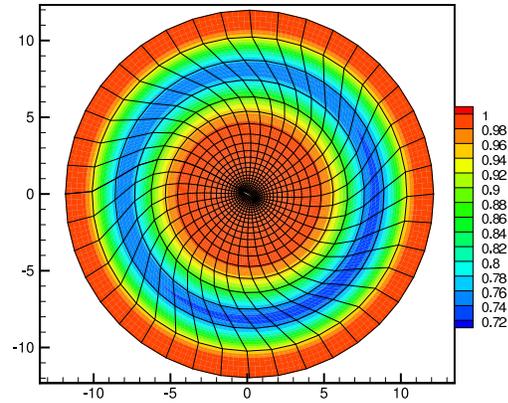,width=0.85\linewidth}(d)
\end{minipage}
\begin{center}
\caption{{\sl Problem 7: (a) Four typical positions of NACA0012 during a complete period. (b) Schematic diagram. (c) Reference grid. (d) Deformed mesh at maximum pitching angle. In (c) and (d) every fifth grid line and skew quality metric are shown around NACA0012.} }
    \label{fig:airfoil_position_grid}
\end{center}
\end{figure}

In this section, we study two airfoils of chord length $c$ pitching and heaving simultaneously subjected to a unit free-stream velocity. The first airfoil is elliptic with a minor to major axes ratio $1:5$ having a pitch centre located at its geometric centre whereas the second airfoil is NACA0012 with a pitch centre at $c/3$ from the leading edge on the axis of symmetry. Following Liang \emph{et al.} \cite{lia_miy_zha_14} for both the airfoils the time dependent rotation motion (pitch) $\alpha(t)$ is prescribed as
\begin{equation}\label{air1}
\alpha(t)=\alpha_{\max}\sin(2\pi f_pt+\pi/2)
\end{equation}
and the translation motion (heave) is given by
\begin{equation}\label{air2}
h(t)=h_{\max}\sin(2\pi f_ht).
\end{equation}
Here we take $\alpha_{\max}=\pm\pi/6$, $h_{\max}=0.25$ and $f_p=f_h=0.4$. Four typical positions of NACA0012 during one complete period are shown Fig. \ref{fig:airfoil_position_grid}(a). A similar motion could be traced for elliptic airfoil as well. 

A schematic diagram depicting the position of NACA0012 airfoil at a given time along with prescribed velocity boundary conditions is shown in Fig. \ref{fig:airfoil_position_grid}(b). For correct flow simulation, it is imperative to prescribe accurate conditions on the surface of the airfoil. In this context, we shall like to provide a sketch of steps required to arrive at streamfunction values on the surface of the airfoil. The motion of any point $P$ on the surface of the airfoil can be described as a result of the motion of the pitch centre $C$ and the motion relative to $C$. From Eq. (\ref{air2}) horizontal and vertical velocity components at $C$ are $0$ and $\dot{h}$ respectively whereas radial and cross-radial components of velocity at $P$ relative to $C$ are $\dot{r}_p=0$ and $r_p\dot{\theta}_p=r_p\dot{\alpha}$ respectively. Thus horizontal and vertical velocity components relative to $C$ of the point $P$ are $r_p\dot{\alpha}\sin(\theta_p)=\dot{\alpha}(Y_p-h)$ and $r_p\dot{\alpha}\cos(\theta_p)=\dot{\alpha}X_p$. Thus resultant velocity at $P$ is $u=\dot{\alpha}(Y_p-h)$, $v=\dot{h}+\dot{\alpha}X_p$. In the transformed $(\xi, \eta)$ coordinate system using $$\psi_{\xi}=(y_{\xi}u-x_{\xi}v),\;\;\;\psi_{\eta}=(y_{\eta}u-x_{\eta}v)$$ it is easy to see that 
\begin{equation}\label{air3}
\psi=\frac{\dot{\alpha}}{2}(Y^2_p-X^2_p)-\dot{\alpha}hY_p-\dot{h}X_p.
\end{equation}

Simulations precede by generating an O-type mesh of radius $\approx 12c$ encircling the airfoil. Undeformed mesh is generated following the conformal transformation 
\begin{equation}\label{air4}
z=\lambda+c_1\frac{a}{\lambda}+c_2\frac{a^2}{\lambda^2}+z_0,\;\;\;\lambda=ae^{i\omega},\;\;\;\omega=\xi+i\eta,\;\;\;a=\frac{1}{2}-c_1
\end{equation}
with $0\le\xi\le1.2\pi$, $-\pi\le\eta\le\pi$ as advocated by Sen \emph{et al.} \cite{sen_kal_gup_13} and presented in Fig. \ref{fig:airfoil_position_grid}(c). For elliptic airfoil we take $c_1=0.2$ and $c_2=0.0$ whereas for NACA0012 the corresponding values are 0.2238 and 0.0165. $z_0$ is appropriately chosen to align the pitch centre with the origin. In this initial mesh, the axis of the airfoil lies along the $x$-axis. Subsequently, mesh adaptation is carried out using IDW remeshing technique. In this context, up to 10\% grid in the immediate vicinity of the airfoil is stipulated to move analogous to the airfoil. Thus conformal nature of the deformed grid is retained adjacent to the airfoil. As this problem involve the pitching of airfoils it is imperative to determine rotation quaternions and interpolate as mentioned earlier. A hybrid version of this strategy is recently used with great success by Apostolatos \emph{et al.} \cite{apo_nay_ble_19} and De Nayer \emph{et al.} \cite{nay_bre_woo_20}. An efficient realization in the current study might be inferred from Fig. \ref{fig:airfoil_position_grid}(d) where a deformed grid corresponding to maximum pitching angle along with associated skew quality metric defined in Knupp \cite{knu_03} is presented. It is heartening to note that restricted skewness is lost in a limited region of the computational domain. For both cases, computations are carried out using $121\times201$ grid with $\delta\tau=0.0025$ for chord length based Reynolds number of 200. Choice of $Re$ is motivated by our desire to compute within the laminar flow regime.

Having obtained boundary conditions for streamfunction on the surface of the airfoil we shall like to document the use of potential flow conditions in the upstream as shown in Fig. \ref{fig:airfoil_position_grid}(b). We impose $\psi=Uy$, $\eta\in\{-\pi,-\pi/4\}\cup\{\pi/4,\pi\}$ and convective boundary condition elsewhere. Vorticity on the surface of the airfoil is computed using the recently developed philosophy of Sen and Sheu \cite{sen_she_17}. Specifically, we employ a cell centre approach to write
\begin{equation}\label{air5}
\omega_{0,j}=-\omega_{1,j}+\frac{2}{J_{0,j}+J_{1,j}}\left(\frac{2}{h}(\psi_{\xi_{0,j}}-\psi_{\xi_{1,j}})+\frac{1}{2k}(\psi_{\eta_{0,j-1}}+\psi_{\eta_{1,j-1}}-\psi_{\eta_{0,j+1}}-\psi_{\eta_{1,j+1}})\right)\;\;\;\forall\;j.
\end{equation}
All remaining vorticity conditions are arrived at using identical philosophy as discussed in the previous section. 

\begin{figure}[h!]
\begin{minipage}[b]{.5\linewidth}\hspace{-1cm}
\centering\psfig{file=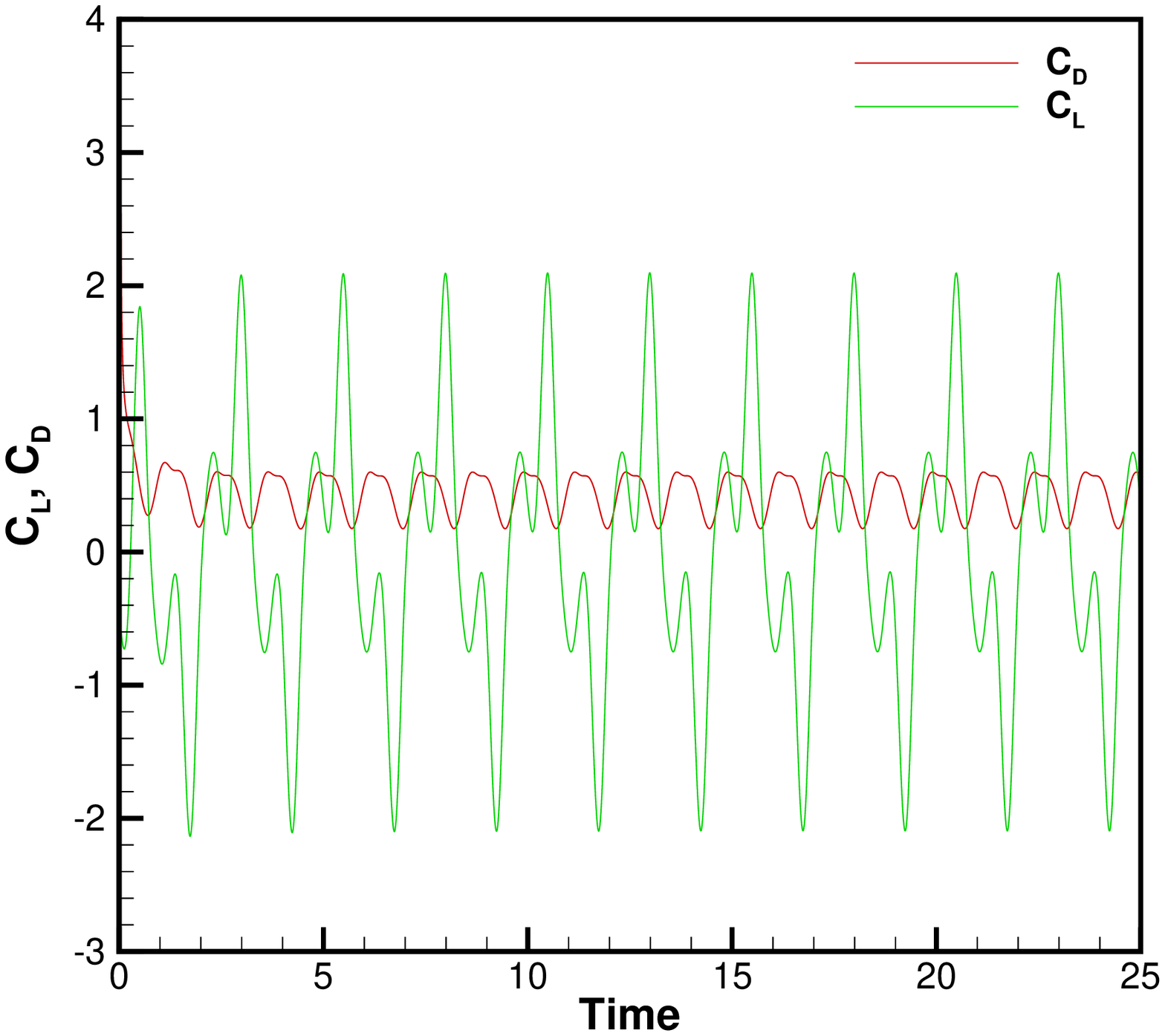,width=0.9\linewidth}
\\(b)
\end{minipage}
\begin{minipage}[b]{.5\linewidth}\hspace{-1cm}
\centering\psfig{file=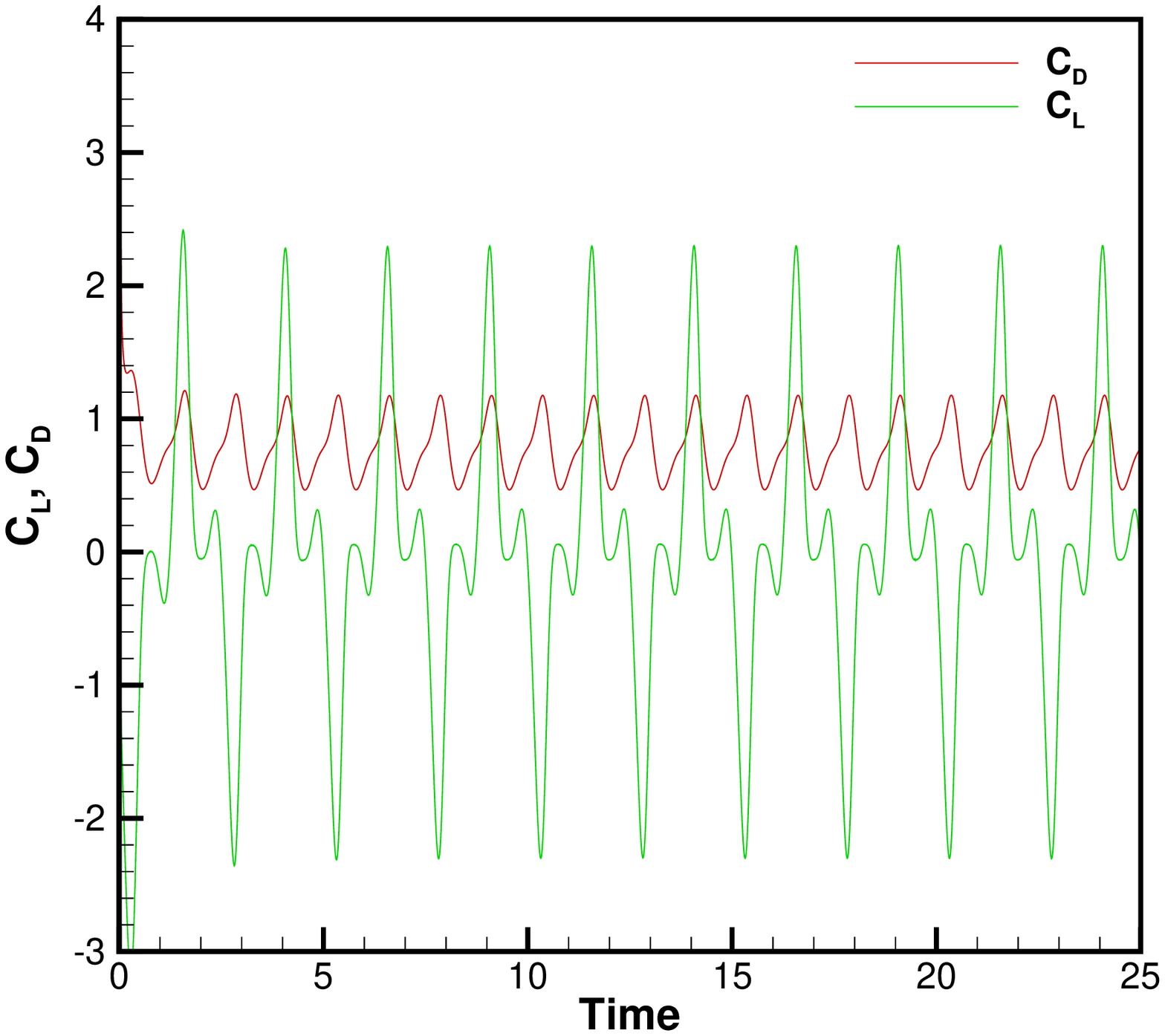,width=0.9\linewidth}
\\(b)
\end{minipage}
\begin{center}
\caption{{\sl Problem 7: Time evolution of drag and lift coefficients: (a) Elliptic airfoil, (b) NACA0012.} }
    \label{fig:airfoil_dl}
\end{center}
\end{figure}

Time variations of $C_D$ and $C_L$ are presented in Fig. \ref{fig:airfoil_dl}. For the flow past blunt bodies, these coefficients at the surface of the body are two important parameters. The time evolution of these two characteristic parameters illustrates the variation of the flow field. Hence we present below details of their calculation. 

From Eq. (\ref{air4}) it is clear that in the undeformed position parametric equation of airfoil is
\begin{subequations}\label{air6}
\begin{empheq}[left=\empheqlbrace]{align}
&x=-c_2+(a+c_1)\cos(\eta)+c_2\cos(2\eta)+x_0, \label{air6_1}\\
&y=(a-c_1)\sin(\eta)-c_2\sin(2\eta)+y_0,\label{air6_2}
\end{empheq}
\end{subequations}
with $-\pi\le\eta\le\pi$. As the airfoil execute pitching and heaving motion its transformed position at any subsequent time $\tau$ concerning the absolute frame of reference might be expressed as
\begin{equation}\label{air7}
\bm X=\mathcal{\bm A}\bm x+\bm H
\end{equation}
where $\bm X=(X, Y)^T$, $\bm x =(x,y)^T$, $\bm H=(0, h(\tau))^T$ and $\mathcal{\bm A}=\left[ {\begin{array}{cc}
   \cos\alpha(\tau) & -\sin\alpha(\tau) \\
   \sin\alpha(\tau) & \cos\alpha(\tau) \\
  \end{array} } \right]$. The unit normal vector to the airfoil at time $\tau$ is thus given by
\begin{equation}\label{air8}
\bm n(\eta, \tau)=(l, m)^T=\left(\frac{dY}{d\eta}, -\frac{dX}{d\eta}\right)^T\scalebox{2.5}{/}\sqrt{\left(\frac{dX}{d\eta}\right)^2+\left(\frac{dY}{d\eta}\right)^2}.
\end{equation}  

Following Eldredge \cite{eld_07} the aerodynamics force on the airfoil is
\begin{equation}\label{air9}
\bm F=\mu\oint_{\mathcal{C}}\left((\bm X-\bm X_c)\times\frac{\partial \bm \omega}{\partial \bm n}-\bm n\times\bm \omega\right)ds+\rho \bar{A} \dot{\bm H}
\end{equation}
where $\mathcal{C}$ is any contour surrounding the airfoil and $s$ is measured around the contour $\mathcal{C}$ in the counterclockwise sense. Here $\mu$ is dynamic viscosity, $\rho$ density, $\bar{A}$ surface area of an airfoil, $\bm X_c$ co-ordinate of the centre of rotation and $\dot{\bm H}$ is the translation velocity of the airfoil. As noted by Eldredge \cite{eld_07} the term outside integration is the inertial term which along with the first term of the integrand account for the pressure effect. On the other hand, the second integrand captures the viscous shear stress effect. We thus arrive at dimensionless components 
\begin{eqnarray}\label{air10}
\frac{\bm F}{\rho U^2 c}&=&\frac{1}{Re}\oint_{\mathcal{C}}\left((\bm X-\bm X_c)\times\frac{\partial \bm \omega}{\partial \bm n}-\bm n\times\bm \omega\right)ds+\frac{\bar{A}}{U^2 c}\dot{\bm H}\nonumber\\
&=&\left(\frac{1}{Re}\oint_{\mathcal{C}}\left((Y-Y_c)\frac{\partial \omega}{\partial \bm n}-m\omega\right)ds,\frac{1}{Re}\oint_{\mathcal{C}}\left(-(X-X_c)\frac{\partial \omega}{\partial \bm n}+l\omega\right)ds+\frac{\bar{A}}{U^2 c}\dot{h}\right)\nonumber\\
&=&(F_{X,a}, F_{Y,a}).
\end{eqnarray}
$F_{X,a}$ and $F_{Y,a}$ are dimensionless aerodynamic force components along co-ordinates axes in an absolute frame of reference. The drag and lift coefficients $C_D$ and $C_L$ respectively are calculated using the forces determined in the moving frame of reference \cite{den_you_03, med_sto_car_11} and are thus given by
\begin{eqnarray}\label{air11}
(C_D, C_L)^T&=&\mathcal{\bm A}^T(F_{X,a}, F_{Y,a})^T\nonumber\\
&=&(F_{X,m}, F_{Y,m}).
\end{eqnarray}
Here $F_{X,m}$ and $F_{Y,m}$ are horizontal and vertical force components respectively in the moving frame attached to the airfoil.

\begin{figure}[htbp]
\begin{minipage}[b]{.55\linewidth}\hspace{-1cm}
\centering\psfig{file=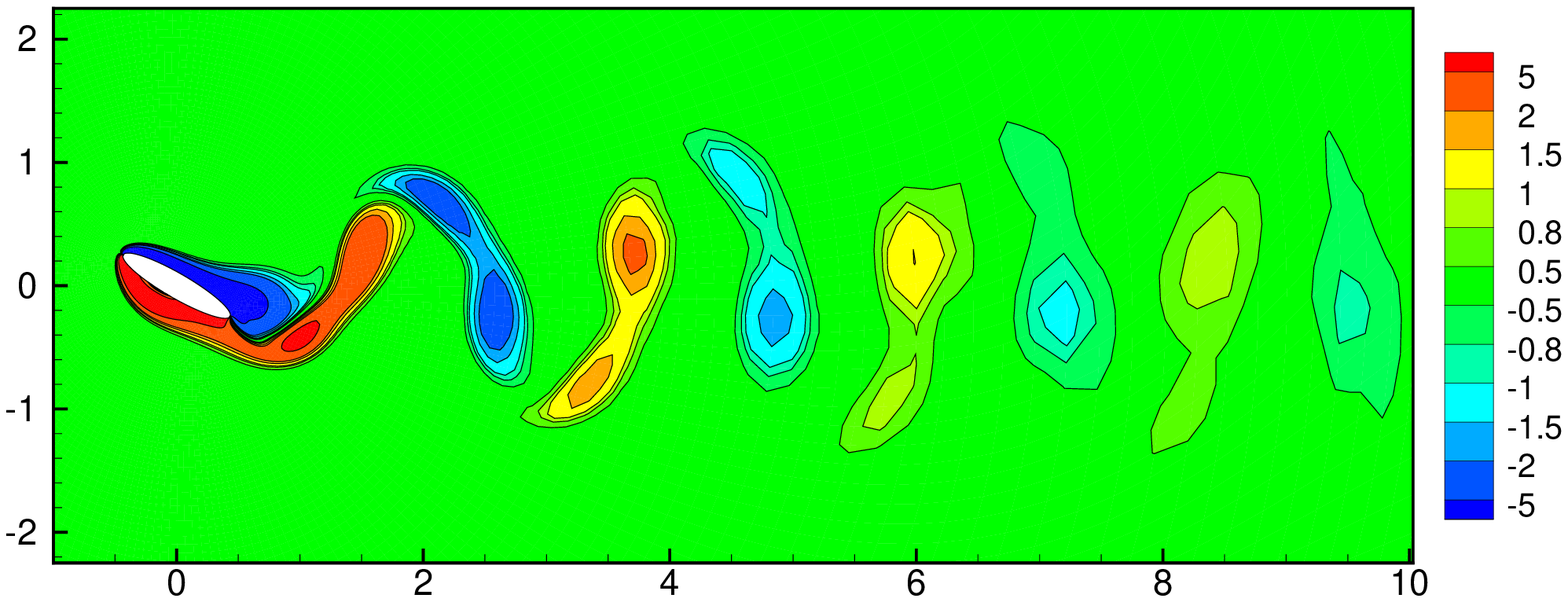,width=0.9\linewidth}
 \\(a) $t_0$
\end{minipage}
\begin{minipage}[b]{.55\linewidth}\hspace{-1cm}
\centering\psfig{file=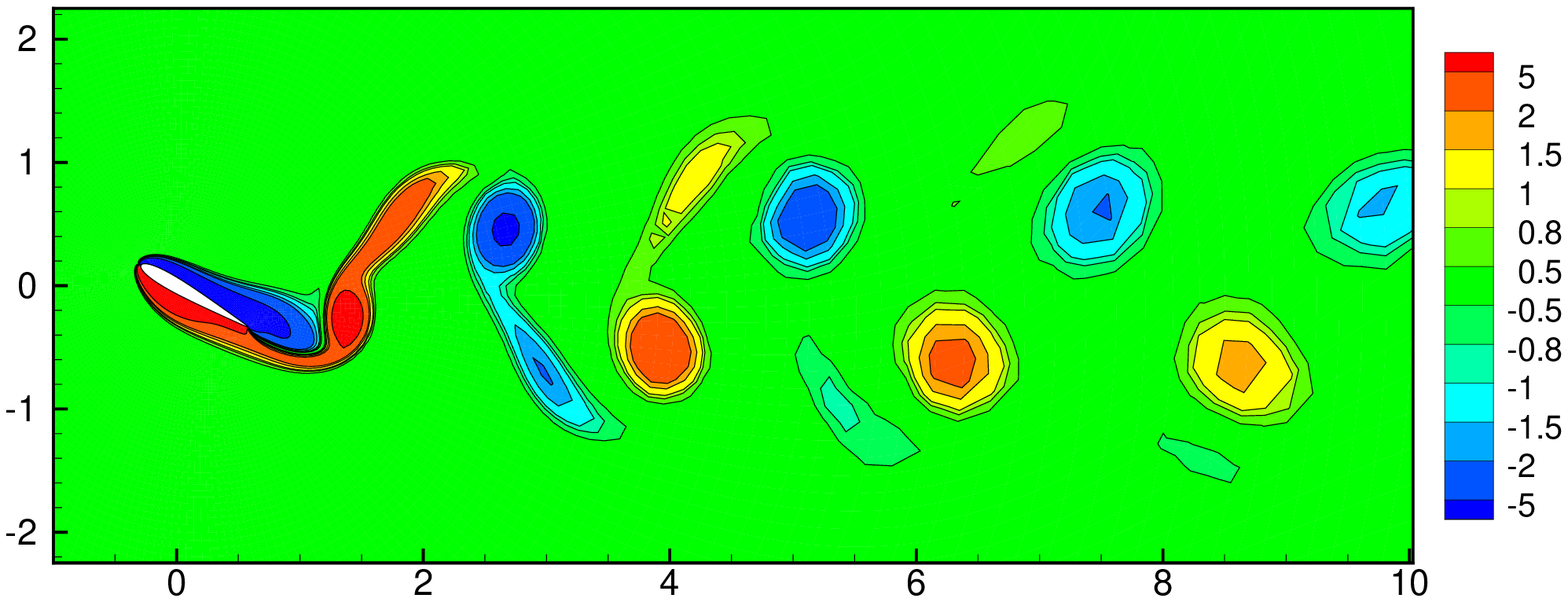,width=0.9\linewidth}
 \\(b) $t_0$
\end{minipage}
\begin{minipage}[b]{.55\linewidth}\hspace{-1cm}
\centering\psfig{file=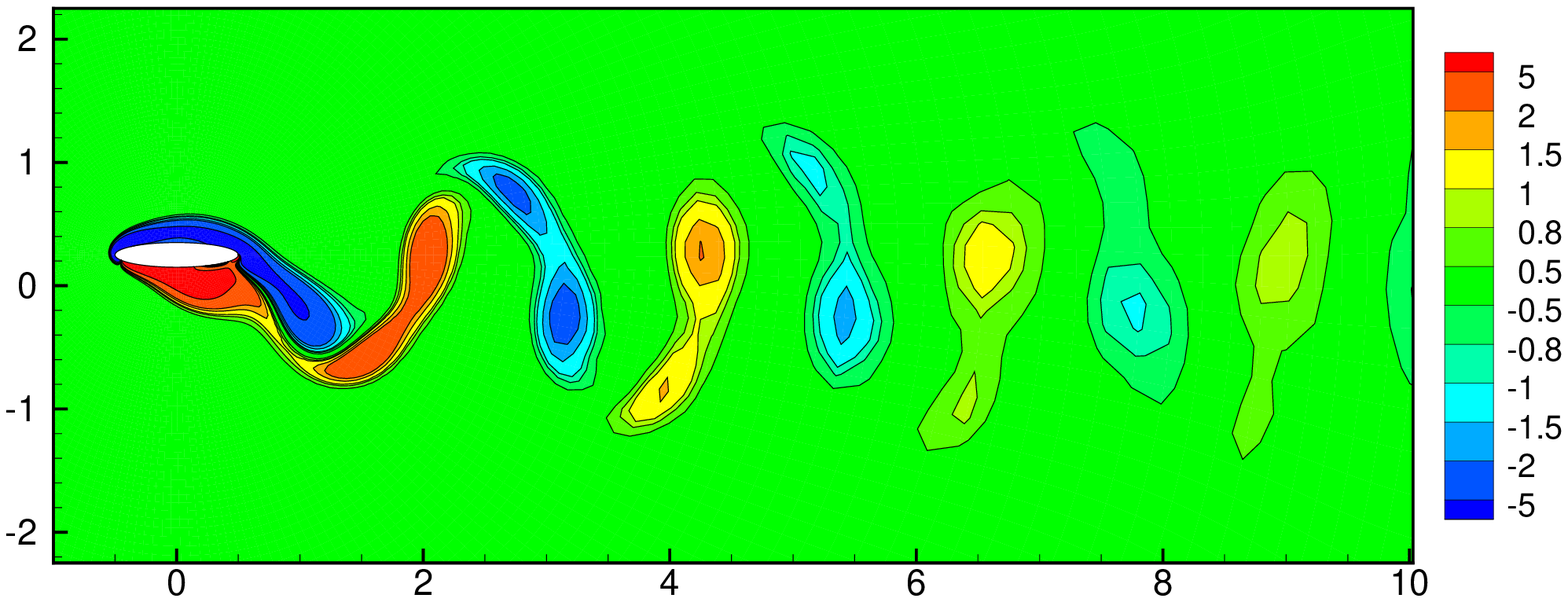,width=0.9\linewidth}
 \\(c) $t_0+T/4$
\end{minipage}
\begin{minipage}[b]{.55\linewidth}\hspace{-1cm}
\centering\psfig{file=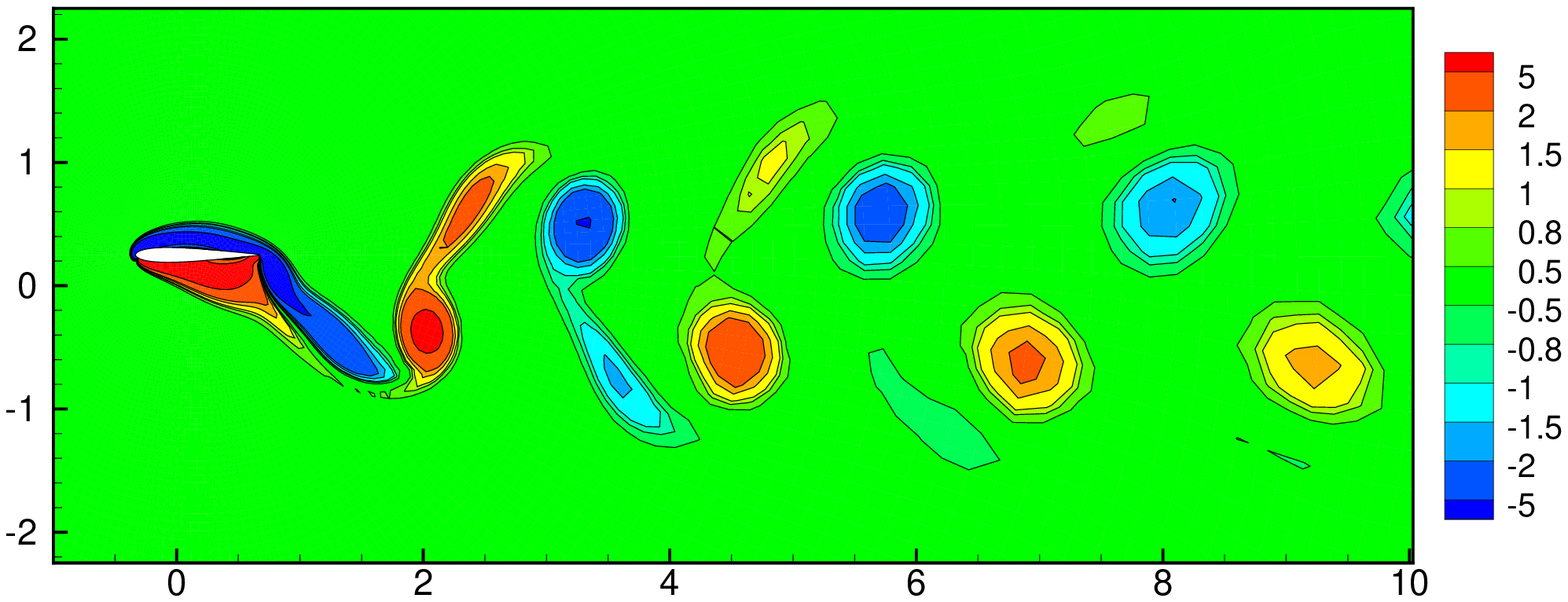,width=0.9\linewidth}
 \\(d) $t_0+T/4$
\end{minipage}
\begin{minipage}[b]{.55\linewidth}\hspace{-1cm}
\centering\psfig{file=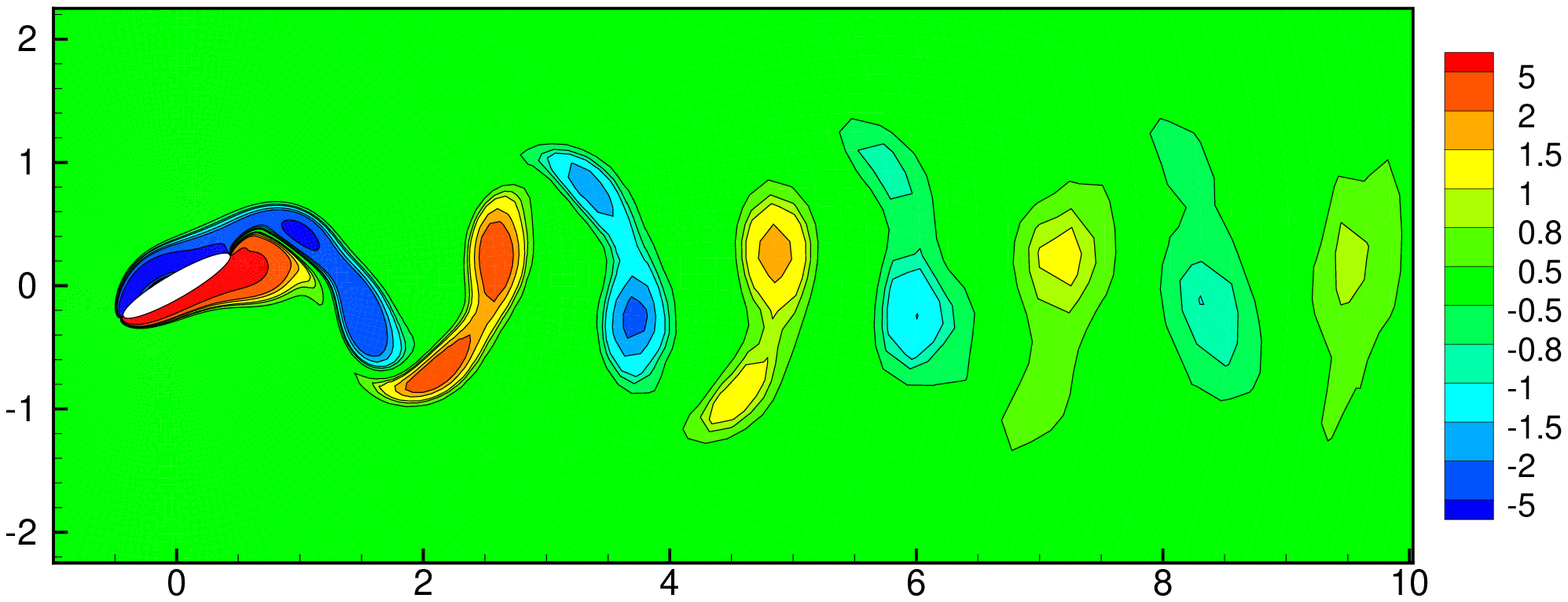,width=0.9\linewidth}
 \\(e) $t_0+2T/4$
\end{minipage}
\begin{minipage}[b]{.55\linewidth}\hspace{-1cm}
\centering\psfig{file=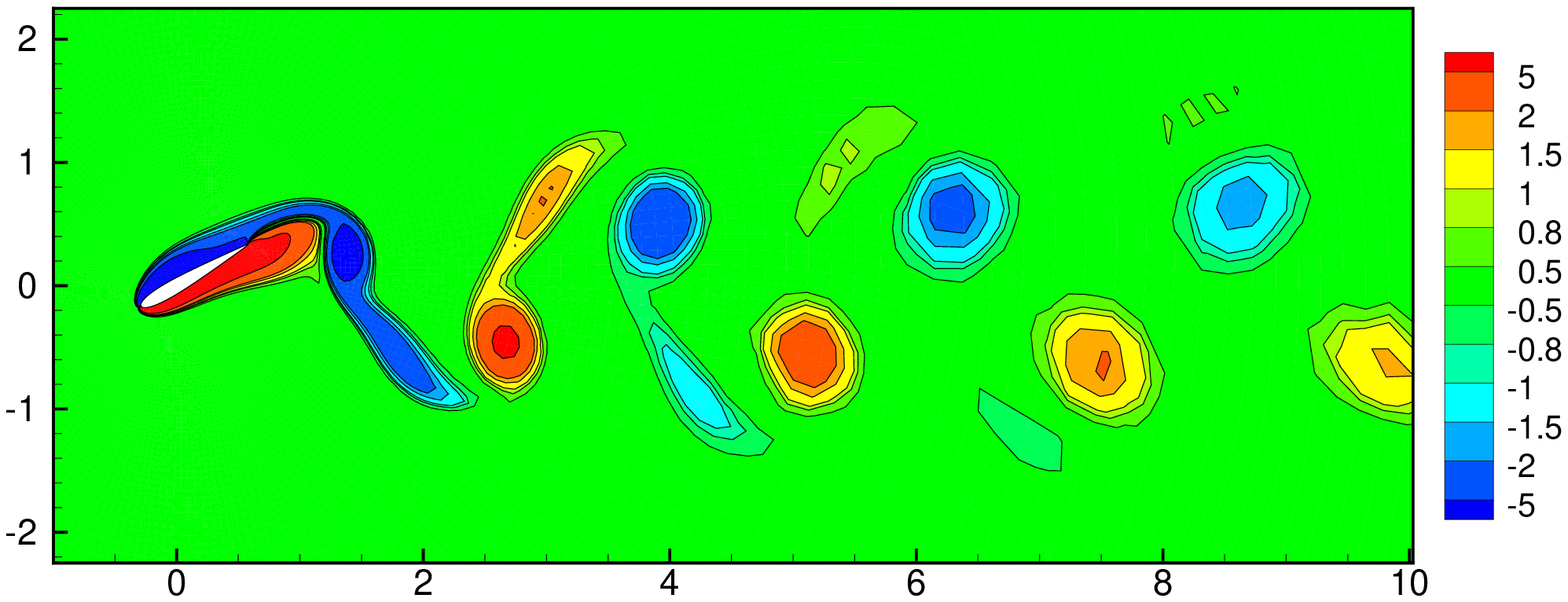,width=0.9\linewidth}
 \\(f) $t_0+2T/4$
\end{minipage}
\begin{minipage}[b]{.55\linewidth}\hspace{-1cm}
\centering\psfig{file=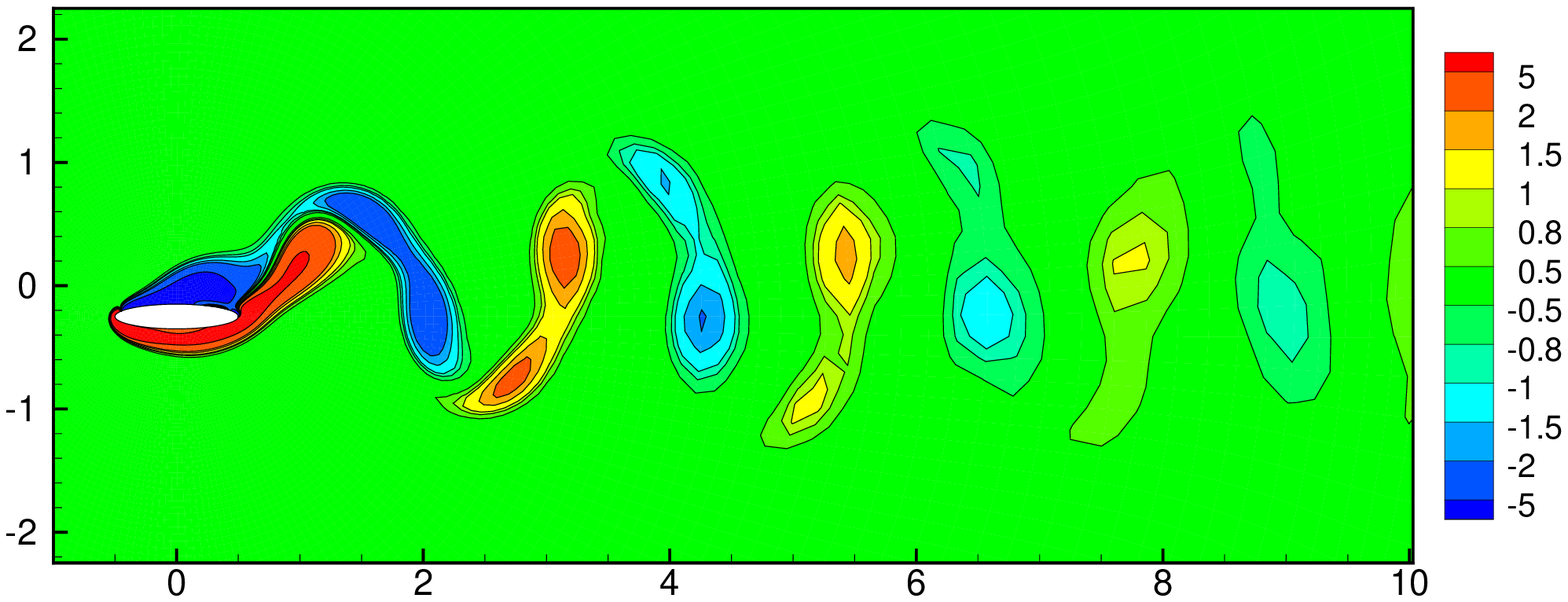,width=0.8\linewidth}
 \\(g) $t_0+3T/4$
\end{minipage}
\begin{minipage}[b]{.55\linewidth}\hspace{-1cm}
\centering\psfig{file=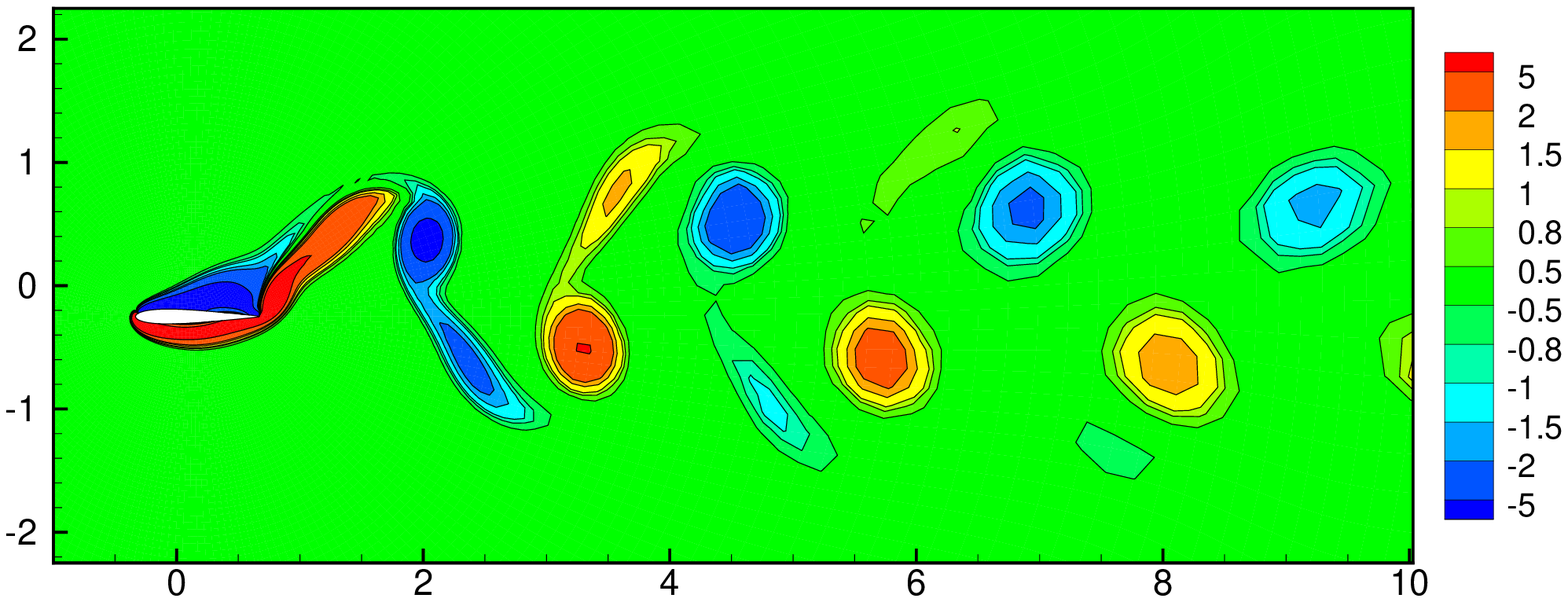,width=0.8\linewidth}
 \\(h) $t_0+3T/4$
\end{minipage}
\begin{center}
\caption{{\sl Problem 7: One cycle of the pitching and heaving airfoil at $Re=200$: elliptic airfoil (left) and NACA0012 (right) at similar stroking phases.} }
    \label{fig:airfoil_vorticity}
\end{center}
\end{figure}

From Fig. \ref{fig:airfoil_dl} we observe that flow field quickly attain periodic nature for both elliptic and NACA0012 airfoil. The rhythmic nature of aerodynamic coefficients is in consonance with pitching and heaving frequency. It is seen that the average drag coefficient of NACA0012 is higher than its elliptic counterpart, whose surface area is twice that of NACA0012. Time averaged $C_D$ values for NACA0012 and elliptic airfoil are 0.783 and 0.435 respectively. On the other hand lift coefficient of NACA0012 has a 15\% higher amplitude at 2.416 compared to 2.096 for the elliptic airfoil.

\begin{figure}[h!]
\begin{minipage}[b]{.55\linewidth}\hspace{-1cm}
\centering\psfig{file=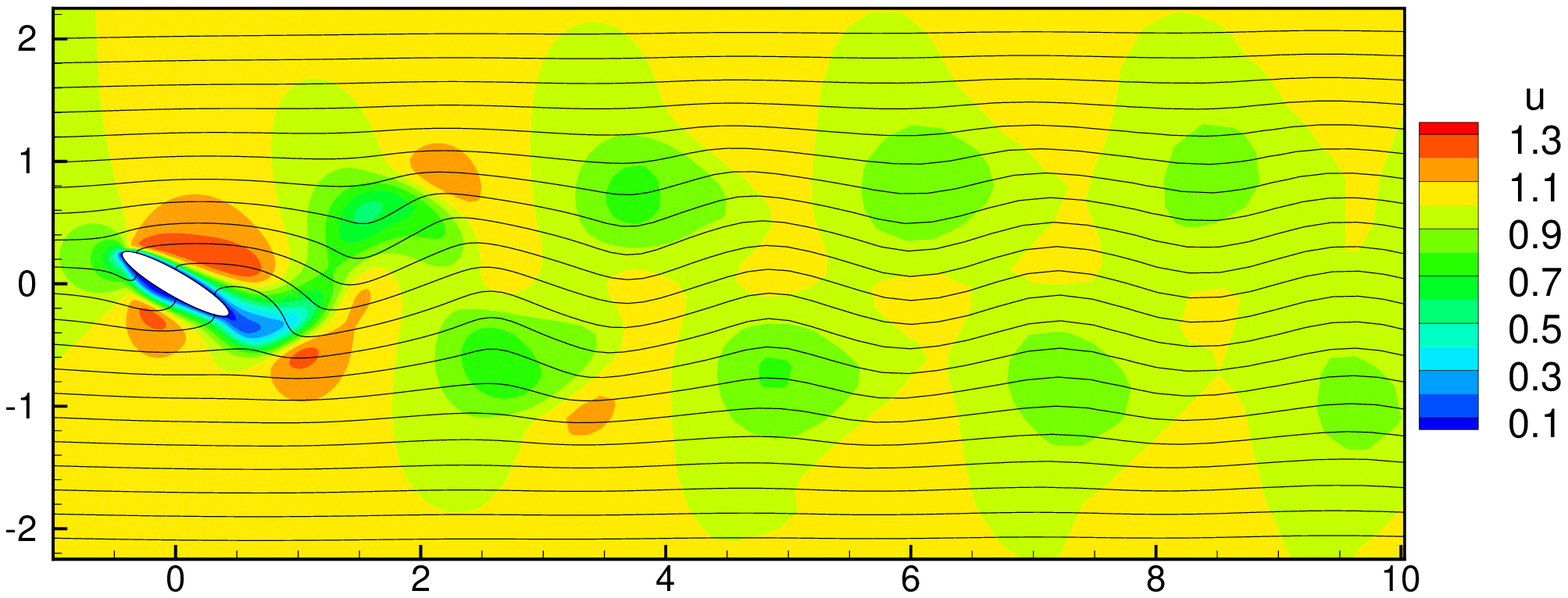,width=0.9\linewidth}
 \\(a) $t_0$
\end{minipage}
\begin{minipage}[b]{.55\linewidth}\hspace{-1cm}
\centering\psfig{file=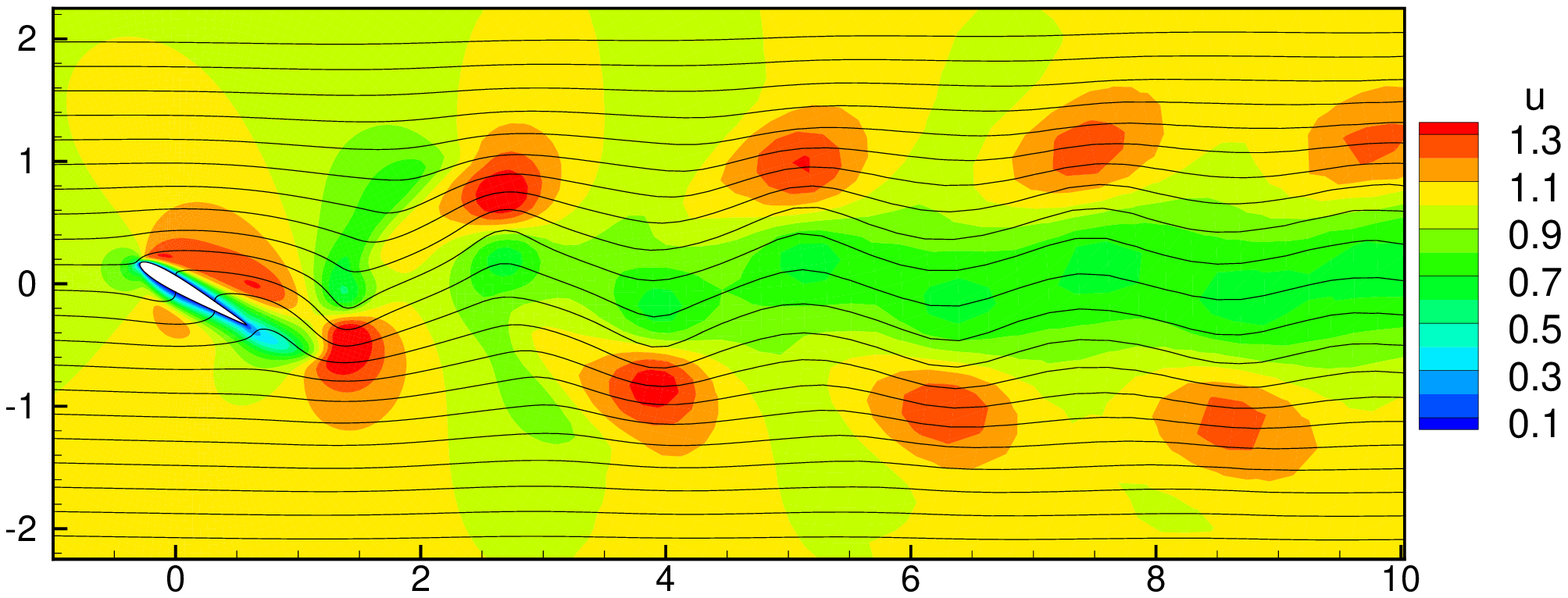,width=0.9\linewidth}
 \\(b) $t_0$
\end{minipage}
\begin{minipage}[b]{.55\linewidth}\hspace{-1cm}
\centering\psfig{file=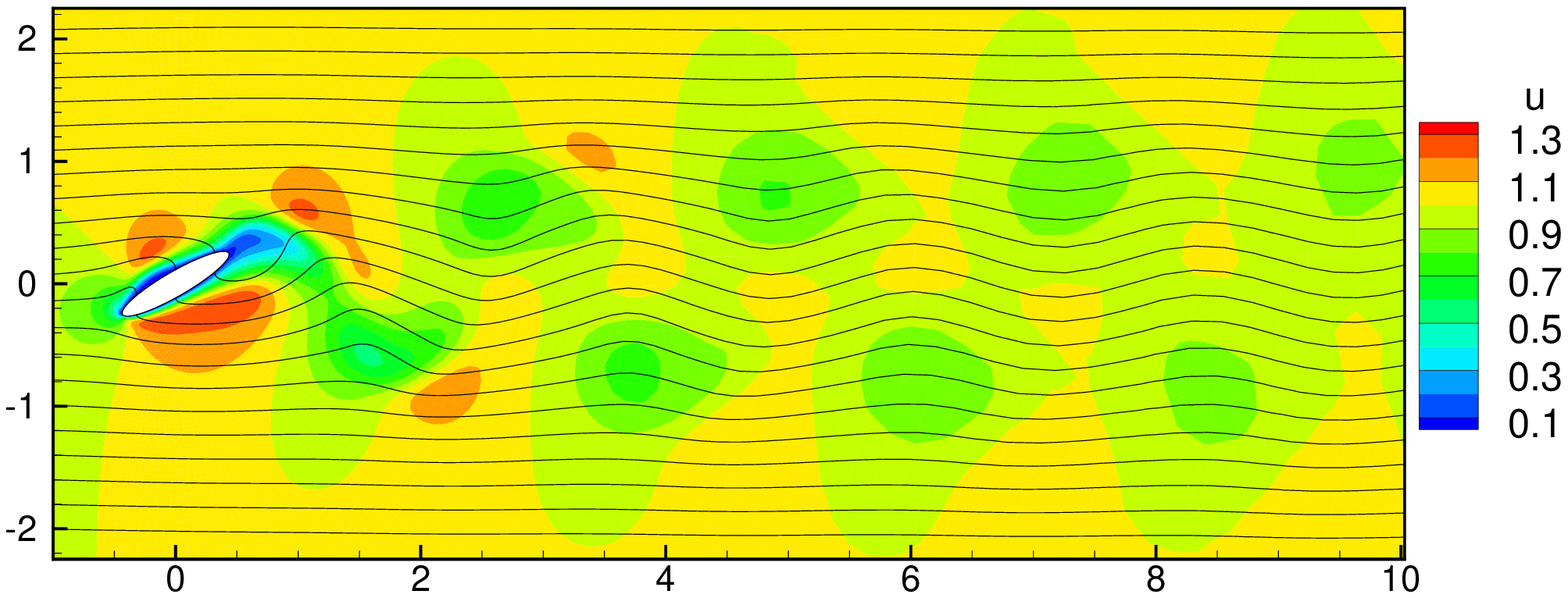,width=0.9\linewidth}
 \\(c) $t_0+T/2$
\end{minipage}
\begin{minipage}[b]{.55\linewidth}\hspace{-1cm}
\centering\psfig{file=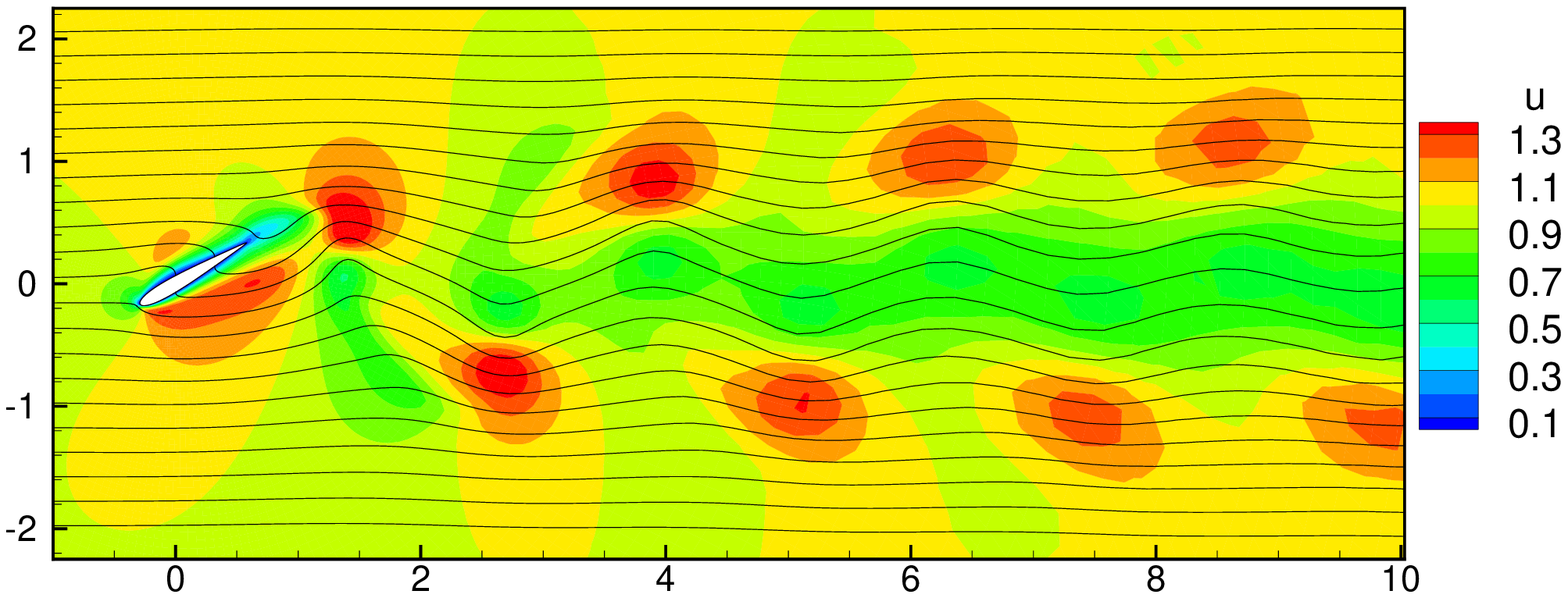,width=0.9\linewidth}
 \\(d) $t_0+T/2$
\end{minipage}
\begin{center}
\caption{{\sl Problem 7: Streamwise velocity $u$ and streamline contours during one cycle of oscillation at $Re=200$: elliptic airfoil (left) and NACA0012 (right) at similar stroking phases.} }
    \label{fig:airfoil_stream}
\end{center}
\end{figure}

We present an excerpt of the vorticity field at four optimal positions of velocity and displacement representing a full cycle of the motion in Fig. \ref{fig:airfoil_vorticity}. Similar stroking phases of both aerofoils are laid out next to each other for sake of comparison. Based on the vorticity field it is obvious that vortices generated at the sharp trailing edge of NACA0012 are stronger \textit{vis-a-vis} elliptic airfoil and relate well to the higher drag value documented in Fig. \ref{fig:airfoil_dl}. Further for NACA0012 shed vortices are found to arrange themselves into two rows based on rotational orientation away from the central line. The downstream vortices behind the elliptic airfoil are elongated, causing the shedding to be weaker and follow $2S$ shedding mode. Fig. \ref{fig:airfoil_stream} shows contour plots of the instantaneous streamlines and streamwise velocity at two different stages of the shedding cycle. These stages correspond to the optimum angle of attack and are alternate mirror-image phases of the flow. Based on the streamwise velocity and the vorticity seen earlier it is obvious that solely in the vicinity of the trailing edge instantaneous phenomena in terms of flow separation and vortex generation take place. For elliptic airfoil, the wake shows reduced unsteadiness and a distinct narrowing \textit{vis-a-vis} NACA0012. Note that vorticity contour comparisons do not provide the same information, as shown in Fig. \ref{fig:airfoil_stream} with the help of streamline contours. 

\section{Conclusions}
A fourth order compact numerical scheme for solving unsteady incompressible flow on time-dependent curvilinear coordinates is presented. The technique developed approximates the convection-diffusion equation on irregular moving geometry and is amenable to both primitive variable and streamfunction-vorticity formulations of the Navier-Stokes (N-S) equations. The method uses compact computational stencil and hence discretization involves neighboring grid points only. In computations of the unsteady convection-diffusion equation and the incompressible N-S equation, the spatial and temporal accuracy of the method is demonstrated. Estimation of grid metrics in this context is rather straightforward as the current study reveals that the satisfaction of GCL is not imperative. In this work, the appropriate implementation of the scheme was achieved in tandem with an efficient grid deformation procedure. A successful simulation of the FSI test cases is a fitting proof of this. Most likely, this is the first attempt at a compact, higher-order computation of non-conservative N-S equations in a continually deforming region. The numerical method described in this study might be useful for the theoretical and numerical investigation of fluid dynamics in the presence of deformable surfaces.

\section*{Acknowledgement}
The first author is thankful to Dr. Guillaume De Nayer, Professur f{\"u}r Str{\"o}mungsmechanik, Helmut-Schmidt-Universit{\"a}t Hamburg, Germany for some encouraging comments. He also acknowledges the use of facilities at the High-Performance Computing Centre, Tezpur University sponsored by DeitY, India in collaboration with C-DAC, India.

%\section*{References}


\begin{thebibliography}{00}

\bibitem{wan_fid_abg_13}
Z.J. Wang, K. Fidkowski, R. Abgrall, F. Bassi, D. Caraeni, A. Cary, H. Deconinck, R. Hartmann, K. Hillewaert, H.T. Huynh, N. Kroll, G. May, P.-O. Persson, B. van Leer and M. Visbal, High-order CFD methods: current status and perspective, International Journal for Numerical Methods in Fluids 72 (2013) 811-845.

\bibitem{vis_gai_02}
M.R. Visbal and D.V. Gaitonde, On the use of higher-order finite-difference schemes on curvilinear and deforming meshes, Journal of Computational Physics 181 (2002) 155-185.

\bibitem{ste_78}
J.L. Steger, Implicit finite-difference simulation of flow about arbitrary two-dimensional geometries, AIAA Journal 16 (1978) 679-686.

\bibitem{gha_far_12}
M. Ghadimi and M. Farshchi, Fourth order compact finite volume scheme on nonuniform grids with multi-blocking, Computers and Fluids 56 (2012) 1-16.

\bibitem{lia_miy_zha_14}
C. Liang, K. Miyaji and B. Zhang, An efficient correction procedure via reconstruction for simulation of viscous flow on moving and deforming domains, Journal of Computational Physics 256 (2014) 55-68.

\bibitem{li_hen_ban_16}
L. Li, W.D. Henshaw, J.W. Banks, D.W. Schwendeman and A. Main, A stable partitioned FSI algorithm for incompressible flow and deforming beams, Journal of Computational Physics 312 (2016) 272-306.

\bibitem{cox_lia_ple_16} 
C. Cox, C. Liang and M.W. Plesniak, A high-order solver for unsteady incompressible Navier-Stokes equations using the flux reconstruction method on unstructured grids with implicit dual time stepping, Journal of Computational Physics 314 (2016) 414-435.

\bibitem{hir_nic_81}
C.W. Hirt and B.D. Nichols, Volume of fluid (VOF) method for the dynamics of free boundaries, Journal of Computational Physics 39 (1981) 201-255.

\bibitem{li_hes_zie_05}
J. Li, M. Hesse, J. Ziegler and A.W. Woods, An arbitrary Lagrangian Eulerian method for moving-boundary problems and its application to jumping over water, Journal of Computational Physics 208 (2005) 289-314.

\bibitem{sjo_yee_vin_14}
B. Sj\"{o}green, H.C. Yee and M. Vinokur, On high order finite-difference metric discretizations
satisfying GCL on moving and deforming grids, Journal of Computational Physics 265 (2014) 211-220.

\bibitem{abe_hag_non_16}
Y. Abe, T. Haga, T. Nonomura and K. Fujii, Conservative high-order flux-reconstruction schemes on moving and deforming grids, Computers and Fluids 139 (2016) 2-16.

\bibitem{per_bon_per_09}
P.-O. Persson, J. Bonet and J. Peraire, Discontinuous Galerkin solution of the Navier-Stokes equations on deformable domains, Computer Methods in Applied Mechanical Engineering 198 (2009) 1585-1595.

\bibitem{wan_prz_94}
Z.J. Wang and A.J. Przekwas, Unsteady flow computation using moving grid with mesh enrichment, AIAA Paper: AIAA-
94-0285, (1994).

\bibitem{boe_sch_bij_07}
A. de Boer, M.S. van der Schoot and H. Bijl, Mesh deformation based on radial basis function interpolation, Computers
and Structures 85 (2007) 784-795.

\bibitem{ren_all_08}
T.C.S. Rendall and C.B. Allen, Unified fluid--structure interpolation and mesh motion using radial basis functions.
International Journal for Numerical Methods in Engineering 74 (2008) 1519-1559.

\bibitem{wit_bij_09}
J.A.S. Witteveen and H. Bijl, Explicit mesh deformation using inverse distance weighting interpolation. 19th AIAA
Computational Fluid Dynamics Conference, AIAA, San Antonio, Texas, AIAA Paper, Vol. 3996, 2009.

\bibitem{luk_col_bla_12}
E. Luke, E. Collins and E. Blades, A fast mesh deformation method using explicit interpolation, Journal of Computational
Physics 231 (2012) 586-601.

\bibitem{sen_nay_bre_17}
S. Sen, G. De Nayer and M. Breuer, A fast and robust hybrid method for block-structured mesh deformation with emphasis on FSI-LES applications, International Journal for Numerical Methods in Engineering 111 (2017) 273-300.

\bibitem{apo_nay_ble_19}
A. Apostolatos, G. De Nayer, G., K-U Bletzinger, M. Breuer and R. W{\"u}chner, Systematic evaluation of the interface description for fluid--structure interaction simulations using the isogeometric mortar-based mapping, Journal of Fluids and Structures 86 (2019) 368-399.

\bibitem{nay_bre_woo_20} 
G. De Nayer, M. Breuer and J.N. Wood, Numerical investigations on the dynamic behavior of a 2-DOF airfoil in the transitional $Re$ number regime based on fully coupled simulations relying on
an eddy-resolving technique, International Journal of Heat and Fluid Flow 85 (2020) 108631.

\bibitem{sen_13}
S. Sen, A new family of (5,5)CC-4OC schemes applicable for unsteady Navier-Stokes equations, Journal of Computational Physics 251 (2013) 251-271.

\bibitem{sen_16}
S. Sen, Fourth order compact schemes for variable coefficient parabolic problems with mixed derivatives, Computers and Fluids 134-135 (2016) 81-89.

\bibitem{sen_she_17}
S. Sen and T. W. H. Sheu, On the development of a nonprimitive Navier-Stokes formulation subject to rigorous implementation of a new vorticity integral condition, Journal of Scientific Computing 72 (2017) 252-290.

\bibitem{che_xie_16}
Y. Chen and X. Xie, Vorticity vector-potential method for 3D viscous incompressible flows in time-dependent curvilinear coordinates, Journal of Computational Physics 312 (2016) 50-81.

\bibitem{tay_23}
G.I. Taylor, On the decay of vortices in a viscous  fluid, Philosophical Magazine 46 (1923) 671-674.

\bibitem{wes_seg_kas_99}
P. Wesseling, A. Segal and C. G. M. Kassels, Computing flows on general three-dimensional nonsmooth staggered grids, Journal of Computational Physics 149 (1999) 333-362.

\bibitem{ge_sot_07}
L. Ge and F. Sotiropoulos, A numerical method for solving the 3D unsteady incompressible Navier–Stokes equations in curvilinear domains with complex immersed boundaries, Journal of Computational Physics 225 (2007) 1782-1809.

\bibitem{che_xie_16}
Y. Chen and X. Xie, Vorticity vector-potential method for 3D viscous incompressible flows in time-dependent curvilinear coordinates, Journal of Computational Physics 312 (2016) 50-81.

\bibitem{ghi_ghi_shi_82}
U. Ghia, K. Ghia and C. Shin, High-Re solutions for incompressible flow using the Navier-Stokes equations and a multigrid method, Journal of Computational Physics 48 (1982) 387-411.

\bibitem{bru_saa_06}
C. H. Bruneau and M. Saad, The 2D lid-driven cavity problem revisited, Computers and Fluids 35 (2006) 326-348.

\bibitem{wil_ros_88}
C. H. K. Williamson and A. Roshko, Vortex formation in the wake of an oscillating cylinder Journal of Fluids and Structures 2 (1988) 355-381.

\bibitem{knu_03}
P. M. Knupp, Algebraic mesh quality metrics for unstructured initial meshes, Finite Elements in Analysis and Design 39 (2003) 217-241.

\bibitem{kal_sen_12}
J. C. Kalita and S. Sen, Triggering asymmetry for flow past circular cylinder at low Reynolds numbers, Computers and Fluids 59 (2012) 44-60.

\bibitem{bri_71}
W. R. Briley, A numerical study of laminar separation bubbles using Navier-Stokes equations, Journal of Fluid Mechanics 47 (1971) 713-736.

\bibitem{tia_lia_yu_11}
Z. F. Tian, X. Liang and P. X. Yu, A higher order compact finite difference algorithm for solving the incompressible Navier-Stokes equations, International Journal for Numerical Methods in Engineering 88 (2011) 511-32.

\bibitem{leo_ste_tho_06}
J. S. Leontini, B. E. Stewart, M. C. Thompson and K. Hourigan, Wake state and energy transitions of an oscillating cylinder at low Reynolds number, Physics of Fluids 18 (2006) 067101 1-9.

\bibitem{wu_ma_zho_07}
J.Z. Wu, H.Y. Ma and M.D. Zhou, Vorticity and Vortex Dynamics, Springer, Berlin, 2007.

\bibitem{for_80}
B. Fornberg, A numerical study of steady viscous flow past a circular cylinder, Journal of
Fluid Mechanics 98 (1980) 819-55.

\bibitem{sen_kal_gup_13}
S. Sen, J.C. Kalita, M.M. Gupta, A robust implicit compact scheme for two-dimensional unsteady
flows with a biharmonic stream function formulation, Computers and Fluids 84 (2013) 141-163.

\bibitem{eld_07}
J.D. Eldredge, Numerical simulation of the fluid dynamics of 2D rigid body motion with the vortex particle method, Journal of Computational Physics 221 (2007) 626-648.

\bibitem{den_you_03}
S.C.R. Dennis and P.J.S. Young, Steady flow past an elliptic cylinder inclined to the stream, Journal of Engineering Mathematics 47 (2003) 101-120.

\bibitem{med_sto_car_11}
W. Medjroubi, B. Stoevesandt, B. Carmo, J. Peinke, High-order numerical simulations of the flow around a heaving airfoil, Computers and Fluids 51 (2011) 68-84.
\end{thebibliography}
\end{document}